\documentclass[12pt]{article}
\usepackage{epsfig}
\usepackage{amssymb,amsmath,amsthm,url}
\usepackage{multirow}
\usepackage{booktabs}
\usepackage[numbers]{natbib}
\usepackage{setspace}

\usepackage{dsfont}

\usepackage{times,amsfonts,eucal}
\usepackage{algorithm}
\usepackage{algorithmic}
\usepackage{graphicx,ifthen}
\usepackage{wrapfig}
\usepackage{latexsym}
\usepackage{color}
\usepackage{mathrsfs}
\usepackage{bm}
\usepackage{bbm}
\usepackage{srctex}

\usepackage{bm}
\newtheorem{theorem}{Theorem}[section]
\newtheorem{lemma}{Lemma}[section]

\newtheorem{cor}{Corollary}[section]

\newtheorem{rem}{Remark}[section]

\newtheorem{assumption}{Assumption}[section]

\numberwithin{equation}{section}
\numberwithin{theorem}{section}
\numberwithin{lemma}{section}
\numberwithin{pro}{section}
\numberwithin{cor}{section}
\numberwithin{definition}{section}
\numberwithin{cons}{section}
\numberwithin{rem}{section}
\numberwithin{exa}{section}
\numberwithin{table}{section}
\numberwithin{figure}{section}
\numberwithin{algo}{section}

\newcommand{\mR}{\mathbb{R}}

\newcommand{\mE}{\mathbb{E}}
\newcommand{\mP}{\mathbb{P}}

\newcommand{\D}{\mathbf{D}}

\newcommand{\Z}{\mathbb{Z}}

\def\beq{\begin{equation}}
\def\eeq{\end{equation}}
\def\bals{\begin{align*}}
\def\eals{\end{align*}}

\def\bal{\begin{align}}
\def\eal{\end{align}}

\allowdisplaybreaks

\begin{document}

\title{Functional data analysis in the Banach space of continuous functions\footnote{This research was partially supported by NSF grants DMS 1305858 and DMS 1407530, by the Collaborative Research Center `Statistical modeling of nonlinear dynamic processes' ({\it Sonderforschungsbereich 823, Teilprojekt A1, C1}) and by the Research Training Group `High-dimensional Phenomena in Probability - Fluctuations and Discontinuity' ({\it RTG 2131}) of the German Research Foundation. Part of the research was done while A.\ Aue was visiting Ruhr-Universit\"at Bochum as a Simons Visiting Professor of the Mathematical Research Institute Oberwolfach.
}}

\author{
 Holger Dette\footnote{Fakult\"at f\"ur Mathematik, Ruhr-Universit\"at Bochum, Universit\"atsstra{\ss}e 150, 44780 Bochum, Germany, emails: \tt{[holger.dette, kevin.kokot]@ruhr-uni-bochum.de}}
\and Kevin Kokot$^\dagger$ \and Alexander Aue\footnote{Department of Statistics, University of California, Davis, CA 95616, USA, email: \tt{aaue@ ucdavis.edu}}
}

\date{\today}
\maketitle
\bibliographystyle{plain}

\begin{abstract}
Functional data analysis is typically conducted within the $L^2$-Hilbert space framework. There is by now a fully developed statistical toolbox allowing for the principled application of the functional data machinery to real-world problems, often based on dimension reduction techniques such as functional principal component analysis. At the same time, there have recently been a number of publications that sidestep dimension reduction steps and focus on a fully functional $L^2$-methodology. This paper goes one step further and develops  data analysis methodology for functional time series in the  space of all continuous functions. The work is motivated by the fact that objects with rather different shapes may still have a small $L^2$-distance and are therefore identified as similar when using an $L^2$-metric. However,  in applications it is often desirable to use metrics reflecting the visualization of the curves in the statistical analysis. The methodological contributions are focused on developing  two-sample and change-point tests as well as confidence bands, as these procedures appear do be conducive to the proposed setting. Particular interest is put on relevant differences; that is, on not trying to test for exact equality, but rather for pre-specified deviations under the null hypothesis.

The procedures are justified through large-sample theory. To ensure practicability, non-standard bootstrap procedures are developed and investigated addressing
 particular features that arise 
  in the problem of testing relevant hypotheses. The finite sample properties are explored through a simulation study and an application to annual temperature profiles.  \medskip 

\end{abstract}

\renewcommand{\baselinestretch}{1.35}

%%%%%%%%%%%%%%%%%%%%%%%%
\section{Introduction}
\label{sec:intro}
%%%%%%%%%%%%%%%%%%%%%%%%

This paper proposes new methodology for the analysis of functional data, in particular for the two-sample and change-point settings. The basic set-up considers a sequence of Banach space-valued time series satisfying mixing conditions. The proposed methodology therefore advances functional data analysis beyond the predominant Hilbert space-based methodology. For the latter case, there exists by now a fully fledged theory. The interested reader is referred to the various monographs Ferraty and Vieu \cite{FerratyVieu2010}, Horv\'ath and Kokoszka \cite{HorvathKokoskza2012}, and Ramsay and Silverman \cite{RamsaySilverman2005} for up-to-date accounts. Most of the available statistical procedures discussed in these monographs are based on dimension reduction techniques such as functional principal component analysis. However, the integral role of smoothness has been discussed at length in Ramsay and Silverman \cite{RamsaySilverman2005} and virtually all functions fit in practice are at least continuous. In such cases dimension reduction techniques can incur a loss of information and fully functional methods can prove advantageous. More recently, Aue et al.\ \cite{AueRiceSonmez2015}, Bucchia and Wendler \cite{BucchiaWendler2015} and  Horv\'ath et al.\ \cite{HorvathKokoszkaRice2014} discussed fully functional methodology in a Hilbert space framework.

Since all functions utilized for practical purposes are at least continuous, and often smoother than that, it might be more natural to develop methodology for functional data in the  space   of continuous functions. This is the approach pursued in the present paper. While it might thus be reasonable to build statistical analysis adopting this  point of view, there are certain difficulties associated with it. Giving up on the theoretically convenient Hilbert space setting means that substantially more effort has to be put into the derivation of theoretical results, especially if one is interested in the incorporation of dependent functional observations. Section~\ref{sec:banach_methods} of the main part of this paper gives an introduction to Banach space methodology and states some basic results, in particular an invariance principle for a sequential process in the space of continuous functions.

The theoretical contributions will be utilized for the development of relevant two-sample and change-point tests in Sections \ref{sec:two-sample} and \ref{sec:change-point}, respectively. Here the usefulness of the proposed approach becomes more apparent as differences between two smooth curves are hard to detect in practice. Additionally, small discrepancies might perhaps not even be of importance in many applied situations. Therefore the ``relevant'' setting is adopted that is not trying to test for exact equality under the null hypothesis, but allows for pre-specified deviations from an assumed null function. For example, if $C(T)$, the space of continuous functions on the compact interval $T$, is equipped with the sup-norm $\|f\|=\sup_{t\in T} | f(t) | $, and $\mu_1$ and $\mu_2$ are the mean functions corresponding to two samples, interest is in hypotheses of the form
\begin{equation} \label{relevant}
H_0\colon \| \mu_1 - \mu_2 \| \leq \Delta
\qquad\mbox{and}\qquad
H_1\colon \| \mu_1 - \mu_2 \| > \Delta,
\end{equation}
where $\Delta \geq 0$ denotes a pre-specified constant. The classical case of testing perfect equality, obtained by the choice $\Delta =0$, is therefore a special case of \eqref{relevant}. However, in applications it might be reasonable to think about this choice carefully and to define precisely the size of change which one is really interested in. In particular, testing relevant hypotheses avoids the consistency problem as mentioned in Berkson \cite{berkson1938}, that is: any consistent test will detect any arbitrary small change in the mean functions if the sample size is sufficiently large. One may also view this perspective as a particular form of a bias-variance trade-off. The problem of testing for a relevant difference between two (one-dimensional) means and other (finite-dimensional) parameters has been discussed by numerous authors in biostatistics (see Wellek \cite{wellek2010} for a recent review), but to the best of our knowledge these testing problems have not been considered in the context of functional data. It turns out that from a mathematical point of view the problem of testing relevant (i.e., $\Delta >0$) hypotheses is substantially more difficult than the classical problem (i.e., $\Delta =0$). In particular, it is not possible to work with stationarity under the null hypothesis, making the derivation of a limit distribution of a corresponding test statistic or the construction of a bootstrap procedure substantially more difficult.

 Section \ref{sec:two-sample}  develops corresponding two-sample tests for the Banach space $C(T)$. Section \ref{sec:change-point} extends these results to the change-point setting (see Aue and Horv\'ath \cite{aueHorvath2013} for a recent review of change-point methodology for time series). Here, one has to deal with the additional complexity of locating the unknown time of change. Several new results for change-point analysis of functional data in $C(T)$ are put forward. A specific challenge here is the fact that the asymptotic null distribution of test statistics for hypotheses of the type \eqref{relevant} depends on the set of extremal points of the unknown difference $\mu_1 - \mu_2 $, and is therefore not distribution free. Most notable for both the two-sample and the change-point problem is the construction of non-standard bootstrap tests for relevant hypotheses to solve this problem. The bootstrap is theoretically validated and then used to determine cut-off values for the proposed procedures.

Another area of application that lends itself naturally to Banach space methodology is that of constructing confidence bands for the mean function of a collection of potentially temporally dependent, continuous functions. There has been recent work by Choi and Reimherr \cite{choiReimherr2016} on this topic in a Hilbert space framework for functional parameters of independent functions based on geometric considerations. Here, results for confidence bands for the mean difference in a two-sample framework are added in Section \ref{sec:conf_bands}. Natural modifications allow for the inclusion of the one-sample case. One of the main differences between the two approaches is that the proposed bands hold pointwise, while those constructed from Hilbert space theory are valid only in an $L^2$-sense. This property is appealing for practitioners, because two mean curves can have a rather different shape, yet the $L^2$-norm of their difference might be very small.

The finite-sample properties of the relevant two-sample and change-point tests and, in particular, the performance of the bootstrap procedures are evaluated with the help of a Monte Carlo simulation study in Section \ref{sec5}. A number of scenarios are investigated, with the outcomes showing that the proposed methodology performs reasonably well. Furthermore, an application to a prototypical data example is given, namely two-sample and cange-point
tests for annual temperature profiles recorded at measuring stations in Australia.

The outline of the rest of this paper is as follows. Section \ref{sec:banach_methods} introduces the basic notions of the proposed Banach space methodology and gives some preliminary results. Section \ref{sec:two-sample} discusses the two-sample problem and Section \ref{sec:change-point} is concerned with change-point analysis. Empirical aspects are highlighted in Section \ref{sec5}. Proofs of the main results can be found in an online supplement to this paper.

%%%%%%%%%%%%%%%%%%%%%%%%
 \section{$C(T)$-valued random variables}
 \label{sec:banach_methods}
%%%%%%%%%%%%%%%%%%%%%%%%

In this section some basic facts are provided about central limit theorems and invariance principles for $C(T)$-valued random variables, where $C(T)$ is the set of continuous functions from $T$ into the real line $\mathbb{R}$. In what follows, unless otherwise mentioned, $C(T)$ will be equipped with the sup norm $\|\cdot\|$, defined by $\|f\|=\sup_{t\in T}|f(t)|$, thus making $(C(T),\|\cdot\|)$ a Banach space. The natural Borel $\sigma$-field $\mathcal{B}(T)$ over $C(T)$ is then generated by the open sets relative to the sup norm $\|\cdot\|$. Measurability of random variables on $(\Omega,\mathcal{A},P)$ taking values in $C(T)$ is understood to be with respect to $\mathcal{B}(T)$. The underlying probability space $(\Omega,\mathcal{A},\mathbb{P})$ is assumed complete. It is further assumed that there is a metric $\rho$ on $T$ such that $(T,\rho)$ is totally bounded. The fact that $T$ is metrizable implies that $C(T)$ is separable and measurability issues are avoided (see Theorem 7.7 in Janson and Kaijser \cite{jankai2015}). Moreover, any random variable $X$ in $C(T)$ is tight (see Theorem 1.3 in Billingsley \cite{billingsley1968}).

Let $X$ be a random variable on $(\Omega,\mathcal{A},P)$ taking values in $C(T)$. There are different ways to formally introduce expectations and higher-order moments of Banach space-valued random variables (see Janson and Kaijser \citep{jankai2015}). The expectation $\mathbb{E}[X]$ of a random variable $X$ in $C(T)$ exists as an element of $C(T)$ whenever $\mathbb{E}[\| X\|] < \infty$. The  $k$th moment exists whenever $\mathbb{E}[\| X \|^k] = \mathbb{E}[ \sup_{t \in T} |X(t)|^k] < \infty $. As pointed out in Chapter 11 of Janson and Kaijser \cite{jankai2015}, $k$th order moments may be computed through pointwise evaluation as $\mathbb{E}   [X (t_1)  \cdots   X(t_k)]$. The case $k=2$ is important as it allows for the computation of covariance kernels in a pointwise fashion.

A sequence of random variables $(X_n\colon n \in \mathbb{N})$ converges in distribution or weakly to a random  variable $X$ in $C(T)$, whenever it is asymptotically tight and its finite-dimensional distributions converge weakly to the finite-dimensional distributions of $X$, that is,
\[
(X_n(t_1), \dots, X_n(t_k)) \Rightarrow (X(t_1), \dots, X(t_k))
\]
for any $t_1,\dots,t_k \in T$ and any $k \in \mathbb{N}$, where the symbol ``$\Rightarrow$'' indicates convergence in distribution in $\mathbb{R}^k$.

A centered random variable $X$ in $C(T)$ is said to be Gaussian if its finite-dimensional distributions are multivariate normal, that is, for any $t_1,\dots,t_k$, $(X(t_1),\dots,X(t_k)) \sim \mathcal{N}_k(0,\Sigma)$, where the $(i,j)$th entry of the covariance matrix $\Sigma$ is given by $\mathbb{E}[X(t_i)X(t_j)]$, $i,j=1,\ldots,k$. The distribution of $X$ is hence completely characterized by its covariance function $k(t,t^\prime) = \mathbb{E}[X(t)X(t^\prime)]$; see Chapter 2 of Billingsley \cite{billingsley1968}.

In general Banach spaces, deriving conditions under which the central limit theorem (CLT) holds is a difficult task, significantly more complex than the counterpart for real-valued random variables. In Banach spaces, finiteness of second moments of the underlying random variables does not provide a necessary and sufficient condition. Elaborate theory has been developed to resolve the issue, resulting in notions of Banach spaces of type 2 and cotype 2 (see the book Ledoux and Talagrand \cite{leTa1991} for an overview). However, the Banach space of continuous functions on a compact interval does not possess the requisite type and cotype properties and further assumptions are needed in order to obtain the CLT, especially to incorporate time series of continuous functions into the framework. To model the dependence of the observations, the notion of $\varphi$-mixing triangular arrays $(X_{n,j}\colon n\in\mathbb{N},~ j=1,\dots,n)$ of $C(T)$-valued random variables is introduced; see Bradley \cite{bradley2005} and Samur \cite{samur1987}. First, for any two $\sigma$-fields $\mathcal{F}$ and $\mathcal{G}$, define
\[
\phi(\mathcal{F},\mathcal{G})
= \sup \big\{ |\mathbb{P}(G|F) - \mathbb{P}(G)| \colon F\in \mathcal{F}, ~G\in\mathcal{G}, ~\mathbb{P}(F)>0 \big\},
\]
where $\mathbb{P}(G|F)$ denotes the conditional probability of $G$ given $F$. Next, denote by $\mathcal{F}^n_{k,k^\prime}$ the $\sigma$-field generated by $(X_{n,j}\colon k\leq j \leq k^\prime)$. Then, define the $\varphi$\textit{-mixing coefficient} as
\[
\varphi (k)
= \sup_{n \in \mathbb{N}, n>k} \max_{k^\prime= 1,\dots,n-k} \phi (\mathcal{F}_{1,k^\prime}^n ,\mathcal{F}_{k^\prime+k,n}^n)
\]
and call the triangular array $(X_{n,j}\colon n\in\mathbb{N}, j=1,\dots,n)$ $\varphi$\textit{-mixing} whenever $\lim_{k\to\infty} \varphi(k) = 0$. The $\varphi$-mixing property is defined in a similar fashion for a sequence of random variables.

In order to obtain a CLT as well as an invariance principle for triangular arrays of $\varphi$-mixing random elements in $C(T)$, the following conditions are imposed.

\begin{assumption}
\label{as:ts}
Throughout  this paper the following conditions  are assumed to hold: \vspace{-.2cm}
\begin{enumerate}\itemsep-.4ex
\item[(A1)]
There is a constant $K$ such that, for all $n\in\mathbb{N}$ and $j=1,\dots,n$, $\mathbb{E}[\|X_{n,j}\|^{2+\nu}] \leq K$ for some $\nu>0$.
\item[(A2)]
Let $\mathbb{E}[X_{n,j}] = \mu^{(j)}$ for any $n\in\mathbb{N}$ and $j=1,\dots,n$.
The distributions of the observations in each row only differ in their means, that is, the centered array $( X_{n,j}-\mu^{(j)} \colon n\in\mathbb{N},~ j=1,\dots,n)$ is rowwise stationary. Additionally, the covariance structure is the same in each row, that is  $$\mathrm{Cov}(X_{n,j}(t),X_{n,j^\prime}(t^\prime)) = \gamma(j-j^\prime,t,t^\prime)$$ for all $n\in\mathbb{N}$ and $j,j^\prime=1,\dots,n$. Note that $\gamma(-j,t,t^\prime) = \gamma(j,t^\prime,t)$.
\item[(A3)]
$(X_{n,j}\colon n\in\mathbb{N},~ j=1,\dots,n)$ is uniformly Lipschitz, that is, there is a real-valued random variable $M$ with $\mathbb{E}[M^2] < \infty$ such that, for any $n\in\mathbb{N}$ and $j=1,\dots,n$, the inequality
\[
|X_{n,j}(t)-X_{n,j}(t^\prime)|\leq M \rho(t,t^\prime) 
\]
holds almost surely for all $t,t^\prime\in T$.
\item[(A4)]
$(X_{n,j}\colon n\in\mathbb{N},~ j=1,\dots,n)$ is $\varphi$-mixing with exponentially decreasing mixing coefficients, that is, there is a constant $a\in [0,1)$ such that $\varphi (k)\leq a^k$ for any $k\in\mathbb{N}$.
\item[(A5)]
For any sequence \((r_n)_{n\in\mathbb{N}}\subset \mathbb{N}\) such that \(r_n\leq n\),
\({r_n}/{n} \to 0\) as \(n\to\infty\), it follows that
\begin{align*}
\frac{1}{\sqrt{n}} \sum_{j=1}^{r_n} (X_{n,j}-\mu^{(j)}) = o_\mP (1) .
\end{align*}
\end{enumerate}
\end{assumption}

Note that these assumptions can be formulated for sequences of random variables $(X_n\colon n\in\mathbb{N})$ in $C(T)$ in a similar way. Condition  (A5)  is satisfied if the distribution of the sums $\sum_{j=1}^k (X_{n,j} - \mu^{(j)})$ is symmetric for any $k=1,\dots,n$ and $n\in\mathbb{N}$ (see the remark after Proposition 3.1 in Samur \cite{samur1984}). Assumptions (A1)--(A4) imply the following CLT which is proved in Section \ref{proofssec2} of the online supplement. {Throughout this paper the symbol $\rightsquigarrow $ denotes weak convergence  in $(C([0,1]))^k $ for  some $k \in \mathbb{N}$.}

\begin{theorem} \label{mixingCLT}
Let $(X_{n,j}\colon n\in\mathbb{N},~ j=1,\dots,n)$ denote a triangular array of random variables in $C(T)$ with expectations $E[X_{n,j}] = \mu^{(j)}$ such that conditions (A1) -- (A4) of Assumption \ref{as:ts} are satisfied. Then,
\begin{align*}
G_n = \frac{1}{\sqrt{n}} \sum_{j=1}^n (X_{n,j} - \mu^{(j)} ) \rightsquigarrow Z
\end{align*}
in $C(T)$, where $Z$ is a centered Gaussian random variable with covariance function
\begin{align}\label{lVar}
C(s,t) = \mathrm{Cov}(Z(s),Z(t)) = \sum_{i=-\infty}^{\infty} \gamma(i,s,t).
\end{align}
\end{theorem}

Assumption (A5) will be used to verify a weak invariance principle for the process $(\hat{\mathbb{V}}_n\colon n\in\mathbb{N})$ given by
\begin{align} \label{39}
\hat{\mathbb{V}}_n(s) = \frac{1}{\sqrt{n}} \sum^{\lfloor sn \rfloor}_{j=1}
(X_{n,j} - \mu^{(j)})
+ \sqrt{n} \Big (s- \frac {\lfloor sn \rfloor}{n} \Big )
\big(X_{\lfloor sn \rfloor+1} - \mu^{(\lfloor sn \rfloor+1)} \big),
\end{align}
useful for the change-point analysis proposed in Section \ref{sec:change-point}. Note that the process $(\hat{\mathbb{V}}_n(s)\colon s \in [0,1])$ is an element of the Banach space $C([0,1], C(T))  = \{ \phi \colon [0,1] \rightarrow C(T) ~|~  \phi  \mbox { is   continuous} \}$, where the norm on this space is given by
\begin{equation} \label{311}
\sup_{s \in [0,1]}  \sup_{t \in T} |  \phi(s,t)|
=\|\phi\|_{C([0,1] \times T )}
\end{equation}
(note that each element of $C \big ([0,1], C(T) \big)$ can equivalently be regarded as an element of $C([0,1] \times T)$). Here and throughout this paper the notation $\|\cdot\|$ is used to denote any of the arising $\sup$-norms as the corresponding space can be identified from the context. The proof of the following result is postponed to Section \ref{proofssec2}  of the online supplement.

{\begin{theorem} \label{WIP}
Let $(X_{n,j}\colon n\in\mathbb{N}, j=1,\dots,n)$ denote an array of $C([0,1])$-valued random variables such that Assumption \ref{as:ts} is satisfied. Then, the weak invariance principle holds, that is,
\begin{align} \label{312}
\hat{\mathbb{V}}_n \rightsquigarrow \mathbb{V}
\end{align}
in $C([0,1] \times T)$, where  $\mathbb{V}$ is a centered Gaussian measure on $C([0,1] \times T)$ characterized by
\begin{align} \label{312a}
\mathrm{Cov} \big (\mathbb{V}(s,t), \mathbb{V}(s^\prime,t^\prime) \big)
= (s \wedge s^\prime) C(t,t^\prime),
\end{align}
and  the long-run covariance function $C$ is given in \eqref{lVar}.
\end{theorem}}

%%%%%%%%%%%%%%%%%%%%%%%%
 \section{The two-sample problem}
 \label{sec:two-sample}
 %%%%%%%%%%%%%%%%%%%%%%%%

From now on, consider the case $T=[0,1]$, as this is the canonical choice for functional data analysis. Two-sample problems have a long history in statistics and the corresponding tests are among the most applied statistical procedures. For the functional setting, there have been a number of contributions as well. Two are worth mentioning in the present context. Hall and Van Keilegom \cite{hallVanKeilegom2007} studied the effect of smoothing when converting discrete observations into functional data. Horv\'ath et al.\ \cite{horvathKokoszkaReeder2013} introduced two-sample tests for $L^p$-$m$ approximable functional time series based on Hilbert-space theory. In the following, a two-sample test is proposed in the Banach-space framework of Section \ref{sec:banach_methods}. To this end, consider two independent samples $X_1, \ldots, X_m$ and $Y_1, \ldots, Y_n$ of $C([0,1])$-valued random variables. Under (A2) in Assumption \ref{as:ts} expectation functions and covariance kernels exist and are denoted by $\mu_1=\mathbb{E} [X_1]$ and $\mu_2=\mathbb{E} [Y_1]$, and $k_1(t,t^\prime)=\mbox{Cov} (X_1(t), X_1(t^\prime))$ and $k_2(t,t^\prime)=\mbox{Cov} (Y_1(t), Y_1(t^\prime))$, respectively. Interest is then in the size of
the maximal deviation
\[
d_\infty=\|\mu_1-\mu_2\| =  \sup_{t \in [0,1] } | \mu_1 (t)-\mu_2 (t) |
\]
between the two mean curves, that is, in testing the hypotheses of a relevant difference
\begin{equation}
\label{H0}
H_0\colon d_\infty
\leq \Delta
\qquad \mbox{versus} \qquad
H_1\colon d_\infty > \Delta,
\end{equation}
where $\Delta \geq 0$ is a pre-specified constant determined by the user of the test. Note again that the ``classical'' two-sample problem $H_0\colon \mu_1=\mu_2$ versus $H_0\colon \mu_1\not=\mu_2$   -- which, to the best of our knowledge, has not been investigated for $C([0,1])$-valued data yet  -- is contained in this setup as the special case $\Delta=0$. Observe also that tests for relevant differences between two finite-dimensional parameters corresponding to different populations have been considered mainly in the biostatistical literature, for example in Wellek \cite{wellek2010}. It is assumed throughout this section that the samples are balanced in the sense that
\begin{equation}
\label{32}
\frac {m}{n+m} \longrightarrow \lambda \in (0,1)
\end{equation}
as $m,n \to \infty$. Additionally, let $X_1\ldots,X_m$ and $Y_1,\ldots,Y_n$ be sampled from independent time series $(X_j\colon j\in\mathbb{N})$ and $(Y_j\colon j\in\mathbb{N})$ that satisfy conditions (A1)--(A4) of Assumption \ref{as:ts}. Under these conditions both functional time series satisfy the CLT and it then follows from Theorem \ref{mixingCLT} that
\begin{align} \label{2convergence1}
\frac{\sqrt{n+m}}{m} \sum_{j=1}^m (X_j - \mu_1)
\rightsquigarrow \frac{1}{\sqrt{\lambda}} ~ Z_1
\qquad\mbox{and}\qquad
\frac{\sqrt{n+m}}{n} \sum_{j=1}^n (Y_j - \mu_2)
\rightsquigarrow \frac{1}{\sqrt{1-\lambda}} ~ Z_2,
\end{align}
where $Z_1$ and $Z_2$ are independent, centered Gaussian processes possessing covariance functions
\[
C_1(t,t^\prime) = \sum_{j=-\infty}^\infty \gamma_1(j,t,t^\prime)
\qquad\mbox{and}\qquad
C_2(t,t^\prime) = \sum_{j=-\infty}^\infty \gamma_2(j,t,t^\prime),
\]
respectively. Here $\gamma_1$ and $\gamma_2$, correspond to the respective sequences $(X_j\colon j\in\mathbb{N})$ and $(Y_j\colon j\in\mathbb{N})$ and are defined in Assumption~\ref{as:ts}. Now, the weak convergence in \eqref{2convergence1} and the independence of the samples imply immediately that
\begin{equation}
\label{33}
Z_{m,n} = \sqrt{n+m} \Big( \frac{1}{m} \sum^m_{j=1} X_j - \frac {1}{n} \sum^n_{j=1} Y_j - (\mu_1 - \mu_2) \Big)
\rightsquigarrow Z
\end{equation}
in $C([0,1])$ as $m,n\to\infty$, where $Z=Z_1/\sqrt{\lambda}+Z_2/\sqrt{1-\lambda}$ is a centered Gaussian process with covariance function
\begin{equation}\label{34}
C(t,t^\prime)=\mathrm{Cov}(Z(t),Z(t^\prime))=\frac{1}{\lambda}C_1(t,t^\prime)+\frac{1}{1-\lambda}C_2(t,t^\prime).
\end{equation}
Under the convergence in \eqref{33} the statistic
\begin{equation}
\label{2statistic}
%\label{35}
\hat d_\infty = %\sup_{t \in [0,1]}
\Big\| \frac {1}{m} \sum^m_{j=1} X_j - \frac {1}{n} \sum^n_{j=1} Y_j \Big\|
\end{equation}
is a reasonable estimator of the maximal deviation $d_\infty = \| \mu_1 - \mu_2\| $, and the null hypothesis in \eqref{H0} is rejected for large values of $\hat d_\infty$. In order to develop a test with a pre-specified asymptotic level, the limit distribution of $\hat d_\infty$ is determined in the following. For this purpose, let
\begin{equation}
\label{35}
\mathcal{E}^{\pm}=  \big \{ t \in [0,1]  \colon \mu_1(t) - \mu_2(t) = {\pm} d_\infty \big  \}
\end{equation}
if $d_\infty >0$, and define $\mathcal{E}^{+} = \mathcal{E}^{-}=[0,1]$ if $d_\infty =0$. Finally, denote by  $\mathcal{E}= \mathcal{E}^+ \cup \mathcal{E}^-$ the set of extremal points  of the difference  $\mu_1 - \mu_2$ of the two mean functions. The first main result establishes the asymptotic distribution of the statistic $\hat d_\infty$.

\begin{theorem}
\label{thm1}
If $X_1,\ldots,X_m$ and $Y_1,\ldots,Y_n$ are sampled from independent time series $(X_j\colon j\in\mathbb{N})$ and $(Y_j\colon j\in\mathbb{N})$ in $C([0,1])$, each satisfying conditions (A1)--(A4) of Assumption \ref{as:ts}, then
\begin{equation} \label{teps}
T_{m,n}=\sqrt{n+m} (\hat d_\infty - d_\infty)
\stackrel{\mathcal{D}}{\longrightarrow}
T({ \cal E})
=\max \Big\{ \sup_{t \in { \cal E}^+} Z(t), \sup_{t \in  { \cal E}^-} - Z(t) \Big\},
\end{equation}
where the  centered Gaussian process  $Z$ is given by \eqref{34} and the sets  $\mathcal{E}^{+}$ and $ \mathcal{E}^{-}$ are defined in
\eqref{35}.
\end{theorem}

It should be emphasized that the limit distribution depends in a complicated way on the set $\cal E$ of extremal points of the difference $\mu_1 -\mu_2$ and is therefore not distribution free, even in the case of i.i.d.\ data. In particular, there can be two sets of processes with corresponding mean functions $\mu_1,\mu_2$ and $\tilde \mu_1,\tilde  \mu_2$ such that  $\|\mu_1 - \mu_2  \| = \| \tilde \mu_1 - \tilde  \mu_2\|$. However, the respective limit distributions in Theorem \ref{thm1} will be entirely different if the corresponding sets of extremal points  ${\cal E} $ and $\tilde {\cal E}$ do not coincide. The proof of Theorem \ref{thm1} is given in Section \ref{subsec:proof:two-sample}  of the  online supplement.  In the case $d_\infty =0$, ${\cal E}^+ ={ \cal E}^- =[0,1]$ and it follows for the random variable  $T({ [0,1]})$ in Theorem \ref{thm1} that
 \begin{equation} \label{310a}
 T = \max_{t \in [0,1] }| Z(t) |.
 \end{equation}
Here  the result is a simple consequence of the  weak convergence  \eqref{33} of the process $Z_{m,n}$ (see Theorem \ref{mixingCLT}) and the continuous mapping theorem.

However,  Theorem \ref{thm1} provides also the distributional properties of the statistic $\hat d_\infty$ in the case $d_\infty>0$. This is required for testing the hypotheses of a relevant difference between the two mean functions (that is, the hypotheses in \eqref{H0} with $\Delta >0$) of primary interest here. In this case the weak convergence of an appropriately standardized version of  $\hat d_\infty$ does {\it not} follow from the weak convergence \eqref{33}, as the process inside the supremum in \eqref{2statistic} is not centered. In fact, additional complexity enters in the proofs because even under the null hypothesis observations cannot be easily centered. For details, refer to Section \ref{subsec:proof:two-sample} of the online supplement.

% %%%%%%%%%%%%%%%%%%%%%%%%
 \subsection{Asymptotic inference}
 \label{asyminf}
%%%%%%%%%%%%%%%%%%%%%%%%%

%%%%%%%%%%%%%%%%%%%%%%%%
\subsubsection{Testing the classical hypothesis $H_0\colon\mu_1\equiv \mu_2$}
 \label{testdiff}
 %%%%%%%%%%%%%%%%%%%%%%%%

Theorem \ref{thm1} also provides the asymptotic distributions of the test statistic $\hat d_\infty$ in the case of two identical mean functions, that is, if $\mu_1 \equiv \mu_2$. This is the situation investigated in Hall and Van Keilegom \cite{hallVanKeilegom2007} and Horv\'ath et al.\ \cite{horvathKokoszkaReeder2013} in Hilbert-space settings. Here it corresponds to the special case $\Delta=0$ and thus $d_\infty=0, \mathcal{E}^\pm=[0,1]$. Consequently,
\[
T_{m,n} \stackrel{\mathcal{D}}{\longrightarrow} T
\qquad(m,n\to\infty),
 \]
where the random variable $T$ is defined in \eqref{310a}. An asymptotic level $\alpha$ test for the classical hypotheses
\begin{equation} \label{hypclass}
H_0\colon \mu_1 = \mu_2 \qquad \mbox{versus} \qquad H_1\colon \mu_1 \neq \mu_2
\end{equation}
may hence be obtained by rejecting $H_0$ whenever
\begin{equation} \label{testclass}
\hat d_\infty > \frac {u_{1- \alpha}}{\sqrt{n+m}},
\end{equation}
where $u_{1- \alpha}$ is the $(1- \alpha)$-quantile of the distribution of the random variable $T$ defined in \eqref{310a}. Using Theorem \ref{thm1} it  is easy to see that the test defined by \eqref{testclass} is consistent and  has asymptotic level $\alpha$.

%%%%%%%%%%%%%%%%%%%%%%%%
 \subsubsection{Confidence bands}
 \label{sec:conf_bands}
%%%%%%%%%%%%%%%%%%%%%%%%

The methodology developed so far  can  easily be applied to the construction of simultaneous asymptotic confidence bands for the difference of the mean functions.  There is a rich literature on confidence bands for functional data in Hilbert spaces. The available work includes Degras~\cite{degras2011}, who dealt with confidence bands for nonparametric regression with functional data; Cao et al.\ \cite{caoYangTodem2012}, who studied simultaneous confidence bands for the mean of dense functional data based on polynomial spline estimators; Cao \cite{cao2014}, who developed simultaneous confidence bands for derivatives of functional data when multiple realizations are at hand for each function, exploiting within-curve correlation; and Zheng et al.\ \cite{zhengYangHardle2014} who treated the sparse case. Most recently Choi and Reimherr \cite{choiReimherr2016} extracted geometric features akin to Mahalanobis distances to build confidence bands for functional parameters.

The results presented here are the first of their kind relating to Banach space-valued functional data. The first theorem uses the limit distribution obtained in Theorem \ref{thm1} to construct asymptotic simultaneous confidence bands for the two-sample case.  A corresponding bootstrap analog will be developed in the next section. Confidence bands for the one-sample case can be constructed in a similar fashion using standard arguments and the corresponding results are consequently omitted.

\begin{theorem}
\label{th:cb-1}
Let the assumptions of Theorem \ref{thm1} be satisfied and, for $\alpha\in(0,1)$, denote by $u_{1- \alpha}$ the  $(1-\alpha)$-quantile of the random variable
%constant such that $\mathbb{P}(
$T$ defined in \eqref{310a} and define the functions
\begin{align*}
\mu^\pm_{m,n} (t) =\frac{1}{m} \sum_{j=1}^m X_j  - \frac{1}{n}  \sum_{j=1}^n Y_j
\pm  \frac{u_{1-\alpha}}{\sqrt{n+m}}.
\end{align*}
Then the set $C_{\alpha,{m,n}} = \big \{ \mu  \in C( [0,1])~\colon~\mu^-_{m,n}(t) \leq  \mu (t) \leq \mu^+_{m,n}(t) ~\mbox{ for all } t \in [0,1] \big\}$ defines a simultaneous asymptotic $(1-\alpha)$ confidence band for $\mu_1-\mu_2$, that is,
\begin{align*}
\lim_{m,n \to\infty} \mathbb{P} ( \mu_1-\mu_2 \in C_{\alpha,{m,n}} )
= 1-\alpha.
\end{align*}
\end{theorem}
Note that, unlike their Hilbert-space counterparts, the simultaneous confidence bands given in Theorem \ref{th:cb-1} (and their bootstrap analogs in Section \ref{sec:conf_bands}) hold for all $t\in[0,1]$ and not only almost everywhere, making the proposed bands more easily interpretable and perhaps more useful for applications.

%%%%%%%%%%%%%%%%%%%%%%%%
\subsubsection{Testing for a relevant difference}
\label{testdiffrel}
%%%%%%%%%%%%%%%%%%%%%%%%

Recall the definition of the random variable $T({\cal E})$ in Theorem \ref{thm1}, then
 the null hypothesis of no relevant difference in \eqref{H0} is rejected at level $\alpha$, whenever the inequality
\begin{equation} \label{221}
\hat d_\infty > \Delta + \frac {u_{1 - \alpha, { \cal E}}}{\sqrt{n+m}}
\end{equation}
holds, where $u_{\alpha, { \cal E}}$ denotes the $\alpha$-quantile of the distribution of $T({ \cal E})$ $(\alpha\in(0,1) )$.  A conservative test avoiding the use  of quantiles depending on the set of extremal points $ { \cal E}$ can be obtained observing the inequality
\begin{equation}
\label{21}
T({ \cal E}) \leq T,
\end{equation}
where the random variable $T$ is defined  in \eqref{310a}. If $u_\alpha$ denotes the $\alpha$-quantile of the distribution of  $T$, then \eqref{21} implies $u_{\alpha, { \cal E}} \leq u_\alpha$ and a conservative asymptotic level $\alpha$ test is given by rejecting the  null hypothesis in \eqref{H0}, whenever the inequality
\begin{equation} \label{22}
\hat d_\infty > \Delta + \frac {u_{1 - \alpha}}{\sqrt{n+m}}
\end{equation}
holds. The properties of the tests \eqref{221} and \eqref{22} depend on the size of the distance $d_\infty$ and will be explained below. In particular, observe the following properties for the test  \eqref{22}:
\begin{itemize}
\item [(a)] {If $d_\infty < \Delta$, Slutsky's theorem yields that
\begin{align*}
\lim_{n,m \to \infty} \mathbb{P} \Big(
& \hat d_\infty > \Delta + \frac {u_{1 - \alpha}}{\sqrt{n+m}} \Big) \nonumber
=\lim_{n,m \to \infty}  \mathbb{P}\big(\sqrt{n+m}(\hat d_\infty-d_\infty)>\sqrt{n+m}(\Delta-d_\infty)+u_{1-\alpha}\big)
=0.
\end{align*}
}
\item [(b)] {If $d_\infty = \Delta$, we have} 
\begin{align}
\limsup_{n,m \to \infty}
\mathbb{P} \Big( \hat d_\infty > \Delta + \frac {u_{1 - \alpha}}{\sqrt{n+m}} \Big)
&=\limsup_{n,m \to \infty} \mathbb{P} \big(\sqrt{n+m} (\hat d_ \infty - d_\infty) > \sqrt{n+m} (\Delta - d_\infty) + u_{1 - \alpha} \big) \nonumber \\
&\leq \lim_{n,m \to \infty}\mathbb{P} \big(\sqrt{n+m} \ ( \hat d_\infty - d_\infty) >  u_{1-\alpha, { \cal E}}\big) =  \alpha
 \label{2levelAlpha}.
\end{align}
\item [(c)]
{If   $d_\infty  > \Delta$, the same calculation as in (a) implies }
\begin{align*}
\lim_{n,m \to \infty} \mathbb{P} \Big(
& \hat d_\infty > \Delta + \frac {u_{1 - \alpha}}{\sqrt{n+m}} \Big) \nonumber
=1,
\end{align*}
proving that the test defined in \eqref{22} is consistent.
\item[(d)]  If the mean functions $\mu_1$ and $\mu_2$ define a boundary point of the hypotheses, that is, $d_\infty = \Delta$ and either $\mathcal{E}^+= [0,1]$ or $\mathcal{E}^-= [0,1]$, then $T(\mathcal{E}) = \max_{t\in[0,1]} Z(t)$  or $T(\mathcal{E}) = \max_{t\in[0,1]} - Z(t)$, and consequently
\[
\lim_{n,m \to \infty} \mathbb{P} \Big ( \hat d_\infty > \Delta + \frac {u_{1- \alpha}}{\sqrt{m+n}} \Big ) = \frac{\alpha}{2}.
\]
\end{itemize}
Using similar arguments it can be shown that the test \eqref{221}  satisfies
\[
\lim_{n,m \to \infty}
\mathbb{P} \Big( \hat d_\infty > \Delta + \frac {u_{1 - \alpha,{\cal E}}}{\sqrt{n+m}} \Big)  =
\begin{cases}
0 & \text{ if }  d_\infty < \Delta. \\
\alpha & \text{ if }  d_\infty = \Delta.  \\
1  & \text{ if }  d_\infty \geq \Delta.
\end{cases}
\]
Summarizing, the tests  for the hypothesis \eqref{H0} of no relevant difference between the two mean functions
defined in \eqref{221}  and \eqref{22}  have asymptotic level $\alpha$ and are  consistent.
However, the discussion given above also shows that the test \eqref{22} is conservative, even
when $\mathcal{E} =  [0,1]$.

%%%%%%%%%%%%%%%%%%%%%%%%%
\subsection{Bootstrap}
\label{subsec:two-sample:bootstrap}
%%%%%%%%%%%%%%%%%%%%%%%%%

In order to use the tests \eqref{testclass},  \eqref{221}  and  \eqref{22}  for  classical and relevant  hypotheses, the quantiles of the distribution of the random variables $T(\mathcal{E}) $ and   $T$ defined in \eqref{teps} and \eqref{310a} need to be estimated, which depend on certain features of the data generating process. The law $T(\mathcal{E}) $ involves the unknown set of extremal points $\mathcal{E}$ of the differences of the mean functions. Moreover, the distributions of $T(\mathcal{E}) $ and $T$ depend on the long-run covariance function \eqref{34}. There are methods available in the literature to consistently estimate the covariance function (see, for example, Horv\'ath et al. \cite{horvathKokoszkaReeder2013}). In practice, however, it is difficult to reliably approximate the infinite sums in \eqref{34} and therefore an easily implementable  bootstrap procedure is proposed in the following.

It turns out that a different and non-standard bootstrap procedure will be required for testing relevant hypotheses than for classical hypotheses (and the construction of confidence bands) as in this case the null distribution depends on the set of extremal points $\mathcal{E}$. The corresponding resampling procedure requires a substantially more sophisticated analysis. Therefore the analysis of bootstrap tests for the classical hypothesis and bootstrap confidence intervals is given first and discussion of bootstrap tests for relevant hypotheses is deferred to Section \ref{bootrel}.

\subsubsection{Bootstrap confidence intervals and tests for the classical hypothesis $H_0: \mu_1=\mu_2$}
\label{sec:conf_bands}

Following  B\"ucher and Kojadinovic \cite{buecher2016} the use of a muliplier block bootstrap is proposed. To be precise, let $(\xi_k^{(1)}\colon k\in\mathbb{N}), \ldots , (\xi_k^{(R)}\colon k\in\mathbb{N})$ and $(\zeta_k^{(1)}\colon k\in\mathbb{N}), \ldots ,(\zeta_k^{(R)}\colon k\in\mathbb{N})$ denote independent sequences of independent standard normally distributed random variables and define the $C([0,1])$-valued processes $\hat B_{m,n}^{(1)}, \ldots , \hat B_{m,n}^{(R)}$ through
\begin{align} \label{2bProcess}
\begin{split}
\hat B_{m,n}^{(r)}(t) = \sqrt{n+m} \Big\{&
	\frac{1}{m} \sum_{k=1}^{m-l_1+1} \frac{1}{\sqrt{l_1}}\Big( \sum_{j=k}^{k+l_1-1} X_{j}(t)
	-\frac{l_1}{n}\sum_{j=1}^m X_{j}(t) \Big) \xi_k^{(r)} \\
&+ \frac{1}{n} \sum_{k=1}^{n-l_2+1} \frac{1}{\sqrt{l_2}}\Big( \sum_{j=k}^{k+l_2-1} Y_{j}(t)
	-\frac{l_2}{n}\sum_{j=1}^n Y_{j}(t) \Big) \zeta_k^{(r)} \Big\}  ~~~~(r=1, \ldots , R)
\end{split}
\end{align}
for $t\in[0,1]$, where $l_1,l_2\in\mathbb{N}$ denote window sizes such that $l_1/m\to 0$ and $l_2/n\to 0$ as $l_1,l_2,m,n\to\infty$. The following result is  a fundamental tool for the theoretical investigations of all bootstrap procedures proposed in this paper and is proved in Section \ref{subsec:proof:two-sample} of the online supplement. 

\begin{theorem} \label{2bTheorem}
Suppose that $(X_j\colon j\in\mathbb{N})$ and $(Y_j\colon j\in\mathbb{N})$ satisfy conditions (A1)--(A4) of Assumption \ref{as:ts} and
 let  $\hat B_{m,n}^{(1)}, \ldots , \hat B_{m,n}^{(R)}$  denote the bootstrap processes  defined by \eqref{2bProcess} such that $l_1 = m^{\beta_1}$, $l_2 = n^{\beta_2}$ with $0<\beta_i<\nu_i/(2+\nu_i)$ and $\nu_i$ given in Assumption \ref{as:ts}, $i=1,2$. Then,
\begin{align*}
(Z_{m,n},\hat{B}_{m,n}^{(1)},\dots,\hat{B}_{m,n}^{(R)})
\rightsquigarrow (Z, Z^{(1)},\dots,Z^{(R)})
\end{align*}
in $C([0,1])^{R+1}$ as $m,n\to\infty$, where $Z_{m,n}$ is defined in \eqref{33} and $Z^{(1)},\dots,Z^{(R)}$ are independent copies of
the centered  Gaussian process $Z$ defined by \eqref{34}.
\end{theorem}

Note that 
Theorem \ref{2bTheorem} holds  under the null hypothesis and alternative. It
 leads to the following  results regarding confidence bands and tests for the classical hypothesis \eqref{hypclass} based on the  the multiplier bootstrap. To this end, note that for the statistics
$$
T_{m,n}^{(r)} =  \|\hat{B}_{m,n}^{(r)} \|,
\qquad r=1 \ldots , R ,
$$
the continuous mapping theorem yields
\begin{align} \label{classicH0:convergence}
\big( \sqrt{n+m} ~ \hat{d}_\infty ,~ T_{m,n}^{(1)},\dots,T_{m,n}^{(R)}\big)
\Rightarrow (T,~ T^{(1)},\dots,T^{(R)}),
\end{align}
where the random variables $T^{(1)},\dots,T^{(R)}$ are independent copies of the statistic  $T$ defined in \eqref{310a}. Now, if $T_{m,n}^{\{\lfloor R(1-\alpha)\rfloor\}}$ is the empirical $(1-\alpha)$-quantile of the bootstrap sample $T_{m,n}^{(1)},\dots,T_{m,n}^{(R)}$, the following results are obtained.
\begin{theorem}
\label{th:cb-2}
Let the assumptions of Theorem \ref{2bTheorem} be satisfied and define the functions
\begin{align*}
\mu^{R,\pm }_{m,n}  (t) =\frac{1}{m} \sum_{j=1}^m X_j  - \frac{1}{n}  \sum_{j=1}^n Y_j
\pm  \frac{T_{m,n}^{\{\lfloor R(1-\alpha)\rfloor\}}}{\sqrt{n+m}}.
\end{align*}
Then,
$\hat{C}_{\alpha,m.n} ^R  =  \{ \mu  \in C( [0,1])\colon\mu^{R,-}_{m,n} (t) \leq  \mu (t) \leq \mu^{R,+}_{m,n} (t)~\mbox{for all}~t \in [0,1]\}$ defines a simultaneous asymptotic $(1-\alpha)$ confidence band for $\mu_1-\mu_2$, that is,
\begin{align*}
\lim_{R \to \infty }
\liminf_{m,n \to\infty} \mathbb{P} ( \mu_1-\mu_2 \in\hat{C}_{\alpha,m.n} ^R   )
\geq  1-\alpha.
\end{align*}
\end{theorem}
This section is concluded with a corresponding statement regarding the bootstrap test for the classical hypotheses in \eqref{hypclass}, which rejects the null hypothesis whenever
\begin{align} \label{bootclasstest}
\hat{d}_{\infty} > \frac{T_{m,n}^{\{\lfloor R(1-\alpha)\rfloor\}}}{\sqrt{n+m}},
\end{align}
where the statistic  $\hat{d}_{\infty}$  is defined in \eqref{2statistic}.

\begin{theorem}\label{bootstrap:classic}
Let the assumptions of Theorem \ref{2bTheorem} be satisfied, then the test \eqref{bootclasstest} has asymptotic level $\alpha$ and is consistent for the hypotheses \eqref{hypclass}. More precisely, under the null hypothesis of no difference in the mean functions,
\begin{align} \label{bootclasstestLVL}
\lim_{R\to\infty} \limsup_{m,n\to\infty} \mathbb{P}\bigg( \hat{d}_{\infty}
	> \ \frac{T_{m,n}^{\{\lfloor R(1-\alpha)\rfloor\}}}{\sqrt{n+m}} \bigg)
	= \alpha,
\end{align}
and, under the alternative,  for any \(R \in \mathbb{N}\),
%we have
%similar calculations as in \eqref{2consistency} imply
\begin{align} \label{bootclasstestCons}
 \liminf_{m,n\to\infty} \mathbb{P}\bigg( \hat{d}_{\infty}
	>  \frac{T_{m,n}^{\{\lfloor R(1-\alpha)\rfloor\}}}{\sqrt{n+m}} \bigg)
	=1.
\end{align}
\end{theorem}

%%%%%%%%%%%%%%%%%%%%%%%%%
\subsubsection{Testing for relevant differences in the mean functions}
\label{bootrel}
%%%%%%%%%%%%%%%%%%%%%%%%%

The problem of constructing an appropriate bootstrap test for the hypotheses of no relevant difference in the mean functions is substantially more complicated. The reason for these difficulties consists in the fact that in the case of relevant  hypotheses the limit distribution of the corresponding test statistic is complicated. In contrast to the problem of testing the classical hypotheses \eqref{hypclass}, where it is sufficient to mimic the distribution of the statistic
$T$  in \eqref{310a} (corresponding to the case $\mu_1\equiv \mu_2$) one requires the distribution of the statistic $T(\mathcal{E} )$, which depends in a sophisticated  way on the set of extreme points of the (unknown) difference $\mu_1 - \mu_2$. Under the null hypothesis $\| \mu_1 - \mu_2 \|  \leq \Delta $ these sets can be very different,
ranging from a singleton to the full interval $[0,1]$. As a consequence the construction of a valid bootstrap procedure requires appropriate consistent estimates of the sets $\mathcal{E}^+$ and $\mathcal{E}^-$ introduced in Theorem \ref{thm1}.

For this purpose, recall the definition of the Haussdorff distance between two sets $A,B \subset \mathbb{R} $, given by
$$
d_H(A,B) =  \max \Big \{ \sup_{x \in A} \inf_{y \in B} |x-y| ,  \sup_{y \in B}  \inf_{x \in A}  |x-y|  \Big\}
$$
and denote by $K([0,1])$ the set of all compact subsets of the  interval $[0,1]$. First, define estimates of the extremal sets  $\mathcal{E}^+$ and $\mathcal{E}^-$ by
\begin{align} \label{estimatedSets}
\hat{\mathcal{E}}_{m,n}^\pm  &= \Big\{ t\in[0,1] \colon   \pm (\hat{\mu}_1(t)-\hat{\mu}_2(t)  )
\geq \hat{d}_\infty - \frac{c_{m,n}}{\sqrt{m+n}} \ \Big\},
%\nonumber \\
%\hat{\mathcal{E}}_{m,n}^- &= \Big\{ t\in[0,1] \ \Big| \ \hat{\mu}_1(t)-\hat{\mu}_2(t)
% \leq -\hat{d}_\infty + \frac{c_{m,n}}{\sqrt{m+n}} \ \Big\}
\end{align}
where $c_{m,n} \sim \log (m+n) $. Our first result shows that the estimated sets $\hat{\mathcal{E}}_{m,n}^+$ and $\hat{\mathcal{E}}_{m,n}^-$ are consistent for $\mathcal{E}^+$ and $\mathcal{E}^-$, respectively.

\begin{theorem} \label{setConvergence}
Let the assumptions of Theorem \ref{2bTheorem} be satisfied, then
\begin{align*}
d_H( \hat{\mathcal{E}}_{m,n}^\pm ,  \mathcal{E}^\pm) \xrightarrow{\mP}  0 ,
\end{align*}
 where the sets $\hat{\mathcal{E}}_{m,n}^\pm$ are defined by
\eqref{estimatedSets}.% and  $\mP^* (\mathcal{A} ) $ denotes the outer probability of the set $\mathcal{A} $.
\end{theorem}

\bigskip

The main implication of Theorem  \ref{setConvergence} consists  in the fact that the random variable
$$ \max_{t\in \hat{\mathcal{E}}_{m,n}^+}
\hat{B}_{m,n}(t) $$
converges weakly to the random variable $\max_{t \in  \mathcal{E}^+} Z(t) $.
Note that $\hat{B}_{m,n} \rightsquigarrow Z $
by Theorem \ref{2bTheorem} and that $d_H ( \hat{\mathcal{E}}_{m,n}^+ , \hat{\mathcal{E}}^+) \to 0 $
in probability by the previous theorem, but the combination of both statements is more delicate and
requires a continuity argument which is given in Section \ref{subsec:proof:two-sample} of the online supplement, where the following result is proved.

\begin{theorem} \label{2jointConvergence}
Let the assumptions of Theorem \ref{2bTheorem} be satisfied and define, for $r=1,\dots,R$,
\begin{align} \label{bootstat}
K^{(r)}_{m,n}
=  \max\Big\{ \max_{t\in \hat{\mathcal{E}}_{m,n}^+}
\hat{B}_{m,n}^{(r)}(t), \ \max_{t\in \hat{\mathcal{E}}_{m,n}^-}
\big(- \hat{B}_{m,n}^{(r)}(t) \big) \Big\}.
\end{align}
Then,
\begin{align} \label{boottestrelweak}
\big(\sqrt{n+m} \;(\hat d_\infty - d_\infty),~ K_{m,n}^{(1)},\dots,K_{m,n}^{(R)}\big)
\Rightarrow (T(\mathcal{E}),~ T^{(1)}(\mathcal{E}),\dots,T^{(R)}(\mathcal{E})),
\end{align}
in $\mathbb{R}^{R+1}$, where $d_\infty = \| \mu_1 - \mu_2 \|$,
 the statistic $\hat d_\infty$ is defined in \eqref{2statistic}
and the variables
$T^{(1)}(\mathcal{E}),\dots,T^{(R)}(\mathcal{E})$ are independent copies of\/ $T(\mathcal{E})$ defined in Theorem \ref{thm1}.
\end{theorem}

Theorem \ref{2jointConvergence} leads to a simple bootstrap test for the hypothesis of no relevant change.  To be precise, let   $K_{m,n}^{\{\lfloor R(1-\alpha)\rfloor\}}$ denote the  empirical $(1-\alpha)$-quantile of the bootstrap sample $K_{m,n}^{(1)},\dots,K_{m,n}^{(R)}$, then the  null hypothesis of no relevant change is rejected at level $\alpha$, whenever
\begin{align}
\label{boottestrel}
\hat{d}_{\infty} > \Delta + \frac{K_{m,n}^{\{\lfloor R(1-\alpha)\rfloor\}}}{\sqrt{n+m}}.
\end{align}
The final result of this section shows that the test \eqref{boottestrel} is  consistent and has asymptotic level $\alpha$.  The proof is obtained by similar arguments as given in the proof of Theorem \ref{bootstrap:classic}, which are omitted for the sake of brevity.

\begin{theorem}\label{two-sample:bootstrap:asymptotic}
Let the assumptions of Theorem \ref{2bTheorem} be satisfied. Then, under the null hypothesis of no relevant difference in the mean functions,
\begin{align*}
\lim_{R\to\infty} \limsup_{m,n\to\infty} \mathbb{P}\bigg ( \hat{d}_{\infty}
	> \Delta + \frac{K_{m,n}^{\{\lfloor R(1-\alpha)\rfloor\}}}{\sqrt{n+m}} \bigg)
	= \alpha,
\end{align*}
and, under the alternative of a relevant difference in the mean functions, for any \(R \in \mathbb{N}\),
%similar calculations as in \eqref{2consistency} imply
\begin{align*}
 \liminf_{m,n\to\infty} \mathbb{P}\bigg( \hat{d}_{\infty}
	> \Delta + \frac{K_{m,n}^{\{\lfloor R(1-\alpha)\rfloor\}}}{\sqrt{n+m}} \bigg)
	=1.
\end{align*}
\end{theorem}

%%%%%%%%%%%%%%%%%%%%%%%%
\section{Change-point analysis} %for $C([0,1])$-valued functional data}
\label{sec:change-point}
%%%%%%%%%%%%%%%%%%%%%%%%

Change-point problems arise naturally in a number of applications (for example, in quality control, economics and finance; see Aue and Horv\'ath \cite{aueHorvath2013} for a recent review). In the functional framework, applications have centered around environmental and climate observations (see Aue et al.\ \cite{aueDubartNorinhoHormann2015,AueRiceSonmez2015}) and intra-day finance data (see Horv\'ath et al.\ \cite{HorvathKokoskza2012}). One of the first contributions in the area are  Berkes et al.\ \cite{berkesGabrysHorvathKokoszka2009} and Aue et al.\ \cite{aueGabrysHorvathKokoszka2009}
who developed change-point analysis  in a Hilbert space setting for independent data. Generalizations to time series of functional data in Hilbert spaces are due to Aston and Kirch \cite{astonKirch2012a,astonKirch2012b}. For Banach-spaces, to the best of our knowledge, the only contributions to change-point analysis available in the literature are due to Ra\v{c}kauskas and Suquet \cite{rackauskasSuquet2004, rackauskasSuquet2006}, who have provided theoretical work analyzing epidemic alternatives for independent functions based on H\"older norms and dyadic interval decompositions.
This section details new results on change-point analysis for $C([0,1])$-valued functional data. The work is the first to systematically exploit a time series structure of the functions as laid out in Section \ref{sec:banach_methods}.

%%%%%%%%%%%%%%%%%%%%%%%%%
\subsection{Asymptotic inference}
\label{sec41}
%%%%%%%%%%%%%%%%%%%%%%%%%

More specifically, the problem of testing for a (potentially relevant) change-point is considered for triangular arrays $(X_{n.j}\colon n\in\mathbb{N}, j=1,\ldots,n)$ of $C([0,1])$-valued random variables satisfying Assumption \ref{as:ts}. Denote by $\mu^{(j)}=\mathbb{E}[X_{n,j}] \in C([0,1])$ the expectation of $X_{n,j}$ and assume as in part (A2) of Assumption \ref{as:ts} that $\gamma(j-j^\prime,t,t^\prime)=\mathrm{Cov}(X_{n,j}(t), X_{n,j^\prime}(t^\prime))$ is the covariance kernel common to all random functions in the sample. Parametrize with $s^* \in (\vartheta,1-\vartheta)$, where $\vartheta\in(0,1)$ is a constant, the location of the change-point, so that the sequence  $(\mu^{(j} )_{j\in \mathbb{N}}$ of mean functions  satisfies
\begin{equation} \label{36a}
\mu_1 = \mu^{(1)} = \cdots = \mu^{(\lfloor ns^* \rfloor)}
\qquad\mbox{and}\qquad
\mu_2 = \mu^{(\lfloor ns^* \rfloor + 1)} = \cdots = \mu^{(n)}.
\end{equation}
Then, for any $n\in\mathbb{N}$, both $X_{n,1},\ldots, X_{n,\lfloor ns^* \rfloor}$ and $X_{n,\lfloor ns^* \rfloor+1}, \ldots, X_{n,n}$ consist of (asymptotically) identically distributed but potentially dependent random functions. Let again $d_\infty=\|\mu_1-\mu_2\|$ denote the maximal deviation between the mean functions before and after the change-point. Interest is then in testing the hypotheses of a relevant change, that is,
\begin{equation} \label{37}
H_0\colon d_\infty \leq \Delta \qquad \mbox{versus} \qquad H_1\colon d_\infty  > \Delta,
\end{equation}
where $\Delta \geq 0 $ is a pre-specified constant.
The relevant change-point test setting may be viewed in the context of a bias-variance trade-off. In the time series setting, one is often interested in accurate predictions of future realizations. However, if the stretch of observed functions suffers from a structural break, then only those functions sampled after the change-point should be included in the prediction algorithm because these typically require stationarity. This reduction of observations, however, inevitably leads to an increased variability that may be partially offset with a bias incurred through the relevant approach: if the maximal discrepancy $d_\infty$ in the mean functions remains below a suitably chosen threshold $\Delta$, then the mean-squared prediction error obtained from predicting with the whole sample might be smaller than the one obtained from using only the non-contaminated post-change sample. In applications to financial data, the size of the allowable bias could also be dictated by regulations imposed on, say, investment strategies (Dette and Wied \cite{detteWied2015} specifically mention Value at Risk as one such example).

Recall the definition of the sequential empirical process in \eqref{39}, where  the argument $s\in[0,1]$ of this process is used to search over all potential change locations.
Note that $(\hat{\mathbb{V}}_n(s,t)\colon (s,t)\in [0,1]^2) $ can be regarded as  an element of the Banach space
%\begin{equation} \label{310}
$  C([0,1]^2)$ (see the discussion before Theorem  \ref{WIP}).
 Define the $C([0,1]^2)$-valued process
\begin{equation} \label{313}
\hat{\mathbb{W}}_n(s,t) = \hat {\mathbb{V}}_n(s,t) - s \hat{\mathbb{V}}_n(1,t),
\qquad s,t \in [0,1],
\end{equation}
then,  under Assumption \ref{as:ts}, Theorem \ref{WIP} and   the continuous mapping theorem show that
  \begin{equation} \label{315}
  \hat {\mathbb{W}}_n \rightsquigarrow \mathbb{W}
  \end{equation}
in $C([0,1]^2)$, where $\mathbb{W}(s,t)= \mathbb{V}(s,t) - s \mathbb{V}(1,t)$.
%and the Gaussian process $\mathbb{V}$ is given by \eqref{312a}.
In particular, $\mathbb{W}$ is a centered Gaussian measure on $C([0,1]^2)$ defined by
\begin{equation} \label{316}
\mathrm{Cov} (\mathbb{W}(s, t), \mathbb{W}(s^\prime,t^\prime))
= (s \wedge s^\prime - s s^\prime)C(t,t^\prime).
\end{equation}
%where the kernel $C$ is defined in \eqref{lVar}.
In order to define a test for the hypothesis of a relevant change-point defined by \eqref{37}
consider the sequential empirical process $(\hat{\mathbb{U}}_{n}\colon n\in\mathbb{N})$ on $C([0,1]^2)$ given by
\begin{equation} \label{318}
\hat {\mathbb{U}}_n(s,t) = \frac {1}{n} \Big( \sum^{\lfloor sn \rfloor}_{j=1}
X_{n,j}(t) + n\Big(s- \frac {\lfloor sn \rfloor}{n}\Big)X_{n, \lfloor sn \rfloor+1}(t)
- s \sum^n_{j=1} X_{n,j}(t) \Big).
\end{equation}
%and note that
%\[
%\mathbb{E} \Big [ \hat {\mathbb{U}}_n \Big(\frac {\lfloor sn \rfloor}{n},t \Big) \Big]
%= \frac{1}{n} \big (\lfloor (s \wedge s^*)n \rfloor
%- \lfloor sn \rfloor \lfloor s^*n \rfloor \big) (\mu_1(t) - \mu_2(t)) ~.
%\]
Evaluating its expected value shows that, in contrast to $\hat {\mathbb{W}}_n$, the process $\hat {\mathbb{U}}_n$ is  typically not centered and the equality
$\sqrt{n} ~ \hat {\mathbb{U}}_n = \hat {\mathbb{W}}_n$ holds only in the case  $\mu_1=\mu_2$.
A straightforward calculation shows
that
\[
\mathbb{E}  \big[ \hat {\mathbb{U}}_n (s,t) \big]
= \big( s \wedge s^* - s s^*  \big) \big (\mu_1 (t) - \mu_2(t) \big)  + o_\mathbb{P}(1)
\]
uniformly in $(s,t) \in [0,1]^2$.
As the function $s \mapsto s \wedge s^* - ss^*$ attains its maximum in the interval $[0,1]$ at the
point $s^*$, the statistic
\begin{equation} \label{teststat}
\mathbb{\hat  M}_n  =  \sup_{s \in [0,1]} \sup_{t \in [0,1]} | \hat {\mathbb{U}}_n (s,t) |
\end{equation}
 is a reasonable estimate of $s^*(1-s^*) ~ d_\infty = s^*(1-s^*) \|\mu_1-\mu_2\|$. It is therefore proposed to reject the null hypothesis in \eqref{37} for large values of the statistic $\mathbb{\hat  M}_n$. The following result specifies the asymptotic distribution of  $\mathbb{\hat  M}_n$.

\begin{theorem} \label{thm2}
Assume  $d_\infty >0$, $s^* \in (0,1)$ and let  $(X_{n,j}\colon n\in\mathbb{N}, j=1,\dots,n)$ be an array of $C([0,1])$-valued random variables
satisfying Assumption \ref{as:ts}. Then
\begin{equation} \label{deps}
\mathbb{D}_n = \sqrt{n} \big( \mathbb{\hat  M}_n   - s^*(1-s^*) d_\infty \big)
\stackrel{\cal D}{\longrightarrow}
D (\mathcal{E})= \max \Big \{ \sup_{t \in \mathcal{E}  ^+} \mathbb{W}(s^*,t),
		\sup_{t \in \mathcal{E} ^-} - \mathbb{W}(s^*,t) \Big \},
\end{equation}
where the statistic $\mathbb{\hat  M}_n $ is defined in \eqref{teststat},
$\mathbb{W}$ is the centered Gaussian measure on $C([0,1]^2)$
characterized by \eqref{316}, $\mathcal{E}=\mathcal{E} ^+\cup\mathcal{E}^-$ and the sets $\mathcal{E}^+$ and  $\mathcal{E} ^-$ are defined  in \eqref{35}.
\end{theorem}

The proof of Theorem \ref{thm2} is given in Section \ref{sec:proofs:change-point} of the online supplement.
The limit distribution of $\mathbb{D}_n$ is rather complicated and depends on the set $\mathcal{E}$ which might be different for functions $\mu_1 - \mu_2$ with the same sup-norm $d_\infty$ but different corresponding set $\mathcal{E}$. It is also worthwhile to mention that the condition $d_\infty >0$ is essential in Theorem \ref{thm2}. In the remaining case $d_\infty =0$ the weak convergence of $\hat {\mathbb{M}}_n$ simply follows from $\sqrt{n} \hat{\mathbb{U}}_n = \hat{\mathbb{W}}_n$, \eqref{315} and the continuous mapping theorem, that is,
\begin{align}
\label{42}
\sqrt{n} \ \mathbb{\hat M}_n
\stackrel{\cal D}{\longrightarrow}
\check {T} = \sup_{(s,t) \in [0,1]^2} |\mathbb{W}(s,t)|
%\qquad (n\to\infty).
\end{align}
whenever $d_\infty =0$.

If $d_\infty >0$, the
   true location of the change-point $s^*$ is unknown and therefore has to be estimated from the available data. The next theorem,
   which is proved  in Section \ref{sec:proofs:change-point} of the online supplement,  proposes one such estimator and specifies its large-sample behavior in form of a rate of convergence.

\begin{theorem} \label{rate}
Assume  $d_\infty >0$, $s^* \in (0,1)$ and let  $(X_{n,j}\colon n\in\mathbb{N}, j=1,\dots,n)$ be  an array of $C([0,1])$-valued random variables
satisfying Assumption \ref{as:ts}, where the random variable $M$ in Assumption (A3) is bounded. Then the estimator
\begin{align} \label{cpEstimator0}
\tilde{s} =  \arg\max_{1\leq k <n} \big\| \hat {\mathbb{U}}_n (\frac kn,\cdot ) \big\|
\end{align}
satsifies
%\begin{align*}
$|\tilde{s}-s^*| = O_{\mathbb{P}}(n^{-1}).$
%\end{align*}
\end{theorem}

\noindent
 Recall that the possible range of change locations is restricted to the open interval $(\vartheta,1-\vartheta)$ and define the modified change-point estimator
\begin{align} \label{cpEstimator}
\hat  s  = \max \big \{ \vartheta, \min \{ \tilde s, 1- \vartheta\} \big\},
\end{align}
where $\tilde{s}$ is given by  \eqref{cpEstimator0}. Since $|\hat{s}-s^*|\leq |\tilde{s}-s^*|$, it follows that $|\hat{s}-s^*| = O_{\mathbb{P}}(n^{-1})$
if $d_\infty >0$, and, if $d_\infty =0$  suppose that $\hat{s}$ converges weakly to a
$[\vartheta,1-\vartheta]$-valued random variable $s_{\max}$.

\begin{cor}  \label{setT}
Let the assumptions of Theorem \ref{rate} be satisfied and define
\begin{align} \label{statistic}
\hat{d}_\infty = \frac{ \mathbb{\hat  M}_n}{\hat{s}(1-\hat{s})}
\end{align}
as an estimator of
$d_\infty$. Then,
$
\sqrt{n}\big( \hat{d}_\infty - d_\infty \big) \Rightarrow
T (\mathcal{E})=
%\frac
{D (\mathcal{E}) }/[{s^*(1-s^*)}]
$,
where $D (\mathcal{E})$ is defined in \eqref{deps}.
\end{cor}

%%%%%%%%%%%%%%%%%%%%%%%%%
%\subsection{Testing for relevant changes}
%\label{subsec:change-point:relevant}
%%%%%%%%%%%%%%%%%%%%%%%%%

\begin{rem} \label{testcp}
{\rm A consistent level $\alpha$ test for the hypotheses  \eqref{37}
is constructed along the lines of the two-sample case discussed in Section \ref{sec:two-sample}.
\begin{itemize}
\item [(a)] Consider first the case $\Delta >0 $, that is, a relevant hypothesis.
If $d_\infty>0$, implying the existence of a change-point $s^*\in(0,1)$, then the inequality
\begin{equation} \label{variableT}
T({ \cal E})\leq T = \frac{1}{s^*(1-s^*)} \sup_{t \in [0,1]} |\mathbb{W}(s^*,t)|
\end{equation}
holds. If $u_{\alpha , \mathcal{E}}$ denotes the quantile of $T({ \cal E})$, then
$
u_{\alpha , \mathcal{E}} \leq u_\alpha
$
for all $\alpha \in (0,1)$. Consequently,
similar arguments as given in Section \ref{testdiffrel} show that the test which rejects the null hypothesis of no relevant change
if
\begin{align} \label{41}
\hat d_\infty > \Delta + \frac {u_{1 - \alpha}}{\sqrt{n}}
\end{align}
% under $H_0\colon d_\infty \leq \Delta$ and $d_\infty >0 $
is  consistent and has  asymptotic level $\alpha$.
Note that an estimator of the long-run covariance function is needed in order to obtain
 the $\alpha$-quantile $u_\alpha$ of the distribution of $T$.  Moreover, the test  \eqref{41} is conservative, even when the set $\mathcal{E}$ of extremal points of the unknown difference $\mu_1 - \mu_2$ is the whole interval $ [0,1]$ (in this case the level is in fact $\alpha/2$ instead of $\alpha$ - 
 see the discussion at the end of  Section \ref{testdiffrel}).
% whenever the set $\mathcal{E}$ of extremal points of the unknown difference $\mu_1 - \mu_2$ is a proper subset of the interval $ [0,1]$.

\item [(b)]In the case of testing the classical hypotheses $H_0 \colon \mu_1 = \mu_2$ versus $H_1\colon \mu_1 \neq \mu_2$, that is
$\Delta=0$,  the test described in \eqref{41} needs to be slightly altered. The asymptotic distribution of $\mathbb{\hat M}_n $ under $H_0$ can be obtained from \eqref{42} and now it can be seen that rejecting $H_0$ whenever
\[
\hat d_\infty  > \frac {\check {u}_{1- \alpha}}{\sqrt{n}} ,
\]
where $\check{u}_{1- \alpha} $ denotes the $ (1-\alpha)$-quantile of the distribution of the random variable $\check T$
defined by \eqref{42}, yields a consistent asymptotic level $\alpha$ test.
\end{itemize}
 }
 \end{rem}

%%%%%%%%%%%%%%%%%%%%%%%%%
\subsection{Bootstrap}
\label{sec42}
%%%%%%%%%%%%%%%%%%%%%%%%%

In order to avoid the difficulties mentioned in the previous remark, a bootstrap procedure is developed and its consistency is shown.
To be precise, denote by
\begin{align*}
\hat{\mu}_1 = \frac{1}{\lfloor \hat{s} n\rfloor}
	\sum_{j=1}^{\lfloor \hat{s}n \rfloor} X_{n,j}
	\qquad \text{and} \qquad
\hat{\mu}_2 = \frac{1}{\lfloor (1-\hat{s}) n\rfloor}
	\sum_{j=\lfloor \hat{s}n \rfloor +1}^{n} X_{n,j}
\end{align*}
estimators for the expectation before and after the change-point.
Let  $(\xi_k^{(1)}\colon k\in\mathbb{N}),\ldots,(\xi_k^{(R)}\colon k\in\mathbb{N})$ denote $R$ independent sequences of independent standard normally distributed random variables
and consider the $C([0,1]^2)$-valued processes  $\hat B_n^{(1)},\ldots , \hat B_n^{(R)}$  defined  by
\begin{align} \label{bProcess}
\begin{split}
\hat{B}_n^{(r)}(s,t) =& \frac{1}{\sqrt{n}} \sum_{k=1}^{\lfloor sn \rfloor} \frac{1}{\sqrt{l}}
\Big( \sum_{j=k}^{k+l-1} \hat{Y}_{n,j}(t)
- \frac{l}{n} \sum_{j=1}^n \hat{Y}_{n,j}(t) \Big) \xi_k^{(r)} \\
&+ \sqrt{n}\Big(s - \frac{\lfloor sn \rfloor}{n} \Big)\frac{1}{\sqrt{l}}
\Big( \sum_{j=\lfloor sn \rfloor +1}^{\lfloor sn \rfloor+l} \hat{Y}_{n,j}(t)
- \frac{l}{n} \sum_{j=1}^n \hat{Y}_{n,j}(t) \Big) \xi_{\lfloor sn \rfloor +1}^{(r)},
\end{split}
\end{align}
where $l\in\mathbb{N}$ is a bandwidth parameter satisfying $l/n\to 0$ as $l,n\to\infty$ and
\[
\hat{Y}_{n,j} = X_{n,j} - (\hat{\mu}_2 - \hat{\mu}_1)
	\mathds{1}\{j > \lfloor \hat{s}n \rfloor \}
\]
for $j=1,\dots,n$ ($n\in\mathbb{N}$). Note that
 it is implicitly assumed that $\hat{B}_n^{(r)}((n-l+1)/n,t) = \hat{B}_n^{(r)}(s,t)$ for any $t\in[0,1]$ and any
$s\in[0,1]$ such that $\lfloor sn \rfloor > n-l+1$.
Next,
define
$$\hat{\mathbb{W}}_n^{(r)}(s,t) = \hat{B}_n^{(r)}(s,t)-s\hat{B}_n^{(r)}(1,t) ~~; ~r = 1,\dots,R.
$$
\begin{theorem} \label{bTheorem}
Let  $\hat{B}_n^{(1)},\ldots,\hat{B}_n^{(R)}$ denote  the  bootstrap processes defined by \eqref{bProcess}, where
$l = n^\beta$  for some $\beta \in (0,1/3)$. Further assume that the underlying  array $(X_{n,j}\colon j=1,\dots,n;  n\in\mathbb{N})$ satisfies
Assumption \ref{as:ts} with the additional requirement that $\nu\geq 2$ in (A1), and  suppose that,
%for any $\varepsilon>0$, there exists
%an $h>0$ such that for $r=1,\dots,R$
%\begin{align} \label{bWIP}
%\limsup_{n\to\infty} \max_{1\leq k\leq hn} \mathbb{P}
%\Big( \Big\| \frac{1}{\sqrt{n}} \sum_{i=1}^{k}\Big(\frac{1}{\sqrt{l}} \sum_{j=i}^{i+l-1} (X_{n,j} - \mu^{(j)})\Big)
%	\xi_i^{(r)} \Big\| >\varepsilon\Big) <1~.
%\end{align}
for any sequence \((r_n\colon n\in\mathbb{N})\subset \mathbb{N}\) such that \(r_n\leq n\),
\({r_n}/{n} \to 0\) as \(n\to\infty\), it follows that
\begin{align} \label{bWIP}
\frac{1}{\sqrt{n}} \sum_{i=1}^{r_n}\Big(\frac{1}{\sqrt{l}}
\sum_{j=i}^{i+l-1} (X_{n,j} - \mu^{(j)})\Big) \xi_i^{(r)} = o_\mP(1).
\end{align}
Then,
%\begin{align} \label{jointConvergence0}
$
(\hat{\mathbb{W}}_n,\hat{\mathbb{W}}_n^{(1)},\dots,\hat{\mathbb{W}}_n^{(R)})
\rightsquigarrow (\mathbb{W}, \mathbb{W}^{(1)},\dots,\mathbb{W}^{(R)})
$
%\end{align}
in $C([0,1]^2)^{R+1}$, where $\hat{\mathbb{W}}_n$  and $\mathbb{W}$ are  defined in \eqref{313} and  \eqref{316}, respectively,
and $\mathbb{W}^{(1)},\dots,\mathbb{W}^{(R)}$ are independent copies of  $\mathbb{W}$.
\end{theorem}

\medskip

The proof of Theorem \ref{bTheorem} is provided in Section \ref{sec:proofs:change-point} of the online supplement.
 Note that condition \eqref{bWIP} is similar to Assumption (A5) and
ensures that the weak invariance principle holds for the bootstrap processes.

We now consider a resampling procedure for the classical hypotheses, that is $\Delta=0$ in \eqref{37}. For that purpose, define, for \(r=1,\dots,R\),
\begin{align} \label{bstatClassic}
\check T_n^{(r) }
= \max \Big \{ \big | \hat{\mathbb{W}}_n^{(r)}(s,t) \big |~ \colon ~s,t \in [0,1] \Big \}.
\end{align}
Then, by the continuous mapping theorem,
\begin{align*}
(\sqrt{n} ~ \hat{\mathbb{M}}_n  ,~ \check T_n^{(1)},\dots,\check T_n^{(R)})
\Rightarrow ( \check T,~  \check T^{(1)},\dots, \check T^{(R)})
\end{align*}
in $\mathbb{R}^{R+1}$, where  $\check T^{(1)},\dots, \check T^{(R)}$ are independent copies  of  the random variable $\check T$ defined in \eqref{42}.
If $\check T_n^{\{\lfloor R(1-\alpha)\rfloor\}}$ is the empirical $(1-\alpha)$-quantile of the bootstrap sample
$\check T_n^{(1)},\check T_n^{(2)} , \dots ,\check  T_n^{(R)}$, the classical null hypothesis  \(H_0: \mu_1 =  \mu_2 \) of no change point
is rejected, whenever  \begin{align} \label{cptestclass}
\hat{\mathbb{M}}_n > \frac{\check T_n^{\{\lfloor R(1-\alpha) \rfloor\}}}{\sqrt{n}} ~.
\end{align}
It follows by similar arguments as given in Section \ref{subsec:proof:two-sample} of the online supplement
that this test is consistent and has asymptotic level $\alpha$ in the sense of Theorem \ref{bootstrap:classic},
that is
\begin{align*}
\lim_{R\to\infty} \limsup_{n\to\infty} \mathbb{P}_{H_0}\bigg(\hat{\mathbb{M}}_n > \frac{\check T_n^{\{\lfloor R(1-\alpha) \rfloor\}}}{\sqrt{n}}  \bigg) = \alpha~,~~~
 \liminf_{n\to\infty} \mathbb{P}_{H_1}\bigg( \hat{\mathbb{M}}_n > \frac{\check T_n^{\{\lfloor R(1-\alpha) \rfloor\}}}{\sqrt{n}} \bigg)
	=1,
\end{align*}
for any \(R \in \mathbb{N}\). The details are omitted for the sake of brevity.

We now continue developing bootstrap methodology for the problem of testing for a relevant change point, that is $\Delta >0$ in \eqref{42}.
It turns out that the theoretical analysis  is substantially more complicated as the null hypothesis defines a set in
in $C([0,1])$.
Similar as in  \eqref{estimatedSets} the estimates of the extremal sets $\mathcal{E}^+$ and
$\mathcal{E}^-$ are defined by
\begin{align} \label{estimatedSetsCP}
\hat{\mathcal{E}}_{n}^\pm  &= \Big\{ t\in[0,1] \colon \pm ( \hat{\mu}_1(t)-\hat{\mu}_2(t) )
 \geq \hat{d}_\infty - \frac{c_{n}}{\sqrt{n}} \ \Big\}
 %\nonumber \\
%\hat{\mathcal{E}}_{n}^- &= \Big\{ t\in[0,1] \ \Big| \ \hat{\mu}_1(t)-\hat{\mu}_2(t)
%\leq -\hat{d}_\infty + \frac{c_{n}}{\sqrt{n}} \ \Big\}
~,
\end{align}
where $c_n \sim \log(n)$ and $\hat d_\infty$ is given in  \eqref{statistic}.
Consider bootstrap analogs
\begin{align} \label{bootana}
T^{(r)}_{n} =  \frac{1}{\hat{s}(1-\hat{s})}\max\big\{ \max_{t\in \hat{\mathcal{E}}_{n}^+}
\hat{W}_{n}^{(r)}(\hat{s},t), \ \max_{t\in \hat{\mathcal{E}}_{n}^-}
\big(- \hat{W}_{n}^{(r)}(\hat{s},t) \big) \big\} ,
\qquad r=1,\dots,R,
\end{align}
 of the statistic $ \sqrt{n}\big( \hat{d}_\infty - d_\infty \big) $  in Corollary \ref{setT}, where
$d_\infty = \| \mu_1 - \mu_2 \|$.
\bigskip

\begin{theorem} \label{jointConvergence}
Let the assumptions of Theorem \ref{bTheorem} be satisfied,
then, if $d_\infty > 0$,
\begin{align*}
(\sqrt{n}(\hat d_\infty - d_\infty),~ T_n^{(1)},\dots,T_n^{(R)})
\Rightarrow (T(\mathcal{E)},~ T^{(1)},\dots,T^{(R)})
\end{align*}
in $\mathbb{R}^{R+1}$, where  $T^{(1)},\dots,T^{(R)}$ are independent copies  of  the random variable $T(\mathcal{E})$ defined in
Corollary \ref{setT}.
\end{theorem}
A test  for the hypothesis of a relevant change-point in time series of continuous functions is now obtained
by rejecting the null hypothesis in \eqref{37}, whenever
\begin{align} \label{bootcptest}
\hat{d}_{\infty} > \Delta + \frac{T_n^{\{\lfloor R(1-\alpha)\rfloor\}}}{\sqrt{n}},
\end{align}
where $T_n^{\{\lfloor R(1-\alpha)\rfloor\}}$ is the empirical $(1-\alpha)$-quantile of the bootstrap sample
$T_n^{(1)},T_n^{(2)} , \dots ,T_n^{(R)}$.
It follows by similar arguments as given in Section \ref{subsec:proof:two-sample} of the online supplement
that this test is consistent and has asymptotic level $\alpha$ in the sense of Theorem \ref{two-sample:bootstrap:asymptotic}, that is
\begin{align*}
\lim_{R\to\infty} \limsup_{n\to\infty} \mathbb{P}_{H_0}\bigg( \hat{d}_{\infty}
	> \Delta + \frac{T_{n}^{\{\lfloor R(1-\alpha)\rfloor\}}}{\sqrt{n}} \bigg) = \alpha
	\end{align*}
	and
	\begin{align*}
\liminf_{n\to\infty} \mathbb{P}_{H_1}\bigg( \hat{d}_{\infty}
	> \Delta + \frac{T_{n}^{\{\lfloor R(1-\alpha)\rfloor\}}}{\sqrt{n}} \bigg)
	=1,
\end{align*}
for any \(R \in \mathbb{N}\). The details are omitted for the sake of brevity.

%%%%%%%%%%%%%%%%%%%%%%%%
\section{Empirical aspects}
\label{sec5}
%%%%%%%%%%%%%%%%%%%%%%%%

In this section, the finite sample properties of the proposed methodology are investigated by means of a small simulation study (see Section \ref{sec51} and \ref{sec52}) and its applicability illustrated in a small data example (see Section \ref{sec53}). All simulation results presented here are based on $1{,}000$ runs, 
the length of the blocks  in the bootstrap procedure is  \(l=2,l_1=2,l_2=2\), and the number of bootstrap replications is chosen as $R=200$ throughout.

%%%%%%%%%%%%%%%%%%%%%%%%
\subsection{Two sample problems}
\label{sec51}
%%%%%%%%%%%%%%%%%%%%%%%%

%%%%%%%%%%%%%%%%%%%%%%%%
\subsubsection{Classical hypotheses}
\label{sec511}
%%%%%%%%%%%%%%%%%%%%%%%%

First, a brief discussion of the bootstrap test \eqref{bootclasstest} for the "classical" hypotheses  \eqref{hypclass} is given. For this problem, Horvath et al.\ \citep{HorvathKokoszkaRice2014} proposed a test in a Hilbert-space framework, and therefore a similar scenario as in this paper is considered. Specifically, the sample sizes are chosen as $m=100$, $n=200$ and the error processes are given by fAR(1) time series (see Horvath et al.\ \citep{HorvathKokoszkaRice2014}).  The left panel of Table \ref{tab1} displays the rejection probabilities of the new test \eqref{bootclasstest} for the mean functions
\begin{equation}\label{meanclass1}
\mu_1 \equiv 0 ~~,~~~\mu_2(t) = a t (1-t)
\end{equation}
for various values of the parameter $a$, while the right panel shows results for the functions
\begin{equation}\label{meanclass2}
\mu_1 \equiv 0 ~~,~~~\mu_2(t) = 0.1 \dfrac{(1-t(1-t))^k}{\int_0^1 (1-t(1-t))^k dt}
\end{equation}
for  different values of $k$. Note that only the model in \eqref{meanclass1} with $a=0$ corresponds to the null hypothesis. A similar approximation of the nominal level as for the test of Horvath et al.\ \citep{HorvathKokoszkaRice2014} is observed as well as reasonable rejection probabilities under the alternative. For the sake of a comparison the results of the test proposed by Horvath et al.\ \citep{HorvathKokoszkaRice2014} are also displayed, using their statistics $U^{(1)}_{100,200}$ and $U^{(2)}_{100,200}$ on page $109$ of their paper (these are the two numbers in brackets). For the models \eqref{meanclass1} the new test is in most cases more powerful than the test proposed by these authors. This superiority is also observed for the models \eqref{meanclass2} if $k=4,5$.  On the other hand, if $k=2,3$, the test of Horvath et al.\ \citep{HorvathKokoszkaRice2014} based  on the statistic $U^{(2)}_{100,200}$ yields the best performance, but the new test is always more powerful than the test based on their statistic  $U^{(1)}_{100,200}$.

\begin{table}[t]
\begin{center}

\begin{tabular}{c|ccc||c|ccc}
                     & \multicolumn{3}{|c||} {\eqref{meanclass1}} &                    & \multicolumn{3}{|c} {\eqref{meanclass2}} \\
\hline
a                    & 1\%          & 5\%          & 10\%         & k                  & 1\%          & 5\%          & 10\%         \\
\hline \hline
\multirow{2}{*}{0}   & 2.7          & 7.4          & 13.7         & \multirow{2}{*}{2} & 27.8         & 56.4         & 79           \\
                     & (1.8, 1.9)   & (6.6, 7.2)   & (12.2, 13.5) &                    & (27.6, 57.8) & (52.6, 77.4) & (63.6, 84.9) \\[4pt]
\multirow{2}{*}{0.4} & 21           & 37.7         & 46.7         & \multirow{2}{*}{3} & 31.3         & 76.4         & 94.9         \\
                     & (19.4, 12.3) & (35.9, 26.5) & (46.7, 36.3) &                    & (27.5, 64.3) & (49.4, 82.1) & (61.2, 88.7) \\[4pt]
\multirow{2}{*}{0.6} & 49.4         & 67.6         & 76.9         & \multirow{2}{*}{4} & 61.2         & 96.8         & 1            \\
                     & (42.1, 29.6) & (62.2, 51.8) & (73.1, 62.5) &                    & (28.2, 71.6) & (52.5, 88.7) & (66.7, 93.8) \\[4pt]
\multirow{2}{*}{0.8} & 74.3         & 87.1         & 91           & \multirow{2}{*}{5} & 90.3         & 1            & 1            \\
                     & (68.6, 53.8) & (85.7, 74.6) & (91.5, 83.1) &                    & (27.8, 78)   & (51.6, 91.7) & (64.3, 95.5)
\end{tabular}

\caption{\it Simulated rejection probabilities of the bootstrap test \eqref{bootclasstest} for the hypotheses   \eqref{hypclass} \label{tab1} (in percent).
The mean functions are given by \eqref{meanclass1} (left part) and by  \eqref{meanclass2} (right part), the sample sizes are $m=100$ and $n=200$ and  the case $a=0$ corresponds
to the null hypotheses.  The numbers in brackets represent the results of the two tests proposed by
 Horvath et al.\  \citep{HorvathKokoszkaRice2014} taken from Table 1 in this reference.}
  \end{center}
\end{table}

%%%%%%%%%%%%%%%%%%%%%%%%
\subsubsection{Confidence bands} \label{sec512}
%%%%%%%%%%%%%%%%%%%%%%%%

In order to investigate the finite sample  properties of the confidence bands proposed in Section \ref{sec:conf_bands}
we investigate  a similar  scenario
 as in Sections 6.3 and 6.4 of Aue et al.\ \cite{aueDubartNorinhoHormann2015}. To be precise let  \(D \in \mathbb{N} \), consider  $B$-spline basis functions
  \(\nu_1,\dots , \nu_D\) (here \(D=21\)) and  the linear space $\mathbb{H}=\mathrm{span}\{ \nu_1,\dots , \nu_D \} $.
Now define  independent  processes  \(\varepsilon_1, \dots, \varepsilon_n \in \mathbb{H} \subset  C([0,1])\) by
\[
  \varepsilon_j = \sum_{i=1}^D N_{i,j} \nu_i,
  \qquad j=1,\dots , n,
\]
 where \(N_{1,j},N_{2,j}, \dots, N_{D,j} \) are independent, normally distributed random
variables with expectation zero and  variance Var$(N_{i,j})=\sigma_i^2=1/i^2$  ($i=1,\ldots , D$; $j=1,\ldots , n$).
%where in each simulation run the numbers $\sigma^2_1, \ldots , \sigma^2_D$ are obtained from a permutation of  $1, 1/2^2, 1/3^2, \ldots , 1/D^2$.
The fMA(1) process is finally given by $ \eta_i = \varepsilon_i + \Theta \varepsilon_{i-1}$, where
the operator $\Theta\colon \mathbb{H} \to \mathbb{H} $ (acting on a finite dimensional space) is defined by
 \(\kappa \Psi\) (here \(\kappa = 0.5\)). The matrix \(\Psi\) is chosen randomly,  
% \(\Psi =  ( \Psi_{ij}  )_{i,i=1}^D \in \mathbb{R}^{D\times D} \), and the entries $ \Psi_{ij}$
that is, a matrix consisting of normally distributed entries with mean zero and standard deviation \( \sigma_i \sigma_j \) is generated 
and then scaled such that the resulting matrix \(\Psi\) has induced norm equal to \(1\).
Finally the
two samples are given by
$$
X_i  =  \mu_1 +  \eta_i^X~~(i=1,\ldots ,m )
\qquad\mbox{and}\qquad
Y_i  =  \mu_2 +  \eta_i^Y~~(i=1,\ldots ,n ),
$$
where $(\eta_i^X\colon i\in \Z)$ and $(\eta_i^Y\colon i\in \Z)$ are independent fMA(1) processes distributed as $(\eta_i\colon i\in \Z)$.
 In Table \ref{tab3} we display
 the simulated
 coverage percentage and the half  width of the confidence band defined in
 Theorem \ref{th:cb-2}, that is
 $$
 \frac{T_{m,n}^{\{\lfloor R(1-\alpha) \rfloor \}}}{\sqrt{m+n}} .
 $$
The two mean functions are given by
\begin{align}  \label{relex1}
\mu_1  (t) = 0,\qquad\mu_2(t) =
\begin{cases}
0.5 t,      & t\in[0,\frac{1}{5}] \\
0.1 ,          & t\in(\frac{1}{5}, \frac{3}{10}] \\
-0.5t+0.25, & t\in(\frac{3}{10}, \frac{7}{10}] \\
-0.1,          & t\in(\frac{7}{10}, \frac{4}{5}] \\
0.5t-0.5    & t\in(\frac{4}{5}, 1]
\end{cases}
\end{align}
(left panel)  and
\begin{align}  \label{relex2}
\mu_1  (t) = 0,\qquad\mu_2(t) =
\begin{cases}
0.4 t,    & t\in[0,\frac{1}{4}] \\
0.1,        & t\in(\frac{1}{4}, \frac{3}{4}] \\
-0.4t+0.4, & t\in(\frac{3}{4}, 1]
\end{cases}
\end{align} (right panel). Note that $d_\infty = 0.1$ in both cases. For example,
for sample sizes $m=50$ and $n=100$  the coverage probability of the $95\%$ uniform confidence band for the difference of the mean
functions in model \eqref{relex2} is $94.1,\%$ and the width is $2\cdot 0.34 = 0.68$. We observe a reasonable approximation of the nominal
level in all cases under consideration.

\begin{table}[h]
{  \centering
\begin{tabular}{c|ccc||ccc}
               & \multicolumn{3}{|c||} {\eqref{relex1}} & \multicolumn{3}{|c} {\eqref{relex2}} \\
\hline
\((m, n)\)     & 1\% & 5\% & 10\% & 1\% & 5\% & 10\% \\
\hline \hline
\((50, 100)\)  & (97.5, 0.44) & (92.9, 0.34) & (88, 0.29)   & (98.2, 0.44) & (94.1, 0.34) & (88.1, 0.29) \\
\((100, 100)\) & (98.3, 0.36) & (94.7, 0.28) & (89.3, 0.24) & (98.9, 0.36) & (95.5, 0.28) & (91.2, 0.24)  \\
\((100, 200)\) & (98.2, 0.31) & (94.5, 0.24) & (90.4, 0.21) & (98.5, 0.31) & (94.2, 0.24) & (89.7, 0.21)
\end{tabular}
\caption{\it\it Simulated coverage probabilities (first number) and half width (second number) of the confidence band for the difference of  the two  mean functions.
The error processes are given by  \(fMA(1)\) processes. Left part: model \eqref{relex1}; right part: model \eqref{relex2}. }
\label{tab3}}
\end{table}

%%%%%%%%%%%%%%%%%%%%%%%%
\subsubsection{Testing for a non relevant difference} \label{sec513}
%%%%%%%%%%%%%%%%%%%%%%%%

In this paragraph the finite sample properties of the
test \eqref{boottestrel} for the relevant hypotheses of the form  \eqref{H0} are investigated, using
 the same scenario as in Section \ref{sec512}, that is $fMA(1)$ time series with mean functions defined by
\eqref{relex1}  and \eqref{relex2} are used.  The  constant \(c_{m,n}\)  in \eqref{estimatedSets} 
is chosen as \(0.1\log(m+n)\).
\begin{figure}[t]
{  \centering
\includegraphics[scale=0.35]{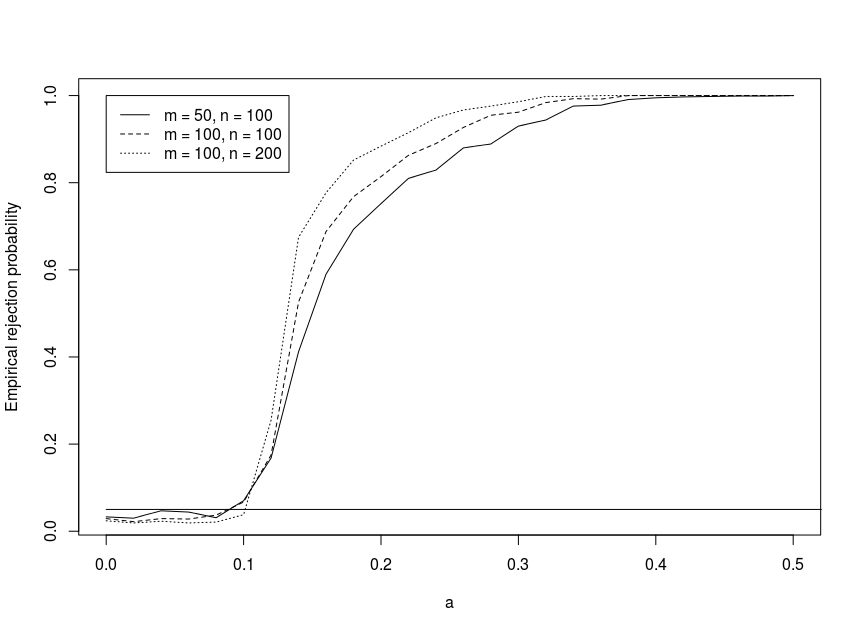}
~~
 \includegraphics[scale=0.35]{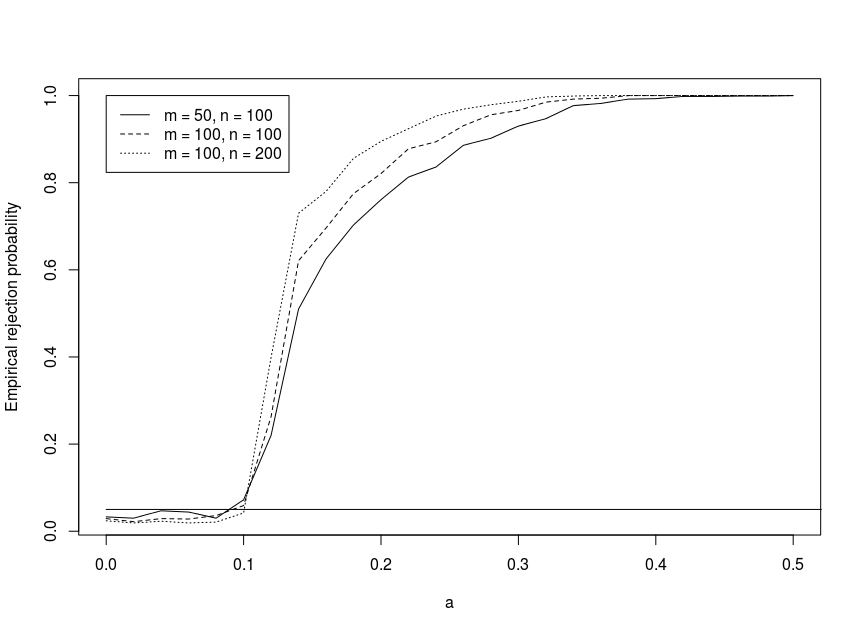}
\caption{\it Simulated rejection probabilities of the test \eqref{boottestrel} for a relevant difference in the maximal distance between mean functions  of size $\Delta =0.1$.
The error processes are given by fMA(1) processes.  Left panel: model \eqref{relex1}; right panel: model \eqref{relex2}.
}
\label{fig1}}
\end{figure}

Figure \ref{fig1} shows the rejection probabilities of the test \eqref{boottestrel} for a relevant difference in the maximal distance between mean functions of size $\Delta=0.1$, where the mean functions are given by \eqref{relex1}
(left panel) and \eqref{relex2} (right panel).
The results correspond to the theoretical properties described in Theorem \ref{two-sample:bootstrap:asymptotic}.
More precisely, if $d_\infty =|a| < \Delta$ (interior of the null hypothesis) the probability of rejection is substantially smaller than $\alpha$
and decreases with increasing sample size.
 At the boundary of the  hypotheses, that is $d_\infty= \Delta$, we observe from  Table \ref{tab2}
 that the simulated level is close to the nominal level $\alpha$. 
%This approximation is slightly more precise for the mean functions \eqref{relex2}.
%, where the set of extremal points of the difference $\mu_1 - \mu_2$ is the complete interval $[0,1]$.
 On the other hand, the power of the test is strictly increasing with $d_\infty > \Delta$
 (see Figure \ref{fig1}).

\begin{table}[h]
{  \centering
\begin{tabular}{c|ccc||ccc}
               & \multicolumn{3}{|c||} {\eqref{relex1}} & \multicolumn{3}{|c} {\eqref{relex2}} \\
\hline
\((n, m)\)     & 1\% & 5\% & 10\% & 1\% & 5\% & 10\% \\
\hline \hline
\((50, 100)\)  & 2.4 & 7   & 13.7 & 2.5 & 7.2 & 13.8 \\
\((100, 100)\) & 1.7 & 6.7 & 11.1 & 1.7 & 5.9 & 11.5 \\
\((100, 200)\) & 1.2 & 3.8 & 9.7  & 1.2 & 4.2 & 10.2
\end{tabular}
\caption{\it\it Simulated nominal level of the test \eqref{boottestrel} for a relevant difference in the mean functions
at the boundary of the null hypothesis, that  is $d_\infty= \Delta =0.1$.
The error processes are given by  fMA(1) processes. Left part: model \eqref{relex1}; right part: model \eqref{relex2}. }
\label{tab2}}
\end{table}

%%%%%%%%%%%%%%%%%%%%%%%%%
\subsection{Change-point inference}  \label{sec52}
%%%%%%%%%%%%%%%%%%%%%%%%%

In this section we investigate the finite sample performance of
 the change-point tests based on the maximal deviation $d_\infty$ proposed in Section \ref{sec:change-point}.
Throughout this section consider a time series of the form
\begin{equation} \label{cpmodel}
X_i  =  \mu_i +  \eta_i^X~~(i=1,\ldots ,n )
\end{equation}
where $( \eta_i^X) _{i \in \mathbb{Z}} $ is the  fMA(1) process defined in Section \ref{sec512} and the sequence of mean functions satisfies \eqref{36a}.
We investigate the problem of testing for a relevant  change of size $\Delta=0.4$ in the mean functions. The corresponding results are depicted in Table \ref{tab4},
where we consider again the model  \eqref{cpmodel} with  mean functions before and after the change point $s^*=0.5$ given by  \eqref{relex1} and \eqref{relex2}. The parameter \(c_n\)  in \eqref{estimatedSetsCP} is chosen as \(0.1 \log(n)\).
The case $a=0.4$ corresponds to the boundary of the hypotheses and here  a rather accurate approximation of the nominal level
is observed. In the interior of the
null hypothesis
(that is $a <0.4$) the rejection probability, for \(n=200, 500\), is strictly smaller than the nominal level and decreasing with an increasing sample size as described
 at the end of Section \ref{testdiffrel} (note that the same arguments also hold for the change point problem). Similarly, under the alternative (i.e. $a>0.4$) the test
shows reasonable rejection probabilities which are increasing with  sample size and $a$.

\begin{table}[h]
{  \centering
\begin{tabular}{c|c||ccc||ccc||ccc }
 & n & \multicolumn{3}{|c||} {100} & \multicolumn{3}{c||}{200} &\multicolumn{3}{|c} {500} \\
\hline
 & $a$         & 1\% & 5\% & 10\% & 1\% & 5\% & 10\%& 1\% & 5\% & 10\%  \\
\hline \hline
\multirow{7}{*}{\eqref{relex1}} 
 & 0.37 & 1.9  & 4.7  & 8.4  & 0.3  & 0.5  & 1.1  & 0    &  0   &  0.1 \\
 & 0.38 & 2.2  & 5.2  & 7.2  & 0.2  & 0.6  & 1.2  & 0    &  0   &  0.1 \\
 & 0.39 & 1.8  & 5    & 9    & 0.4  & 1.1  & 3.4  & 0.2  &  0.5 &  1.2 \\
 & 0.4  & 2.8  & 9.3  & 17.7 & 1.3  & 5.1  & 10.4 & 0.9  &  4.2 &  8.6 \\
 & 0.41 & 6.1  & 14.0 & 24.7 & 5    & 15   & 26.2 & 12.1 & 29.8 & 44.4 \\
 & 0.42 & 09.3 & 25.3 & 40.4 & 14.8 & 36.9 & 58.1 & 45.9 & 84.2 & 95.4 \\
 & 0.43 & 19.3 & 42.2 & 62.9 & 42.2 & 76.4 & 91.0 & 91.8 & 99.4 & 99.7 \\
\hline \hline
\multirow{7}{*}{\eqref{relex2}}
 & 0.37 & 1.9  & 4.6  & 8.2  & 0.3  & 0.5  & 1.1  & 0    & 0    & 0    \\
 & 0.38 & 2.1  & 4.6  & 7.2  & 0.1  & 0.6  & 1.2  & 0    & 0    & 0.1  \\
 & 0.39 & 2.0  & 5.2  & 8.7  & 0.3  & 1.1  & 3.1  & 0.1  & 0.2  & 0.8  \\
 & 0.4  & 2.3  & 7.8  & 16.3 & 1.5  & 5.4  & 11.6 & 0.7  & 4.2  & 9.7  \\
 & 0.41 & 6.7  & 17.4 & 32.6 & 7.9  & 21.3 & 37.3 & 18.0 & 43.8 & 64.9 \\
 & 0.42 & 14.6 & 35.8 & 54.9 & 27.7 & 62.1 & 81.9 & 76.1 & 96.0 & 99.5 \\
 & 0.43 & 32.7 & 63.9 & 78.3 & 68.1 & 91.8 & 96.5 & 98.1 & 99.7 & 99.8 \\
\end{tabular}
\caption{\it\it Simulated nominal level (in percent) of the test \eqref{bootcptest} for  the  hypotheses  \eqref{37} of a  relevant change  in
the maximal deviation of the mean functions
in model \eqref{cpmodel}, where $\Delta=0.4$.
The case $a=0.4$ corresponds  to the boundary of the hypotheses, $a < 0.4 $ to the null hypothesis and $a>0.4$ to the alternative.
The error processes are given by fMA(1) processes. Upper part: model \eqref{relex1}; lower part: model \eqref{relex2}.
 }
\label{tab4}}
\end{table}

%%%%%%%%%%%%%%%%%%%%%%%%%
\subsection{Data example}  \label{sec53}
%%%%%%%%%%%%%%%%%%%%%%%%%

To illustrate the proposed methodology, two applications to annual temperature profiles are reported in this section. Data of this kind were recently used in Aue and van Delft \cite{aueVandelft} and van Delft et al.~\cite{vanDelftEtAl} in the context of stationarity tests for functional time series and earlier in Fremdt et al.~\cite{fremdtEtAl} in support of methodology designed to choose the dimension of the projection space obtained with fPCA. For all examples, functions were generated from daily values through representation in a Fourier basis consisting of 49 basis functions, where reasonable deviations from this preset do not qualitatively change the outcome of the analyses to follow.

%%%%%%%%%%%%%%%%%%%%%%%%%
\subsubsection{Two-sample tests}  \label{sec531}
%%%%%%%%%%%%%%%%%%%%%%%%%

For the two-sample testing problem, annual temperature profiles were obtained from daily temperatures recorded at measuring stations in Cape Otway, a location close to the southernmost point of Australia, and Sydney, a city on the eastern coast of Australia. This led to $m=147$ respectively $n=153$ functions for the two samples. Differences in the temperature profiles are expected due to different climate conditions, so the focus of the relevant tests is on working out how big the discrepancy might be. The data considered here is part of the larger data set considered, for example, in Aue and van Delft \cite{aueVandelft} and van Delft et al.~\cite{vanDelftEtAl}.

To set up the test for the hypotheses \eqref{H0}, the statistic in \eqref{2statistic} was computed, resulting in the value
\[
\hat d_\infty=5.73.
\]
To see whether this is significant, the proposed bootstrap methodology was applied. To estimate the extremal sets in \eqref{35}, the estimators in \eqref{estimatedSets} were utilized with $c_{m,n}=0.1\log(m+n)=0.570$. The resulting bootstrap quantiles are reported in the second row of Table \ref{tab:ts-temp}. Also reported in this table are the results of the bootstrap procedure in \eqref{boottestrel} for various levels $\alpha$ and relevance $\Delta$. Note that the maximum difference in mean the functions is achieved at $t=0.99$, towards the end of December and consequently during the Australian summer. The results show that there is strong evidence in the data to support the hypothesis that the maximal difference is at least $\Delta=5.4$, but that there is no evidence that the maximal difference is even larger than $\Delta=5.6$. Several intermediate values of $\Delta$ led to weaker support of the alternative. The left panel of Figure \ref{fig:cp-before-after} displays the difference in mean functions graphically.

\begin{table}[t]
\vspace{.5cm}
\begin{center}
\begin{tabular}{c@{\qquad}llll}
\hline
$\Delta$ & 99\% & 97.5\% & 95\% & 90\%  \\
\hline
$q$  & 5.138 & 4.201 & 3.757 & 3.009  \\
\hline
5.4  & TRUE  & TRUE  & TRUE  & TRUE    \\
5.45 & FALSE & TRUE  & TRUE  & TRUE    \\
5.5  & FALSE & FALSE & TRUE  & TRUE    \\
5.55 & FALSE & FALSE & TRUE  & TRUE   \\
5.6  & FALSE & FALSE & FALSE & FALSE   \\
\hline
\end{tabular}
\end{center}
\caption{\it Summary of the bootstrap two sample procedure for relevant hypotheses with varying $\Delta$ for the annual temperature curves. The label TRUE refers to a rejection of the null, the label FALSE to a failure to reject the null.}
\label{tab:ts-temp}
\end{table}

%%%%%%%%%%%%%%%%%%%%%%%%%
\subsubsection{Change-point tests}  \label{sec532}
%%%%%%%%%%%%%%%%%%%%%%%%%

Following Fremdt et al.~\cite{fremdtEtAl}, annual temperature curves were obtained from daily minimum temperatures recorded in Melbourne, Australia. This led to 156 annual temperature profiles ranging from 1856 to 2011 to which the change-point test for the relevant hypotheses in \eqref{37} was applied based on the rejection decision in \eqref{41}. To compute the test statistic $\hat d_\infty$ in \eqref{statistic}, note that the estimated change-point in \eqref{cpEstimator} was $\hat s=0.62$ (corresponding to the year 1962). This gives
\[
\hat d_\infty = 1.765.
\]
To see whether this value leads to a rejection of the null, the multiplier bootstrap procedure was utilized with bandwidth parameter $l=1$, leading to the rejection rule in \eqref{bootcptest}. In order to apply this procedure, first the extremal sets $\hat{\mathcal{E}}^+$ and $\hat{\mathcal{E}}^-$ in \eqref{estimatedSetsCP} were selected, choosing $c_n=0.1\log n=0.504$. This yielded the bootstrap quantiles reported in the second row of Table \ref{tab:cp-temp}.

\begin{figure}[t]
{  \centering
\includegraphics[scale=0.35]{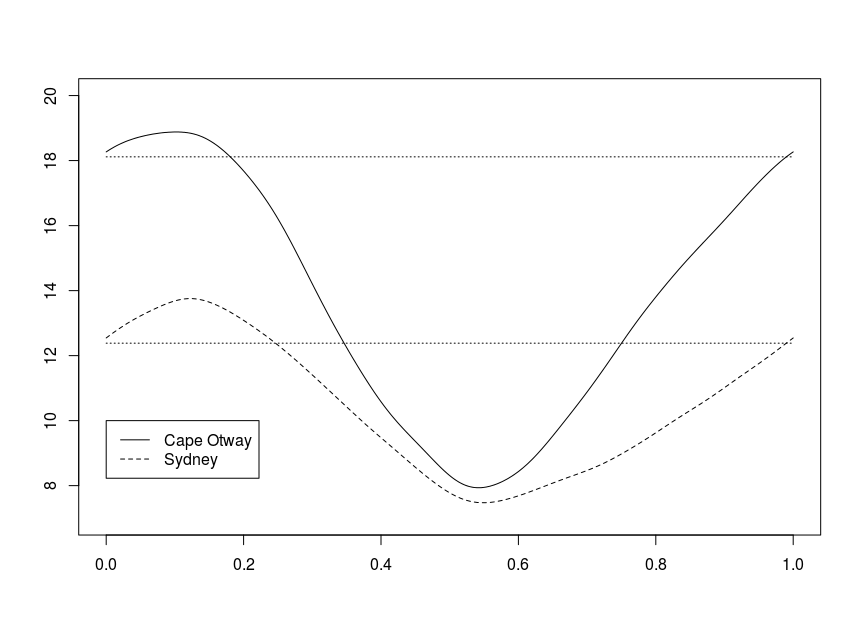}
~~
 \includegraphics[scale=0.35]{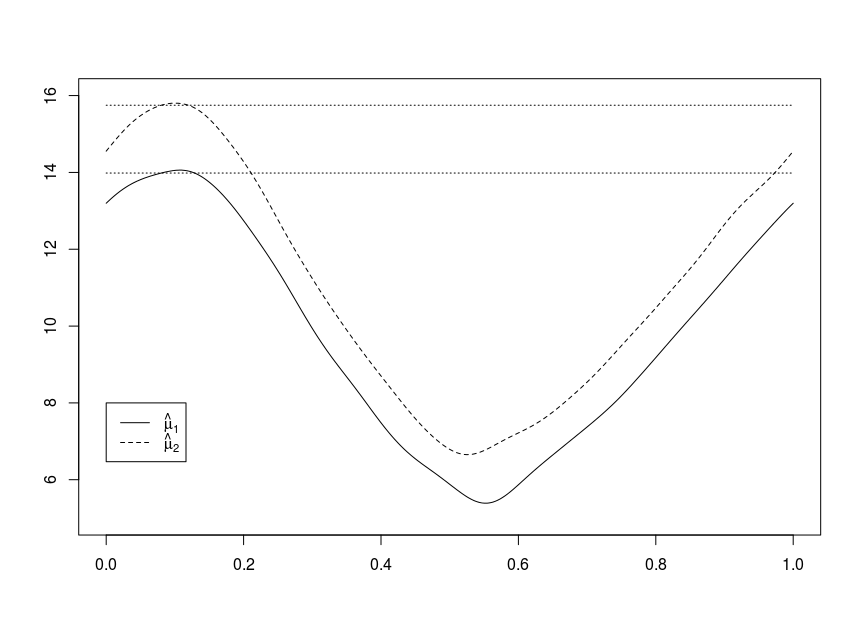}
\caption{\it Mean functions for the Australian temperature data. Left panel: Estimated mean functions of the Cape Otway and Sydney series for the two-sample case. Right panel: Estimated mean functions before and after the estimated change-point for the Melbourne temperature series.
}
\label{fig:cp-before-after}}
\end{figure}

Several values for $\Delta$, determining which deviations are to be considered relevant, were then examined. The results of the bootstrap testing procedure are summarized in Table \ref{tab:cp-temp}. It can be seen that the null hypothesis of no relevant change was rejected at all considered levels for the smaller choice $\Delta=1.2$. On the other extreme, for $\Delta=1.4$, the test never rejected. For the intermediate values $\Delta=1.25,1.3,1.35$, the null was rejected at the 2.5\%, 5\% and 20\% level, at the 5\% and 10\% level, and at the 10\% level, respectively. Estimating the mean functions before and after $\hat s$ (1962) shows that the maximum difference of the mean functions is approximately $1.765$, lending further credibility to the conducted analyses. The right panel of Figure \ref{fig:cp-before-after} displays both mean functions for illustration. It can be seen that the mean difference is maximal during the Australian summer (in February), indicating that the mean functions of minimum temperature profiles have been most drastically changed during the hottest part of the year. The results here are in agreement with the findings put forward in Hughes et al.~\cite{hughesEtAl}, who reported that average temperatures in Antarctica have risen due to increases in minimum temperatures.

In summary, the results in this section highlight that there is strong evidence in the data for an increase in the mean function of Melbourne annual temperature profiles, with the maximum difference between ``before'' and ``after'' mean functions being at least 1.25 degrees centigrade. There is weak evidence that this difference is at least 1.35 degrees centigrade, but there is no support for the relevant hypothesis that it is even larger than that.

\begin{table}[t]
\vspace{.5cm}
\begin{center}
\begin{tabular}{c@{\qquad}llll}
\hline
$\Delta$ & 99\%          & 97.5\%          & 95\%         & 90\%  \\
\hline
$q$ & 6.632 & 6.278 & 5.603 & 4.697 \\
\hline
1.2  & TRUE  & TRUE  & TRUE  & TRUE  \\
1.25 & FALSE & TRUE  & TRUE  & TRUE  \\
1.3  & FALSE & FALSE & TRUE  & TRUE  \\
1.35 & FALSE & FALSE & FALSE & TRUE  \\
1.4  & FALSE & FALSE & FALSE & FALSE \\
\hline
\end{tabular}
\end{center}
\caption{\it Summary of the bootstrap change-point procedure for relevant hypotheses with varying $\Delta$ for the annual temperature curves. The label TRUE refers to a rejection of the null, the label FALSE to a failure to reject the null.}
\label{tab:cp-temp}
\end{table}

\bigskip
\medskip

{\bf Acknowledgements }
The authors thank Martina
 Stein, who typed  parts  of this manuscript with considerable technical expertise, Stanislav Volgushev
 for helpful discussions about the proof of Theorem \ref{setConvergence} and Josua G\"osmann and Daniel Mei{\ss}ner for helpful discussions about the proof of Lemma \ref{maxSetConvergence}.
 %The authors are also grateful to three referees and the associate editor. Their constructive comments on an earlier version yielded to substantial improvement of an earlier version of this paper.

\bibliography{adk-banach-13112017}

\newpage
%%%%%%%%%%%%%%%%%%%%%%%%%
\section{Online supplement: Proofs}
\label{sec:proofs}
%%%%%%%%%%%%%%%%%%%%%%%%%

%%%%%%%%%%%%%%%%%%%%%%%%
\subsection{Convergence of suprema of non-centered processes }
\label{subsec:banach_methods:main}
%%%%%%%%%%%%%%%%%%%%%%%%

This section collects some preliminary results useful in the proofs of the main results. Van der Vaart and Wellner \cite{wellner1996} formulate the concept of weak convergence in a slightly more general way, not restricting to sequences of random variables. They state their theory for \textit{nets} of random variables $(X_\alpha\colon \alpha \in A)$, where $A$ is a \textit{directed set}, that is, a non-empty set equipped with a partial order $\leq$ with the additional property that, for any pair $a,b\in A$, there is $c\in A$ such that $a\leq c$ and $b\leq c$. This is a natural extension of the case considered in this paper and it will sometimes be used in the proofs to follow.
The weak convergence result presented next is central to establishing limit results for the two-sample and change-point tests.

\begin{theorem}
\label{thm0}
Let $(X_{\alpha}\colon\alpha\in A)$ denote a net of random variables taking values in $C(T)$ and let $\mu \in C(T)$. If $r\colon A \to \mathbb{R}_+$ is defined such that $a_{\alpha} = log(r_{\alpha})/ \sqrt{r_{\alpha}} = o(1)$ and if $Z_\alpha=\sqrt{r_{\alpha}}(X_{\alpha} - \mu)$ converges weakly to a Gaussian random variable $Z$ in $C(T) $, then
\[
D_\alpha
= \sqrt{r_{\alpha}} \big(\|X_\alpha\| - \|\mu\|\big)
\Rightarrow D({\mathcal{E}})
= \max\Big\{ \sup_{t \in \mathcal{E}^+} Z(t), \sup_{t \in \mathcal{E}^-} -Z(t) \Big\}
\]
in $\mathbb{R}$, where $\Rightarrow$ denotes convergence in distribution, $\mathcal{E}=\mathcal{E}^+\cup\mathcal{E}^-$ the set of extremal points of $\mu$, divided into $\mathcal{E}^{\pm} = \{ t \in T \colon \mu(t) = \pm \|\mu\| \}$, tacitly adopting the convention $\mathcal{E}^{\pm} =T$ if $\mu \equiv 0$.
\end{theorem}

\medskip
\noindent
{\bf Proof of Theorem \ref{thm0}.}
 First notice that, if $\mu$ is the zero function, the result is implied by the continuous mapping theorem. Hence, only the case $\|\mu\|>0$ is considered in the following for which some arguments from Raghavachari \cite{raghavachari1973} are applied.   To show that $D_\alpha\stackrel{\cal D}{\longrightarrow}D(\mathcal{E})$, introduce
\begin{align*}
D_{\alpha}(\mathcal{E}) &= \sqrt{r_{\alpha}} \left(\sup_{t \in \mathcal{E}} |X_{\alpha}(t)| - \|\mu\|\right)
\end{align*}
and note that the assertion of  Theorem \ref{thm0} directly follows from the following lemmas.
\begin{lemma}
\label{lem:proofs:weak_conv:1}
Under the assumptions of Theorem \ref{thm0}, it holds that
$
D_{\alpha}(\mathcal{E})\Rightarrow D({\mathcal{E}}).
$
\end{lemma}
\begin{lemma}
\label{lem:proofs:weak_conv:2}
Under the assumptions of Theorem \ref{thm0}, it holds that  $
R_{\alpha} = D_{\alpha} - D_{\alpha}(\mathcal{E}) = o_{\mathbb{P}}(1).
$
\end{lemma}

\medskip
\noindent
{\bf Proof of Lemma \ref{lem:proofs:weak_conv:1}.}
Define the random variable
\begin{align*}
\tilde{D}_{\alpha}(\mathcal{E})
= \max \Big\{ \sup_{t \in \mathcal{E}^+} Z_{\alpha}(t), \sup_{t \in \mathcal{E}^-} -Z_{\alpha}(t) \Big\}.
\end{align*}
From the continuous mapping theorem it follows that $\tilde{D}_{\alpha}(\mathcal{E}){\Rightarrow} D({\mathcal{E}})$. Recall that $d_\infty=\|\mu\|$ and rewrite $D_{\alpha}(\mathcal{E})$ as
\begin{align*}
D_{\alpha}(\mathcal{E})
&= \sqrt{r_{\alpha}} \left(\sup_{t \in \mathcal{E}} |X_{\alpha}(t)| - d_{\infty} \right) = \sqrt{r_{\alpha}} \max \left \{ \sup_{t \in \mathcal{E}^+} |X_{\alpha}(t)|-d_{\infty}, \ \sup_{t \in \mathcal{E}^-} |X_{\alpha}(t)|-d_{\infty}\right\} \\
&=  \sqrt{r_{\alpha}} \max \left \{ \sup_{t \in \mathcal{E}^+} X_{\alpha}(t) - d_{\infty}, \ \sup_{t \in \mathcal{E}^-} -X_{\alpha}(t)-d_{\infty} \right\}
+ o_{\mathbb{P}}(1) \\
&= \sqrt{r_{\alpha}} \max \left\{ \sup_{t \in \mathcal{E}^+} (X_{\alpha}(t) - \mu(t)), \ \sup_{t \in \mathcal{E}^-} (-X_{\alpha}(t) + \mu(t)) \right\} + o_{\mathbb{P}}(1) \\
&= \tilde{D}_{\alpha}(\mathcal{E}) + o_{\mathbb{P}}(1).
\end{align*}
The assertion follows.

\medskip\medskip
\noindent
{\bf Proof of Lemma \ref{lem:proofs:weak_conv:2}.}
First observe that the assumed weak convergence of $Z_\alpha=\sqrt{r_\alpha}(X_\alpha-\mu)$ to $Z$ and the continuous mapping theorem imply that $\|Z_\alpha\|{\Rightarrow}\|Z\|$. Consequently, Slutsky's lemma yields that
\begin{align} \label{firstCondition}
\lim_{{\alpha} \to \infty}
\mathbb{P}\left(\|X_{\alpha}-\mu\|
> \frac{a_{\alpha}}{2} \right) = 0,
\end{align}
noting that $a_\alpha=\log(r_\alpha)/\sqrt{r_\alpha}=o(1)$ by assumption. The proof of the lemma is now given in two steps.
\\
{\it Step 1:}\/ First, define the sets
\begin{align*}
\mathcal{E}_{\alpha}^{\pm} &= \{ t \in T \colon |\pm d_{\infty} - \mu(t)| \leq a_{\alpha} \}  = \{ t \in T \colon \pm \mu(t) \geq d_{\infty} - a_{\alpha} \}
\end{align*}
on which the function $\pm\mu$ is within $a_\alpha$ from its extremal value $d_\infty$, and let $\mathcal{E}_{\alpha} = \mathcal{E}_{\alpha}^+ \cup \mathcal{E}_{\alpha}^-$. It will be shown in the following that $D_\alpha$ can be replaced with $D_\alpha(\mathcal{E}_\alpha)$ in the definition of $R_\alpha$ without changing its asymptotic behavior. Here, $D_\alpha(\mathcal{E}_\alpha)$ is defined as $D_\alpha(\mathcal{E})$, using $\mathcal{E}_\alpha$ in place of $\mathcal{E}$. A corresponding definition is used for $D_\alpha(T\!\setminus\!\mathcal{E}_\alpha)$. Then,
\begin{align}
\label{calEn}
0 \leq R_{\alpha}
&= D_\alpha - D_\alpha(\mathcal{E})
=\max\big\{D_\alpha(\mathcal{E}_\alpha)- D_\alpha(\mathcal{E}),D_\alpha(T\!\setminus\!\mathcal{E}_\alpha)- D_\alpha(\mathcal{E})\big\}.
\end{align}
The second term in the maximum on the right-hand side of \eqref{calEn} is negligible as the following considerations show. Note that
\begin{align*}
D_\alpha(T\!\setminus\!\mathcal{E}_\alpha)- D_\alpha(\mathcal{E})
=& \sqrt{r_{\alpha}} \Big( \sup_{t \in T\setminus \mathcal{E}_{\alpha}}
|X_{\alpha}(t) - \mu(t) +\mu(t)|
 - \sup_{t \in \mathcal{E}} |X_{\alpha}(t)| \Big)  \\
\leq & \sqrt{r_{\alpha}} \Big( \sup_{t \in T\setminus \mathcal{E}_{\alpha}}
\left( |Y_{\alpha}(t)| + |\mu(t)| \right)
 - \sup_{t \in \mathcal{E}} |X_{\alpha}(t)|\Big),
\end{align*}
where $Y_\alpha=X_\alpha-\mu$ is the centered version of $X_\alpha$. The definition of $\mathcal{E}_\alpha$ yields that
\begin{align*}
\sqrt{r_{\alpha}}
\Big( \sup_{t \in T\setminus \mathcal{E}_{\alpha}}
\left( |Y_{\alpha}(t)| + |\mu(t)| \right)
 - \sup_{t \in \mathcal{E}} |X_{\alpha}(t)| \Big )
&  <  \sqrt{r_{\alpha}} \Big ( \sup_{t \in T\setminus \mathcal{E}_{\alpha}}
|Y_{\alpha}(t)| + d_{\infty} - a_{\alpha}
 - \sup_{t \in \mathcal{E}} |X_{\alpha}(t)| \Big) \\
&\leq \sqrt{r_{\alpha}} \sup_{t \in T} |Y_{\alpha}(t)| - \log(r_{\alpha})
 - \sqrt{r_{\alpha}}  \Big( \sup_{t \in \mathcal{E}} |X_{\alpha}(t)| - d_{\infty} \Big).
\end{align*}
Observe next that
\begin{align*}
\lim_{\alpha \to \infty} &\mathbb{P} \Big(
\sqrt{r_{\alpha}} \sup_{t \in T} |Y_{\alpha}(t)| - \log(r_{\alpha})
 - \sqrt{r_{\alpha}} \Big( \sup_{t \in \mathcal{E}} |X_{\alpha}(t)| - d_{\infty} \Big ) > 0
 \Big ) \\
\leq & \lim_{\alpha \to \infty} \mathbb{P} \Big(
\Big | \sqrt{r_{\alpha}} \sup_{t \in T} |Y_{\alpha}(t)| - \sqrt{r_{\alpha}}
\Big ( \sup_{t \in \mathcal{E}} |X_{\alpha}(t)| - d_{\infty} \Big ) \Big |
> \log(r_{\alpha}) \Big) \\
\leq & \lim_{\alpha \to \infty} \left\{ \mathbb{P} \Big(
\sqrt{r_{\alpha}} \sup_{t \in T} |Y_{\alpha}(t)|
> \frac{\log(r_{\alpha})}{2} \Big)
+ \mathbb{P} \Big(
 \sqrt{r_{\alpha}} \Big | \sup_{t \in \mathcal{E}} |X_{\alpha}(t)| - d_{\infty} \Big |
> \frac{\log(r_{\alpha})}{2} \Big)  \right\}   ~=~ 0.
\end{align*}
The first term on the right-hand side converges to $0$ because of \eqref{firstCondition} and, recalling the weak convergence shown in Lemma \ref{lem:proofs:weak_conv:1}, it follows by similar arguments that the second term converges to $0$. On the other hand, the first term in \eqref{calEn} is always nonnegative as $\mathcal{E}$ is a subset of $\mathcal{E}_{\alpha}$. Therefore,
\begin{align*}
0 \leq R_{\alpha} = D_\alpha(\mathcal{E}_\alpha)-D_\alpha(\mathcal{E})+o_{\mathbb{P}}(1)
= \sqrt{r_{\alpha}}\Big( \sup_{t \in \mathcal{E}_{\alpha}} |X_{\alpha}(t)|
- \sup_{t \in \mathcal{E}} |X_{\alpha}(t)| \Big ) + o_{\mathbb{P}}(1)
\end{align*}
holds and it suffices to evaluate $D_\alpha(\mathcal{E}_\alpha)-D_\alpha(\mathcal{E})$.
\\
{\it Step 2:}\/  Define $R_\alpha^\pm=D_\alpha(\mathcal{E}_\alpha^\pm)-D_\alpha(\mathcal{E}^\pm)$  and observe that then
%\begin{align*}
$0 \leq R_{\alpha} %&= \max \Big\{\sqrt{r_{\alpha}}
%\Big ( \sup_{t \in \mathcal{E}_{\alpha}^{+}} |X_{\alpha}(t)|
%- \sup_{t \in \mathcal{E}^{}} |X_{\alpha}(t)|\Big ), \\
%& \hspace{48pt} \sqrt{r_{\alpha}}
%\Big ( \sup_{t \in \mathcal{E}_{\alpha}^{-}} |X_{\alpha}(t)|
%- \sup_{t \in \mathcal{E}^{}} |X_{\alpha}(t)| \Big ) \Big \} + o_{\mathbb{P}}(1) \\
%&
\leq \max \big \{R_{\alpha}^+, R_{\alpha}^- \big \} + o_{\mathbb{P}}(1)$.
%\end{align*}
For the completion of the proof it is necessary to show that $R_{\alpha}^{\pm} = o_{\mathbb{P}}(1)$. To this end,
write
\begin{align*}
0 \leq R_{\alpha}^+ &= \sqrt{r_{\alpha}}
\Big ( \sup_{t \in \mathcal{E}_{\alpha}^{+}} |X_{\alpha}(t)|
- \sup_{t \in \mathcal{E}^{+}} |X_{\alpha}(t)| \Big )  =\sqrt{r_{\alpha}}
\Big ( \sup_{t \in \mathcal{E}_{\alpha}^{+}} X_{\alpha}(t)
- \sup_{t \in \mathcal{E}^{+}} X_{\alpha}(t) \Big ) + o_{\mathbb{P}}(1) \\
& = \sqrt{r_{\alpha}}
\Big( \sup_{t \in \mathcal{E}_{\alpha}^{+}} \big ( X_{\alpha}(t) - d_{\infty}\big)
- \sup_{t \in \mathcal{E}^{+}} \big( X_{\alpha}(t) - d_{\infty} \big)\Big)
+ o_{\mathbb{P}}(1) \\
%& = \sqrt{r_{\alpha}}
%\Big( \sup_{t \in \mathcal{E}_{\alpha}^{+}} \big ( X_{\alpha}(t) - d_{\infty} \big  )
%- \sup_{t \in \mathcal{E}^{+}} \big  ( X_{\alpha}(t) - \mu(t) \big  ) \Big)
%+ o_{\mathbb{P}}(1) \\
& \leq \sqrt{r_{\alpha}}
\Big( \sup_{t \in \mathcal{E}_{\alpha}^{+}} \big  ( X_{\alpha}(t) - \mu(t) \big  )
- \sup_{t \in \mathcal{E}^{+}}\big  ( X_{\alpha}(t) - \mu(t) \big  ) \Big)
+ o_{\mathbb{P}}(1) \\
& = \sqrt{r_{\alpha}}
\Big( \sup_{t \in \mathcal{E}_{\alpha}^{+}} Y_{\alpha}(t)
- \sup_{t \in \mathcal{E}^{+}} Y_{\alpha}(t) \Big)
+ o_{\mathbb{P}}(1).
\end{align*}
Define
$\mathcal{E}^+(\gamma) = \{s \in T \colon
\exists \ t \in \mathcal{E}^+ \text{ with } \rho(t,s) < \gamma \}$ and
$\delta_{\alpha} = 2 \inf\{\gamma > 0 \colon
\mathcal{E}_{\alpha}^+ \subset \mathcal{E}^+ (\gamma) \}$.
Since $\mathcal{E}_{\alpha}^+ \subset \mathcal{E}^+(\delta_{\alpha})$ the above expression can be bounded by
\begin{align*}
 \sqrt{r_{\alpha}}
  \Big( \sup_{t \in \mathcal{E}_{\alpha}^{+}} Y_{\alpha}(t)
- \sup_{t \in \mathcal{E}^{+}} Y_{\alpha}(t)\Big)
+ o_{\mathbb{P}}(1) & \leq \sqrt{r_{\alpha}}
\Big ( \sup_{t \in \mathcal{E}^+(\delta_{\alpha})} Y_{\alpha}(t)
- \sup_{t \in \mathcal{E}^{+}} Y_{\alpha}(t) \Big )
+ o_{\mathbb{P}}(1) \\
&\leq \sqrt{r_{\alpha}} \sup_{\rho(s,t)< \delta_{\alpha}} |Y_{\alpha}(s) - Y_{\alpha}(t)|
+ o_{\mathbb{P}}(1).
\end{align*}
Because of equicontinuity, it remains to show that
$\lim_{\alpha \to \infty} \delta_{\alpha} = 0$. Now, the sequence $(\delta_{\alpha}\colon\alpha \in A)$ decreases in $\alpha$ and
$\delta_{\alpha} > 0$ so that $\lim_{\alpha \to \infty} \delta_{\alpha}$ exists.
By construction, $\mathcal{E}_{\alpha}^+ \subset \mathcal{E}^+(\delta_{\alpha})$ but
$\mathcal{E}_{\alpha}^+ \not \subset \mathcal{E}^+(\delta_{\alpha}/4)$. There is hence a subsequence
$(s_{\alpha}\colon\alpha \in A) \subset \mathcal{E}_{\alpha}^+$ such that
$\rho(s_{\alpha}, t) \geq \delta_{\alpha}/4$ for all $t \in \mathcal{E}^+$ and all
$\alpha \in A$. The compactness of $T$ implies that $(s_{\alpha}\colon\alpha \in A)$
contains a convergent subsequence
$(s_{\alpha_{\beta}}\colon\beta \in A)$. It follows thus that
$\lim_{\beta \to \infty} s_{\alpha_\beta} = s \in T$ and
$d_{\infty} = \lim_{\beta \to \infty} \mu(s_{\alpha_\beta}) = \mu(s)$
because $\mu$ is continuous. Consequently $s \in \mathcal{E}^+$ but on the other hand
$\rho(s_{\alpha_\beta},s) \geq \delta_{\alpha}/4$ so that
$\lim_{\alpha \to \infty}\delta_{\alpha}$ has to be $0$.
The lemma is now proven because similar arguments imply that $R_{\alpha}^- = o_{\mathbb{P}}(1)$.

%%%%%%%%%%%%%%%%%%%%%%%%
\subsection{Proofs of the results in Section \ref{sec:banach_methods}}
\label{proofssec2}
%%%%%%%%%%%%%%%%%%%%%%%%

\begin{proof}[\bf Proof of Theorem \ref{mixingCLT}]
Note that the results of Section 1.5 in Van der Vaart and Wellner \cite{wellner1996} also hold, if  the space
$\ell^\infty$ is replaced by $C(T)$. Consequently, the assertion follows from the
 convergence of the finite dimensional distributions of $(G_n\colon n\in\mathbb{N})$
 and the existence of  a metric $\rho$ such that $(G_n\colon n\in\mathbb{N})$ is uniformly $\rho$-equicontinuous in probability.

In a first step, convergence of finite-dimensional distributions is verified.
Let $q\in\mathbb{N}$ and $t_1,\dots,t_q \in T$. By the Cram{\'e}r-Wold device, the convergence of the finite-dimensional distributions
follows from
\begin{align} \label{barZdef}
%\begin{split}
Z_n = \sum_{j=1}^q c_j G_n(t_j)
&= \sum_{j=1}^q c_j  \frac{1}{\sqrt{n}} \sum_{i=1}^n \big(X_{n,i}(t_j) - \mu^{(i)}(t_j) \big)
% \\ &= \frac{1}{\sqrt{n}} \sum_{i=1}^n \sum_{j=1}^q c_j \big(X_{n,i}(t_j) - \mu^{(i)}(t_j) \big)
= \frac{1}{\sqrt{n}} \sum_{i=1}^n Z_{i,n} \Rightarrow \bar{Z}=\sum_{j=1}^q c_j Z(t_j)~,
% \end{split}
\end{align}
where $Z_{i,n} = \sum_{j=1}^q c_j (X_{n,i}(t_j) - \mu^{(i)}(t_j) ) $  and $c_1,\dots,c_q\in\mathbb{R}$ denote arbitrary constants.

In order to show the weak convergence of $Z_n$ to $\bar Z$, a blocking technique similar to Lemma A.1 in  B\"ucher and Kojadinovic \cite{buecher2016} is utilized.
%More precisely, we split the sum in \eqref{Zndef} into blocks where each block is composed of a big followed by a small subblock. The idea is that the variance of $Z_n$ is asymptotically the same as the variance of the sum containing only the big subblocks. But at the same time the length of each small subblock is large enough such that we can find a sequence of independent random variables whose sum is asymptotically the same as the sum over the big subblocks.
Each block consists of a big subblock followed by a small subblock.
For $1/(2+2\nu)<\eta_b < \eta_s <1/2$, where $\nu>0$ is the same as in Assumption (A1),
define the length of the small and the big subblocks as
$s_n = \lfloor n^{1/2-\eta_s} \rfloor$ and $b_n = \lfloor n^{1/2-\eta_b} \rfloor$, respectively.
In total there are therefore $k_n = \lfloor n/(b_n+s_n) \rfloor$ blocks. Now consider the random variables
\begin{align*}
B_{j,n} = \sum_{i=(j-1)(b_n+s_n)+1}^{(j-1)(b_n+s_n)+b_n} Z_{i,n} \qquad\mbox{and}\qquad
S_{j,n} = \sum_{i=(j-1)(b_n+s_n)+b_n+1}^{j(b_n+s_n)} Z_{i,n},\qquad j=1,\dots,k_n,
\end{align*}
corresponding to  the sums of the $j$-th big and small subblock, respectively.
It directly follows that
\begin{align} \label{reprZN}
Z_n = \frac{1}{\sqrt{n}} \sum_{j=1}^{k_n} B_{j,n} + \frac{1}{\sqrt{n}} \sum_{j=1}^{k_n} S_{j,n}
+\frac{1}{\sqrt{n}} R_n ~,
\end{align}
where $R_n = \sum_{i=k_n(b_n+s_n)+1}^n Z_{i,n} $ denotes the sum of the
remaining terms $Z_{i,n}$ after the last small subblock. Now the variance
of $Z_n$ can be written as
\begin{align} \label{varZnCLT}
\begin{split}
\text{Var}(Z_n) = & ~ \text{Var}\Big (\frac{1}{\sqrt{n}}\sum_{j=1}^{k_n} B_{j,n}\Big )
+ \frac{2}{n} \sum_{j,j'=1}^{k_n} \mathbb{E}[B_{j,n} S_{j',n}]
+ \frac{2}{n} \sum_{j=1}^{k_n} \mathbb{E}[B_{j,n} R_n]  \\
& +\frac{1}{n} \sum_{j,j'=1}^{k_n} \mathbb{E}[S_{j,n} S_{j',n}]
+ \frac{2}{n} \sum_{j=1}^{k_n} \mathbb{E}[S_{j,n} R_n]
+ \frac{1}{n}\mathbb{E}[R_n^2].
\end{split}
\end{align}
It will be shown that each term except the first on the right-hand side of \eqref{varZnCLT}
converges to zero. As a consequence, it follows that
\begin{align}\label{h1a}
|Z_n - 1/\sqrt{n} \sum_{j=1}^{k_n} B_{j,n}| &= |1/\sqrt{n} \sum_{j=1}^{k_n} S_{j,n}+1/\sqrt{n} R_n| = o_\mathbb{P}(1) , \\
\mbox{Var}(Z_n) & =\text{Var}\bigg(\sum_{j=1}^{k_n} B_{j,n}\bigg) + o(1).
\label{h2a}
\end{align}
First show that the second term on the right side of \eqref{varZnCLT} converges to zero.
To this end,
\begin{align*}
\mathbb{E}[B_{j,n} S_{j',n}] =
\sum_{i=(j-1)(b_n+s_n)+1}^{(j-1)(b_n+s_n)+b_n} \sum_{i'=(j'-1)(b_n+s_n)+b_n+1}^{j'(b_n+s_n)}
\mathbb{E}[Z_{i,n} Z_{i',n}] ~.
\end{align*}
Note that $\sigma(X(t))\subset \sigma(X)$ for any
random variable $X$ in $C([0,1])$ and $t\in[0,1]$, where $\sigma(X)$ denotes the
$\sigma$-field generated by $X$ (see Problem 1.7.1 in Van der Vaart and Wellner \cite{wellner1996}). Hence,
$\phi(\sigma(X(t)),\sigma(Y(s))) \leq \phi(\sigma(X),\sigma(Y))$ for any
$C([0,1])$-valued random variables $X,Y$ and any $s,t\in[0,1]$.
Using this fact and formula (3.17) in  Dehling and  Philipp \cite{dehling2002}
leads to the bound
\begin{align} \label{bound1CLT}
|\mathbb{E}[Z_{i,n} Z_{i',n}]|
& \leq \sum_{l,l'=1}^q |c_l c_{l'}| |\text{Cov}(X_i (t_l), X_{i'}(t_{l'}))| \\
& \leq 2 \varphi(|i-i'|)^{1/2} \sum_{l,l'=1}^q |c_l c_{l'}| ~
\mathbb{E}[X_i(t_l)^2]^{1/2} \mathbb{E}[X_{i'}(t_{l'})^2]^{1/2} \nonumber
\lesssim \varphi(|i-i'|)^{1/2},
\end{align}
where the last inequality follows from Assumption (A1)
and the symbol $\lesssim $ means less or equal up to a constant independent of $n$.
This gives
\begin{align} \label{blockIneq}
|\mathbb{E}[B_{j,n} S_{j,n}]|
&\lesssim  \sum_{i=(j-1)(b_n+s_n)+1}^{(j-1)(b_n+s_n)+b_n}
\sum_{i'=(j-1)(b_n+s_n)+b_n+1}^{j(b_n+s_n)} \varphi(|i-i'|)^{1/2}
%\\
% &\leq \text{const} \times \sum_{i=1}^{b_n} \sum_{i'=b_n+1}^{b_n + s_n} \varphi(|i-i'|)^{1/2} \nonumber
\lesssim   \sum_{i=1}^{b_n +s_n -1} i \varphi(i)^{1/2} < \infty .
\end{align}
The last inequality holds, since Assumption (A4) yields
$ \sum_{i=1}^\infty i \varphi(i)^{1/2} \leq \sum_{i=1}^\infty i a^{i/2} < \infty $.
Similarly it can be shown that
$|\mathbb{E}[B_{j,n} S_{j-1,n}]| < \infty$. For $j'>j+1$ and $j>j'+1$ there is at least one
big subblock between the observations and since $\varphi(\cdot)$ is monotonically
decreasing, it follows that
\begin{align*}
\vert \mathbb{E} [B_{j,n} S_{j',n} ] \vert
\lesssim  \sum_{i=(j-1)(b_n+s_n)+1}^{(j-1)(b_n+s_n)+b_n}
\sum_{i'=(j'-1)(b_n+s_n)+b_n+1}^{j'(b_n+s_n)} \varphi(|i-i'|)^{1/2}
= O(b_n s_n \varphi(b_n)^{1/2} ) ~.
\end{align*}
Overall, it follows that
\begin{align*}
\frac{2}{n} \sum_{j,j'=1}^{k_n} \mathbb{E}[B_{j,n} S_{j',n}]
= O(n^{-1} k_n) + O(n^{-1} k_n^2 b_n s_n \varphi(b_n)^{1/2}) = O(b_n^{-1}) + O(n b_n^{-1} s_n a^{b_n /2}) =o(1),
\end{align*}
using that
$n^{-1}k_n = (b_n + s_n)^{-1} = O(b_n^{-1})$  and  that $s_n/b_n \to 0$ and $na^{b_n/2}\to 0$.

For the third  term in \eqref{varZnCLT}, proceed in a similar way as in \eqref{blockIneq} to get
\begin{align*}
\mathbb{E}[B_{k_n,n} R_n]
& \lesssim  \sum_{i=(k_n-1)(b_n+s_n)+1}^{(k_n-1)(b_n+s_n)+b_n}
\sum_{i'= k_n(b_n+s_n)+1}^n \varphi(|i-i'|)^{1/2}
 \lesssim   \sum_{i=1}^{b_n + (n-k_n(s_n+b_n)) -1} i \varphi(i)^{1/2} < \infty .
\end{align*}
In the case of $j<k_n$, there is again at least one big subblock between the observations in
$\mathbb{E}[B_{j,n}R_n]$ and therefore it follows that
\begin{align*}
\vert \mathbb{E} [B_{j,n} R_{n} ] \vert
&  \lesssim   \sum_{i=(j-1)(b_n+s_n)+1}^{(j-1)(b_n+s_n)+b_n}
\sum_{i'= k_n(b_n+s_n)+1}^n \varphi(|i-i'|)^{1/2} %= O(b_n (n-k_n(b_n+s_n)) \varphi(b_n)^{1/2} )
= O(b_n^2 \varphi(b_n)^{1/2}),
\end{align*}
since $n-k_n(s_n+b_n) \leq s_n +b_n = O(b_n)$.
Altogether, the calculations above yield
\begin{align*}
\frac{2}{n} \sum_{j=1}^{k_n} \mathbb{E}[B_{j,n} R_n]
%= \frac{2}{n} \sum_{j=1}^{k_n-1} \mathbb{E}[B_{j,n} R_n]
% + \frac{2}{n} \mathbb{E}[B_{k_n,n} R_n]
=O(n^{-1} k_n b_n^2 \varphi(b_n)^{1/2})
+ O(n^{-1})
= O(b_n a^{b_n /2}) + O(n^{-1}),
\end{align*}
which converges to zero.

Now consider the fourth term in \eqref{varZnCLT} and  use \eqref{bound1CLT} to conclude
\begin{align} \label{smallIneq}
\mathbb{E}[S_{j,n} S_{j',n}]
&\lesssim   \sum_{i=(j-1)(b_n+s_n)+b_n+1}^{j(b_n+s_n)}
	\sum_{i'=(j'-1)(b_n+s_n)+b_n+1}^{j'(b_n+s_n)} \varphi(|i-i'|)^{1/2}.
\end{align}

For $j=j'$, it holds that
\begin{align} \label{bound2.1CLT}
\mathbb{E}[S_{j,n}^2]
\lesssim  \sum_{i,i'=1}^{s_n} \varphi(|i-i'|)^{1/2}
\lesssim  \sum_{i=0}^{s_n-1} (s_n -i) \varphi(i)^{1/2}
\lesssim  s_n \sum_{i=0}^{\infty} a^{i/2}\lesssim   s_n
\end{align}
and therefore  $\mathbb{E}[S_{j,n}^2] = O(s_n)$. Since there is always at least
one big subblock between two small subblocks, it follows that
$\mathbb{E}[S_{j,n}S_{j',n}] = O(s_n^2 \varphi(b_n)^{1/2})$ for $j\neq j'$.
Hence,
\begin{align*}
\begin{split}
\frac{1}{n} \sum_{j,j'=1}^{k_n} \mathbb{E}[S_{j,n}S_{j',n}]
= O(n^{-1} k_n s_n ) + O(n^{-1} k_n^2 s_n^2 \varphi(b_n)^{1/2})
= O(b_n^{-1} s_n) + O(n b_n^{-2} s_n^2 a^{b_n / 2}),
\end{split}
\end{align*}
which converges to zero, since $b_n^{-1} s_n \to 0$ and $n a^{b_n / 2} \to 0$ as $n\to\infty$.

For the fifth term in \eqref{varZnCLT}, use similar arguments as for the third term and for the last term, use the same arguments as in \eqref{bound2.1CLT} to get
\begin{align*}
\frac{1}{n} \mathbb{E}[R_n^2] = O(n^{-1} (n-k_n(b_n + s_n)) ) = O(n^{-1}b_n) \to 0.
\end{align*}
% where we again used the fact that $n-k_n(b_n + s_n) < b_n+s_n =O(b_n)$.

From \eqref{h1a} it follows that it suffices to show the convergence
$n^{-1/2} \sum_{j=1}^{k_n} B_{j,n} \Rightarrow \bar Z$
in order to establish \eqref{barZdef}. For that purpose let
$\psi_{j,n}(t) = \exp (itn^{-1/2} B_{j,n})$ and define
$ \mathbb{E} \big[ \prod_{j=1}^{k_n} \psi_{j,n}(t)\big]$
as the characteristic function of $n^{-1/2} \sum_{j=1}^{k_n} B_{j,n}$.
Let $B'_{1,n},\dots,B'_{k_n,n}$  denote  independent random variables  such that
$B_{j,n}$ and $B'_{j,n}$ are equally distributed $(j=1,\ldots , k_n)$ and define
$ \prod_{j=1}^{k_n} \mathbb{E} [\psi_{j,n}(t)]$ as the
characteristic function of $n^{-1/2} \sum_{j=1}^{k_n} B'_{j,n}$. Then
\begin{align*}
\Big | \mathbb{E}  \Big [ \prod_{j=1}^{k_n} \psi_{j,n}(t) \Big ]
- \prod_{j=1}^{k_n} \mathbb{E}  \Big  [\psi_{j,n}(t) \Big ]  \Big |  &  \leq \sum_{i=1}^{k_n}
\Big | \prod_{j=1}^{i-1} \mathbb{E} \big  [\psi_{j,n}(t) \big ] \Big |
~ \Big | \mathbb{E} \Big [ \prod_{j=i}^{k_n} \psi_{j,n}(t) \Big ]
- \mathbb{E} \big [\psi_{i,n}(t)\big ]
 \mathbb{E} \Big  [ \prod_{j=i+1}^{k_n} \psi_{j,n}(t) \Big ] \Big  |  \\
 & \lesssim  k_n \max_{1\leq i\leq k_n-1}
\phi  \Big ( \sigma\left( \psi_{i,n}(t) \right),
\sigma \Big ( \prod_{j=i+1}^{k_n} \psi_{j,n}(t) \Big  ) \Big )\lesssim   k_n \varphi(s_n)  = o(1),
\end{align*}
where Lemma 3.9 of Dehling and Philipp \cite{dehling2002} was used for the the second inequality, while the third inequality follows,
since there are always $s_n$ observations between two big subblocks.
Hence, it suffices to show that $n^{-1/2} \sum_{j=1}^{k_n} B'_{j,n}$ converges in distribution
to $\bar Z$ in order to establish \eqref{barZdef}. For that purpose, utilize the Lindeberg--Feller central limit theorem
for triangular arrays. It is first shown that Var$(n^{-1/2} \sum_{j=1}^{k_n} B'_{j,n})$
converges to the variance of $\bar Z$.

Recall that the random variables $B'_{1,n}, \dots, B'_{k_n,n}$ are independent and have
the same distributions as $B_{1,n},\dots,B_{k_n,n}$, respectively. Thus,
\begin{align}\label{varasy}
\text{Var}\Big(\frac{1}{\sqrt{n}} \sum_{j=1}^{k_n} B'_{j,n} \Big )
&= \frac 1n \sum_{j=1}^{k_n} \text{Var}(B'_{j,n})
% = n^{-1} \sum_{j=1}^{k_n} \text{Var}(B_{j,n})
= \text{Var} \Big (\frac{1}{\sqrt{n}} \sum_{j=1}^{k_n} B_{j,n} \Big)
- \frac{1}{n} \sum_{\substack{j,j'=1 \\ j\neq j'}}^{k_n} \mathbb{E}[B_{j,n}B_{j',n}].
\end{align}
It is already known that Var$(n^{-1/2} \sum_{j=1}^{k_n} B_{j,n})=$ Var$(Z_n) + o(1)$. Recall the
calculations in \eqref{smallIneq} and the subsequent discussion to note that it can be shown in
a similar way that $\mathbb{E}[B_{j,n}B_{j',n}] = O(b_n^2 \varphi(s_n)^{1/2})$,
for $j\neq j'$. Consequently,
\begin{align*}
\frac{1}{n} \sum_{\substack{j,j'=1 \\ j\neq j'}}^{k_n} \mathbb{E}[B_{j,n}B_{j',n}]
= O(n^{-1} k_n^2 b_n^2 \varphi(s_n)^{1/2}).
\end{align*}
Since $ n^{-1}k_n = O((b_n +s_n)^{-1}) = O(b_n^{-1}) $,
$O(n^{-1} k_n^2 b_n^2 \varphi(s_n)^{1/2}) = O(n a^{s_n/2}) \to 0$, this means that
Var$(n^{-1/2} \sum_{j=1}^{k_n} B'_{j,n})=$Var$(Z_n) + o(1)$.

Now consider the variance of $\bar Z$. Recall the definition of $\bar Z$ in \eqref{barZdef}
and note that
its variance satisfies
\begin{align} \label{varBarZ}
\mathrm{Var}(\bar Z)
= \mathrm{Var} \Big (\sum_{l=1}^q c_l Z(t_l)\Big )
= \sum_{l,l'=1}^q c_l c_{l'} \mathrm{Cov}(Z(t_l),Z(t_{l'}))
= \sum_{l,l'=1}^q c_l c_{l'}  C(t_l,t_{l'}).
\end{align}
%We have to show that the variance of $Z_n$, which we defined in \eqref{Zndef}, converges to the
%expression above.
The variance of $Z_n$ can be expressed as
\begin{align*}
\mathrm{Var}(Z_n)
&= \mathrm{Var} \Big (\frac{1}{\sqrt{n}} \sum_{i=1}^n Z_{i,n} \Big  )
= \frac{1}{n} \sum_{i,i'=1}^n \mathrm{Cov}(Z_{i,n},Z_{i',n})
=  \frac{1}{n} \sum_{l,l'=1}^q c_l c_{l'} \sum_{i,i'=1}^n \mathrm{Cov}(X_{n,i}(t_l),X_{n,i'}(t_{l'})) \\
&= \frac{1}{n} \sum_{l,l'=1}^q c_l c_{l'} \sum_{i,i'=1}^n \gamma(i-i',t_l,t_{l'})
=  \frac{1}{n} \sum_{l,l'=1}^q c_l c_{l'} \sum_{i=-n}^n (n-|i|) \gamma(i,t_l,t_{l'}) = \mathrm{Var}(\bar Z)  + o(1)
\end{align*}
by the dominated convergence theorem.  Consequently,
$\mathrm{Var}(n^{-1/2} \sum_{j=1}^{k_n} B'_{j,n})=\mathrm{Var}(\bar Z) + o_{\mathbb{P}}(1)$ follows from the previous discussion.

Finally, verify the Lindeberg condition for the random variables $B_{j,n}'$.
Using  H\"older's inequality with $p=1+\nu /2$ and $q=(2+\nu)/\nu$
(here $\nu$ is the same as in Assumption (A1)) and Markov's
inequality yields
\begin{align*}
\mathbb{E}\big[(B_{j,n}')^2
	\mathds{1} \{|B_{j,n}'| > \sqrt{n}\delta\} \big] &=
 \mathbb{E}\big[B_{j,n}^2
	\mathds{1} \{|B_{j,n}| > \sqrt{n}\delta\} \big]   \leq  \mathbb{E} \big[ |B_{j,n}|^{2+\nu}\big]^{2/(2+\nu)}
	\mathbb{P}\big(|B_{j,n}| > \sqrt{n}\delta\big)^{\nu/(2+\nu)} \\
%&= \frac{1}{n} \sum_{j=1}^{k_n} \mathbb{E} \big[ |B_{j,n}|^{2+\nu}\big]^{2/(2+\nu)}
%	\mathbb{P}\big(|B_{j,n}|^{2+\nu} > n^{(2+\nu)/2}\delta^{2+\nu}\big)^{\nu/(2+\nu)} \\
&\leq  \mathbb{E} \big[ |B_{j,n}|^{2+\nu}\big]^{2/(2+\nu)}
	\mathbb{E}\big[|B_{j,n}|^{2+\nu}\big]^{\nu/(2+\nu)} (n^{1/2}\delta)^{-\nu} =  \mathbb{E} \big[ |B_{j,n}|^{2+\nu}\big] (n^{1/2}\delta)^{-\nu} .
\end{align*}
 Minkowski's inequality gives
\begin{align*}
\mathbb{E} \big[ |B_{j,n}|^{2+\nu}\big]^{1/(2+\nu)}
%&= \mathbb{E} \Big[ \Big| \sum_{i=(j-1)(b_n+s_n)+1}^{(j-1)(b_n+s_n)+b_n} Z_{i,n} \Big|^{2+\nu} \Big ]^{1/(2+\nu)}
	 \leq \sum_{i=(j-1)(b_n+s_n)+1}^{(j-1)(b_n+s_n)+b_n}
	\mathbb{E}\big[|Z_{i,n}|^{2+\nu}\big]^{1/(2+\nu)} = O(b_n),
\end{align*}
where the last estimate follows from
\begin{align*}
\max_{1\leq i\leq n} \mathbb{E}\big[|Z_{i,n}|^{2+\nu}\big]^{1/(2+\nu)}
%&= \max_{1\leq i\leq n} \mathbb{E}\Big[
%	\Big|\sum_{l=1}^q c_l (X_{n,i}(t_l)-\mu^{(i)}(t_l))
%	\Big|^{2+\nu}\Big]^{1/(2+\nu)} \\
&\leq \max_{1\leq i\leq n} \sum_{l=1}^q |c_l|
	\mathbb{E}\big[|X_{n,i}(t_l)-\mu^{(i)}(t_l)|^{2+\nu}\big]^{1/(2+\nu)} < \infty
\end{align*}
by  Assumption (A1). Combining these estimates gives the Lindeberg condition, that is,
\begin{align*}
\frac{1}{n} \sum_{j=1}^{k_n} & \mathbb{E}\big[(B_{j,n}')^2
	\mathds{1} \{|B_{j,n}'| > \sqrt{n}\delta\} \big]
= O(n^{-1} k_n b_n^{2+\nu} n^{-\nu/2}) = O(b_n^{1+\nu}  n^{-\nu/2})
= O(n^{1/2-\eta_b(1+\nu)}) =o(1),
\end{align*}
since by assumption $\eta_b > 1/(2+2\nu)$.

The proof is completed showing  that  the  process $(G_n)_{n\in\mathbb{N}}\subset C([0,1])$ is uniformly
$\rho$-equicontinuous in probability with respect to the metric $\rho (s,t) =|s-t|$.
With the notation  $Y_{n,j} := X_{n,j} -\mu^{(j)}$ we obtain 
\begin{align} \nonumber 
\Vert G_n (s) - G_n (t)\Vert_2^2
%=&  \mathbb{E}\bigg[ \Big\vert \frac{1}{\sqrt{n}} \sum_{j=1}^n
%	\big(X_{n,j}(s) - \mu^{(j)}(s) - X_{n,j}(t)+\mu^{(j)}(t) \big) \Big\vert^2 \bigg] \\
%=&  \mathbb{E}\bigg[ \Big\vert \frac{1}{\sqrt{n}} \sum_{j=1}^n
%	\big(Y_{n,j}(s) - Y_{n,j}(t)\big) \Big\vert^2 \bigg] \\
=& \frac{1}{n} \sum_{j,j'=1}^n \mathbb{E} \bigg[
	\big(Y_{n,j}(s) - Y_{n,j}(t) \big)
	\big(Y_{n,j'}(s) - Y_{n,j'}(t) \big) \bigg]  \\
 \leq  &\frac{2}{n} \sum_{j=0}^{n-1}
	(n-j) \mathbb{E}\big[ \big(Y_{n,1}(s) - Y_{n,1}(t) \big)
	\big(Y_{n,1+j}(s) - Y_{n,1+j}(t) \big) \big] ~
 \label{tight0CLT}
\end{align}
(note that each row of the array $\{Y_{n,j} ~\colon~ n\in\mathbb{N},~ j=1,\dots,n \}$ is
stationary). 
Using Assumption (A3), straightforward calculations yield, for any $i=1,\dots, n$,
%\begin{align*}
%| & Y_{n,i}(s) - Y_{n,i}(t) |^2 \\
%&= | X_{n,i}(s) - X_{n,i}(t) + (\mu^{(i)}(s) - \mu^{(i)}(t)) |^2 \\
%&\leq | X_{n,i}(s) - X_{n,i}(t)|^2
%	+ 2| X_{n,i}(s) - X_{n,i}(t)| ~ |\mu^{(i)}(s) - \mu^{(i)}(t) |
%	+ |\mu^{(i)}(s) - \mu^{(i)}(t) |^2 \\
%&\leq M^2 | s-t|^2 + 2 M \mathbb{E}[M] ~ | s-t |^2 + \mathbb{E}[M]^2 ~ | s-t |^2 \\
%\end{align*}
%and therefore we have
\begin{align*}
\mathbb{E} \big[ | Y_{n,i}(s) - Y_{n,i}(t)|^2 \big]^{1/2}
	\lesssim | s-t | ~.
\end{align*}
The inequality above and (3.17) in Dehling and Philipp \cite{dehling2002} together with Assumption (A4) imply,
for any $j=1,\dots,n$,
\begin{align*}
\mathbb{E} & \big[ (Y_{n,1}(s) - Y_{n,1}(t) )
	(Y_{n,1+j}(s) - Y_{n,1+j}(t) ) \big] %\\
%&\leq | \mathrm{Cov} (Y_{n,1}(s) - Y_{n,1}(t), Y_{n,1+j}(s) - Y_{n,1+j}(t) )| \\
\lesssim |s-t|^2 ~ \varphi(j)^{1/2}
\lesssim |s-t|^2 ~ a^{j/2} ~.
\end{align*}
Therefore it follows that 
\begin{align*}
\Vert G_n (s) - G_n (t)\Vert_2^2
\lesssim |s-t|^2 \sum_{j=0}^{\infty} a^{j/2}
	\lesssim |s-t|^2 ,
\end{align*}
and we obtain from Theorem 2.2.4 in Van der Vaart and Wellner \cite{wellner1996}
\begin{align*}
\Big\Vert \sup_{\rho(s,t)\leq \delta} \vert G_n (s)
	- G_n (t) \vert \Big\Vert_2
&\lesssim \int_0^\eta \sqrt{D(\nu,\rho)}~d\nu + \delta D(\eta,\rho) \\
&\lesssim \int_0^\eta \frac{1}{\sqrt{\nu}} ~d\nu
	+  \frac{\delta}{\eta}
	=  2 \sqrt{\eta} +  \frac{\delta}{\eta}  ~.
\end{align*}
where $D(\eta,\rho) = \text{c}
 \lceil \frac{1}{\eta} \rceil  $ is the packing number with respect to the metric $\rho(s,t) = |s-t|$.
Markov's inequality now yields, for any $\varepsilon>0$,
\begin{align*}
\mP \bigg( \sup_{\rho(s,t)\leq \delta} \big\vert G_n (s)
	- G_n (t) \big\vert > \varepsilon \bigg)
\leq \frac{1}{\varepsilon^2} \Big\Vert \sup_{\rho(s,t)\leq \delta} \vert G_n (s)
	- G_n (t) \vert \Big\Vert_2^2
\lesssim \frac{1}{\varepsilon^2}
	\bigg[ 2 \sqrt{\eta} + \frac{\delta}{\eta} \bigg]
\end{align*}
and, since $\eta >0$ is arbitrary, we have 
\begin{align*}
\lim_{\delta\searrow 0} ~ \limsup_{n\to\infty} ~ \mP \bigg(
	\sup_{\rho(s,t)\leq \delta} \big\vert G_n (s)
	- G_n (t) \big\vert > \varepsilon \bigg) = 0~.
\end{align*}
 This means that $(G_n)_{n\in\mathbb{N}}$ is asymptotically uniformly
$\rho$-equicontinuous in probability and the assertion in Theorem \ref{mixingCLT} follows.
\end{proof}

\begin{proof}[\bf Proof of Theorem \ref{WIP}]
Theorem \ref{mixingCLT} implies that
$n^{-1/2} \sum_{j=1}^n (X_{n,j}-\mu^{(j)})\rightsquigarrow Z$
in $C([0,1])$, where $Z$ is centered Gaussian with covariance function $\mathrm{Cov}(Z(s),Z(t)) = C(s,t)$
for $s,t\in[0,1]$. Moreover, Assumption (A4) together with Remark 3.6.4 in Samur \cite{samur1987}
yield, for any $\varepsilon >0$,
\begin{align} \label{samurCondition1}
\lim_{n\to\infty} n \mathbb{P} \bigg(\frac{1}{\sqrt{n}} \Vert X_{n,1}-\mu^{(1)} \Vert > \varepsilon\bigg) = 0.
\end{align}
The assertion now follows from Corollary 3.5 in Samur \cite{samur1987}.
\end{proof}

%%%%%%%%%%%%%%%%%%%%%%%%
\subsection{Proofs of the results  of  Section \ref{sec:two-sample}}
\label{subsec:proof:two-sample}
%%%%%%%%%%%%%%%%%%%%%%%%

\begin{proof}[\bf Proof of Theorem \ref{thm1}]
Note that $\mathbb{N}\times\mathbb{N}$ is a directed set and $[0,1]$ is compact. Now by the weak convergence in \eqref{33}, the claim follows directly from Theorem \ref{thm0}.
\end{proof}

\begin{proof}[\bf Proof of Theorem \ref{th:cb-1}]
Note that
$$
\{ \mu_1 - \mu_2 \in C_{\alpha,m,n} \} = \Big \{ \sup_{t \in [0,1]} \sqrt{n+m} \ \Big | \ \frac {1}{m} \sum^m_{j=1} \widetilde{X_j}(t) - \frac {1}{n} \sum^n_{j=1} \widetilde{Y_j}(t) \  \Big | \ \leq u_{1- \alpha} \Big \},
$$
where $\widetilde{X_j} = X_j - \mu_1$ and $\widetilde {Y_j} = Y_j - \mu_2$. Therefore, it follows from the discussion in Section \ref{testdiff} (applied to the random variables $\widetilde{X_j}$ and $\widetilde{Y_j}$) that
$$
\lim_{m,n \to \infty} \mathbb{P} (\mu_1 - \mu_2 \in C_{\alpha,m,n}) = \lim_{m,n \to \infty} \mathbb{P} (T_{m,n} \leq u_{1- \alpha}) = 1 -  \alpha.
$$
This is the assertion.
\end{proof}

\begin{proof}[\bf Proof of Theorem \ref{2bTheorem}]
Using the notations
\begin{align*}
V_m = \frac{\sqrt{n+m}}{m} \sum_{j=1}^m (X_j - \mu_1)
\qquad\mbox{and}\qquad
W_n = \frac{\sqrt{n+m}}{n} \sum_{j=1}^n (Y_j - \mu_2)
\end{align*}
write $Z_{m,n} = V_m + W_n$. Next  define
\begin{align*}
U_m^{(r)}  &= \frac{\sqrt{n+m}}{m} \sum_{i=1}^{m-l_1+1}
	\frac{1}{\sqrt{l_1}}\Big( \sum_{j=i}^{i+l_1-1} X_{j}(t)
	-\frac{l_1}{m}\sum_{j=1}^m X_{j}(t) \Big) \xi_i^{(r)} , \\
G_n^{(r)}  &= \frac{\sqrt{n+m}}{n} \sum_{i=1}^{n-l_2+1} \frac{1}{\sqrt{l_2}}\Big( \sum_{j=i}^{i+l_2-1}
	Y_{j}(t) -\frac{l_2}{n}\sum_{j=1}^n Y_{j}(t) \Big) \zeta_i^{(r)} ~,
\end{align*}
then
%\begin{align} \label{vectorSplit}
$(Z_{m,n},\hat{B}_{m,n}^{(1)},\dots,\hat{B}_{m,n}^{(R)})
= (V_m, U_m^{(1)}, \dots, U_m^{(R)}) + (W_n, G_n^{(1)}, \dots, G_n^{(R)}),
$
%\end{align}
and with similar but easier arguments as in the proof of Theorem \ref{bTheorem}, it can be shown that
\begin{align*}
(V_m, U_m^{(1)}, \dots, U_m^{(R)})
	&\rightsquigarrow \tfrac{1}{\sqrt{\lambda}} (Z_1, Z_1^{(1)},\dots,Z_1^{(R)}), \\
(W_n, G_n^{(1)}, \dots, G_n^{(R)})
	&\rightsquigarrow \tfrac{1}{\sqrt{1-\lambda}} (Z_2, Z_2^{(1)},\dots,Z_2^{(R)})
\end{align*}
in $C([0,1])^{R+1}$. Since the two vectors  are independent, it directly
follows that
$
(Z_{m,n},\hat{B}_{m,n}^{(1)},\dots,\hat{B}_{m,n}^{(R)})
\rightsquigarrow (Z, Z^{(1)},\dots,Z^{(R)})
$
in $C([0,1])^{R+1}$.
\end{proof}

\begin{proof}[\bf Proof of Theorems \ref{th:cb-2} and \ref{bootstrap:classic}]
The statement  \eqref{bootclasstestLVL} of  Theorem \ref{bootstrap:classic} is a direct consequence of
 Proposition F.1 in the online supplement of B\"ucher and Kojadinovic \cite{buecher2016}
 (note that the continuity of the random variable \(T\) is implied by the results in
Gaenssler et al. \cite{Gaenssler2007}).

 Now Theorem  \ref{th:cb-2}
 follows from the well-known relation between confidence sets and tests for simple hypotheses (see for example Lehmann
 \cite{lehmann1986}, p.\ 214),  observing  that
 $$
 C_{\alpha,m,n} = \bigg\{ f \in C([0,1])\colon \big\| \hat  \mu_1 - \hat \mu_2 - f  \big\|   \leq \frac{T_{m,n}^{\{\lfloor R(1-\alpha)\rfloor\}}}{\sqrt{n+m}} \bigg \}
 $$
 is defined by the acceptance region of an asymptotic level $\alpha$  test for the hypotheses $H_{0,f}\colon  \mu_1-\mu_2 \equiv f$ versus  $H_{0,f}\colon \mu_1-\mu_2 \not \equiv f$,
 which rejects the null hypothesis, whenever
$$
 \big\| \hat  \mu_1 - \hat \mu_2 - f  \big\|   >  \frac{T_{m,n}^{\{\lfloor R(1-\alpha)\rfloor\}}}{\sqrt{n+m}}.
$$
The remaining statement regarding consistency follows from similar arguments as given in  B\"ucher et al. \cite{buecher2017}. 
Under the alternative, Theorem \ref{thm1} yields, for any \(K \in \mathbb{N}\),
\begin{align} \label{bootclasstestCons1}
\lim_{m,n \to \infty} \mathbb{P}(\sqrt{n+m} ~ \hat{d}_\infty \geq K )
  = \lim_{m,n \to \infty} \mathbb{P}(\sqrt{n+m}(\hat{d}_\infty-d_\infty)
  + \sqrt{n+m} ~ d_\infty \geq K ) = 1.
\end{align}
Moreover, under the alternative, Theorem \ref{2bTheorem} together with the continuous mapping theorem
yield, for any \(R \in \mathbb{N}\),
\begin{align} \label{bootclasstestCons2}
\lim_{K\to\infty} \lim_{m,n\to\infty} \mathbb{P}(T_{m,n}^{\{\lfloor R(1-\alpha)\rfloor\}} > K)
  \leq \lim_{K\to\infty} \lim_{m,n\to\infty} \mathbb{P} (\max_{r=1, \dots , R} T_{m,n}^{(r)} > K) = 0.
\end{align}
The assertion in \eqref{bootclasstestCons} now follows from \eqref{bootclasstestCons1} and
\eqref{bootclasstestCons2}.
\end{proof}

\begin{proof}[\bf Proof of Theorem \ref{setConvergence}]
Straightforward calculations show that
\begin{align*}
\D_{m,n} &=
 \big |\hat{d}_\infty - d_\infty \big  |  +
\sup_{t\in [0,1]}  \big | \hat{\mu}_1(t) - \hat{\mu}_2(t) - (\mu_1(t) - \mu_2(t)) \big | \\
& \geq \big  |\hat{d}_\infty - d_\infty \big |
+ \sup_{t\in \mathcal{E}^+} \big | \hat{\mu}_1(t) - \hat{\mu}_2(t) - d_\infty \big |  \\
&=  \sup_{t\in \mathcal{E}^+} \big(  \big |\hat{d}_\infty - d_\infty\big  |
+\big  |   d_\infty -(\hat{\mu}_1(t) - \hat{\mu}_2(t))  \big | \big) \\
&
\geq \sup_{t\in \mathcal{E}^+} \big  |\hat{d}_\infty -(\hat{\mu}_1(t) - \hat{\mu}_2(t)) \big | \\
& \geq  \hat{d}_\infty - \inf_{t\in \mathcal{E}^+} \big  (\hat{\mu}_1(t) - \hat{\mu}_2(t) \big ),
\end{align*}
and  \eqref{33}  and Theorem \ref{thm1} yield (observing the order of the sequence $c_{m,n}$) that
$
\mP \big( \frac{c_{m,n}}{\sqrt{m+n}} \geq \D_{m,n}  \big)
\to 1
$
as $m,n \to \infty$. Therefore,
\begin{align*}
 \mP \big( \mathcal{E}^+ \subset \hat{\mathcal{E}}^+_{m,n} \big) =
 \mP \Big( \inf_{t\in \mathcal{E}^+} (\hat{\mu}_1(t) - \hat{\mu}_2(t))
\geq \hat{d}_\infty - \frac{c_{m,n}}{\sqrt{m+n}} \Big)  \geq \mP \Big( \frac{c_{m,n}}{\sqrt{m+n}} \geq \D_{m,n}  \Big)
\to 1
\end{align*}
as $m,n \to \infty$. This means that, for any $\varepsilon > 0$,
\begin{align} \label{haus1}
\mP \big( \sup_{x\in\mathcal{E}^+} \inf_{y\in \hat{\mathcal{E}}^+_{m,n}} |x-y|
> \varepsilon \big) \to 0
\end{align}
as $m,n \to \infty$.
On the other hand, the inequality
\begin{align*}
& \sup_{x_{m,n} \in \hat{\mathcal{E}}^+_{m,n}} |\mu_1(x_{m,n}) - \mu_2(x_{m,n})
- d_\infty| \leq
  \sup_{x_{m,n} \in \hat{\mathcal{E}}^+_{m,n}} |\mu_1(x_{m,n}) - \mu_2(x_{m,n})
- \hat{d}_\infty| + |\hat{d}_\infty - d_\infty| \\
&\leq   \sup_{x_{m,n} \in \hat{\mathcal{E}}^+_{m,n}}  \Big(
|\mu_1(x_{m,n}) - \mu_2(x_{m,n})
- (\hat{\mu}_1(x_{m,n}) - \hat{\mu}_2(x_{m,n}))|
+ |\hat{d}_\infty - d_\infty| +  \frac{c_{m,n}}{\sqrt{m+n}} \Big) \nonumber
\end{align*}
shows that
\begin{align*}
& \mP \Big( \sup_{x_{m,n} \in \hat{\mathcal{E}}^+_{m,n}} |\mu_1(x_{m,n}) - \mu_2(x_{m,n})
- d_\infty| \leq  2 \frac{c_{m,n}}{\sqrt{m+n}} \Big) \nonumber \\
&\geq \mP \Big( \sup_{t \in [0,1]} |\mu_1(t) - \mu_2(t)
- (\hat{\mu}_1(t) - \hat{\mu}_2(t))|
+ |\hat{d}_\infty - d_\infty|
\leq \frac{c_{m,n}}{\sqrt{m+n}} \Big) \to 1
\end{align*}
as $m,n \to \infty$. Therefore,
\begin{align}
 \label{probBound}
  \sup_{x_{m,n} \in \hat{\mathcal{E}}^+_{m,n}} |\mu_1(x_{m,n}) - \mu_2(x_{m,n})
- d_\infty| \to 0
\end{align}
(outer) almost surely, which implies
\begin{align} \label{haus2}
\sup_{x_{m,n}\in\hat{\mathcal{E}}^+_{m,n}} \inf_{y\in \mathcal{E}^+} |x_{m,n}-y|  \to 0
\end{align}
(outer) almost surely  as $m,n \to \infty$. To see this, assume the contrary. Then there would exist a sequence
$(x_{m,n})$ in $\hat{\mathcal{E}}^+_{m,n}$ with
$\inf_{y\in \mathcal{E}^+} |x_{m,n}-y| \geq \varepsilon$. Now consider  the closed set
$\mathcal{U}_\varepsilon = \{ t \in \hat{\mathcal{E}}^+_{m,n} \ \colon  \ |t-s| \geq \varepsilon,
~\forall s \in \mathcal{E}^+ \} $. Then
 $x_{m,n} \in \mathcal{U}_\varepsilon $ and
$$\max \{ \mu_1(t) - \mu_2(t) \ \colon  \ t \in \mathcal{U}_\varepsilon \} < d_\infty,$$
which contradicts \eqref{probBound} and proves \eqref{haus2}.
Combining \eqref{haus1} and \eqref{haus2}, for any $\varepsilon > 0 $,
\begin{align*}
\mP \big( d_H (\hat{\mathcal{E}}^+_{m,n} , \mathcal{E}^+) > \varepsilon \big)
= \mP \Big( \max \big\{ \sup_{x\in\hat{\mathcal{E}}^+_{m,n}} \inf_{y\in \mathcal{E}^+}
|x-y|, \sup_{x\in\mathcal{E}^+} \inf_{y\in \hat{\mathcal{E}}^+_{m,n}} |x-y| \Big\}
> \varepsilon \big) \to 0
\end{align*}
as $m,n \to \infty$.  In  the same way it can be shown that
$
\mP \big( d_H (\hat{\mathcal{E}}^-_{m,n} , \mathcal{E}^-) > \varepsilon \big) \to 0
$
as $m,n \to \infty$,  and the assertion of of Lemma  \ref{setConvergence} follows.
\end{proof}

\begin{proof}[\bf Proof of Theorem \ref{2jointConvergence}] The proof is a direct
consequence of the following lemma, which might be of own interest
and will be proved at the end of this section.
\end{proof}

\begin{lemma} \label{maxSetConvergence}
Let $((X_\alpha, X^{(1)}_\alpha, \dots X^{(R)}_\alpha)\colon\alpha \in A)$
be a net of random variables in $C([0,1])^{R+1}$ and 
$((M^{(1)}_\alpha, M^{(2)}_\alpha) \colon$ $ \alpha \in A)$ a net of random elements in the
set $K([0,1])$ of all compact subsets of the interval $[0,1]$. Furthermore, let
$((X_\alpha, X^{(1)}_\alpha, \dots X^{(R)}_\alpha)\colon\alpha \in A)$ converge weakly to
$(X, X^{(1)},\dots, X^{(R+1)}) $ in $C([0,1])^{R+1}$, where
$X^{(1)},\dots, X^{(R+1)}$ are independent copies of $X$, and
$((M^{(1)}_\alpha,M^{(2)}_\alpha) \colon\alpha \in A)$ converge in
probability to the non random  sets $(M^{(1)},M^{(2)})$ in $K([0,1])^2$, that is
$$
\mP \big (
\max \big \{ d_H (M^{(1)}_\alpha,M^{(1)}), d_H (M^{(2)}_\alpha,M^{(2)}) \big\} > \varepsilon \big ) \to 0
$$
for all $\varepsilon >0$. Then the random variables
\begin{align*}
Y_\alpha &= \max\Big\{\max_{t\in M^{(1)}} X_\alpha(t),
\max_{t\in M^{(2)}} (-X_\alpha(t))\Big\} \quad\mbox{and}\quad
Y^{(r)}_\alpha = \max\Big\{\max_{t\in M^{(1)}_\alpha} X^{(r)}_\alpha(t),
\max_{t\in M^{(2)}_\alpha} (-X^{(r)}_\alpha(t))\Big\} ,
\end{align*}
($r = 1,\dots,R$)  satisfy
\begin{align*}
(Y_\alpha, Y^{(1)}_\alpha, \dots, Y^{(R+1)}_\alpha)
\Rightarrow (Y, Y^{(1)}, \dots, Y^{(R)})
\end{align*}
in $\mR^{R+1}$ where
$Y = \max\big\{\max_{t\in M^{(1)}} X(t),
\max_{t\in M^{(2)}} (-X(t))\big\} $
and $Y^{(1)},\dots,Y^{(R)}$ are independent copies of $Y$.
\end{lemma}

Theorem \ref{2bTheorem}, Lemma \ref{setConvergence} and
Lemma \ref{maxSetConvergence} show that
\begin{align*}
\big(\tilde{D}_{m,n}(\mathcal{E}),~ K_{m,n}^{(1)},\dots,K_{m,n}^{(R)}\big)
\Rightarrow (T(\mathcal{E}),~ T^{(1)}(\mathcal{E}),\dots,T^{(R)}(\mathcal{E})),
\end{align*}
in $\mathbb{R}^{R+1}$, where
\begin{align*}
\tilde{D}_{m,n}(\mathcal{E}) = \sqrt{n+m} ~
	\max \Big\{\sup_{t\in\mathcal{E}^+} Z_{m,n},
	 \sup_{t \in\mathcal{E}^-} (-Z_{m,n})\Big\}
\end{align*}
and the random variable $Z_{m,n}$ is defined by \eqref{33}. The proofs of Lemmas \ref{lem:proofs:weak_conv:1} and  \ref{lem:proofs:weak_conv:2} yield that
\begin{align*}
\sqrt{n+m}  (\hat d_\infty - d_\infty)
	= \tilde{D}_{m,n}(\mathcal{E}) + o_{\mathbb{P}}(1),
\end{align*}
thus completing the proof of Theorem \ref{2jointConvergence}.

\begin{proof}[\bf Proof of Lemma \ref{maxSetConvergence}]
Only the convergence $Y^{(1)}_\alpha \Rightarrow Y^{(1)}$ is shown, since the general case follows from
very similar albeit notationally more complex arguments.

Let $\mathcal{S}$ denote the space $C([0,1]) \times K([0,1])^2 $ equipped with the metric $d$ defined by
$$ d((f,A,B), (g,A',B')) = \max \big\{ \|f-g\|_\infty , d_H(A,A'), d_H(B,B') \big\}$$
for any $(f,A,A'), (g,B,B') \in S$. Slutsky's
theorem then yields
$(X^{(1)}_\alpha , M^{(1)}_\alpha, M^{(2)}_\alpha ) \rightsquigarrow (X, M^{(1)}, M^{(2)})$
in $\mathcal{S}$. Now consider the function
$$
\mathcal{H}\colon
 \mathcal{S}  \to \mR, \;\, (f,A,A')  \mapsto \max \{ \max_{t\in A} f(t), \max_{t\in A'} (-f(t)) \}
$$
and show that it is continuous. For this purpose it is sufficient to prove that the
maps
$$
h\colon \mathcal{S}  \to \mR, \;\, (f,A,A')  \mapsto \max_{t\in A} f(t)
\qquad\mbox{and}\qquad
h^-\colon \mathcal{S}  \to \mR \;\, (f,A,A')  \mapsto \max_{t\in A'} (-f(t))
$$
are continuous (since the maximum of two continuous functions is  continuous).
Exemplarily, focus on the map $h$ (the corresponding statement for $h^-$ follows from closely analogous arguments). Let   $\varepsilon >0 $ be arbitrary and $(f,A,A'),(g,B,B') \in \mathcal{S}$ such that
$$d((f,A,A'), (g,B,B')) < \min\{\delta , \varepsilon/2\},$$
where $\delta >0$ is chosen
such that, for any $t,s \in T$ with $|t-s| < \delta$, the inequality
$|g(t)-g(s)| < \varepsilon/2$ holds (note that  $g$ is uniformly continuous on $[0,1]$). Then,
\begin{align*}
|h(f,A,A') - h(g,B,B')| = |\max_{t\in A} f(t) - \max_{t\in B} g(t)|
\leq |\max_{t\in A} f(t) - \max_{t\in A} g(t)| + |\max_{t\in A} g(t) - \max_{t\in B} g(t)|,
\end{align*}
where the first term can be bounded by
$\max_{t\in A} |f(t)-g(t)|\leq \|f-g\|_\infty < \varepsilon/2$. For the second term, first
consider the case where $\max_{t\in A} g(t) \geq \max_{t\in B} g(t)$ and obtain
\begin{align*}
|\max_{t\in A} g(t) - \max_{t\in B} g(t)| = \max_{t\in A} g(t) - \max_{t\in B} g(t)
= g(t_1) - \max_{t\in B} g(t),
\end{align*}
where $t_1 \in {\arg \max}_{t\in A} g(t)$. Since $d_H(A,B) < \delta$, there is a
$t_2 \in B$ such that $|t_1-t_2| < \delta$. Thus,
$$ g(t_1) - \max_{t\in B} g(t) \leq g(t_1) - g(t_2) < \varepsilon/2 ~.$$
For the second case, use the same arguments to show that
$|\max_{t\in A} g(t) - \max_{t\in B} g(t)| < \varepsilon/2$. This yields
\begin{align*}
|h(f,A,A') - h(g,B,B')| < \varepsilon/2 + \varepsilon/2 = \varepsilon,
\end{align*}
verifying that $h$ is continuous.

The discussion at the beginning of this proof shows that
$$Y^{(1)}_\alpha = \max\big\{\max_{t\in M^{(1)}_\alpha} X^{(1)}_\alpha(t),
\max_{t\in M^{(2)}_\alpha} (-X^{(1)}_\alpha(t))\big\}
\Rightarrow Y^{(1)} . $$
The proof is complete.
\end{proof}

%%%%%%%%%%%%%%%%%%%%%%%%
\subsection{Proofs of the results in Section \ref{sec:change-point}}
\label{sec:proofs:change-point}
%%%%%%%%%%%%%%%%%%%%%%%%

\begin{proof}[\bf Proof of Theorem \ref{thm2}]
Recall the notations \eqref{39}, \eqref{313} and \eqref{318} to obtain
\begin{align*}
&\sqrt{n} \left(\hat {\mathbb{U}}_n (s,t) - \big( s \wedge s^* - s s^*   \big)  \big (\mu_1 (t) - \mu_2(t) \big) \right) \\
%& ~~= \sqrt{n} \left(\hat {\mathbb{U}}_n (s,t) - \mathbb{E} [\hat {\mathbb{U}}_n (s,t)] \right)
%	+ \sqrt{n} \left(\mathbb{E} [\hat {\mathbb{U}}_n (s,t)]
%	- \big( s \wedge s^* - s s^*   \big)  \big (\mu_1 (t) - \mu_2(t) \big) \right) \\
& ~~= \ \hat {\mathbb{W}}_n (s,t)
	+ \sqrt{n} \left(E[\hat {\mathbb{U}}_n (s,t)]
	- \big( s \wedge s^* - s s^*   \big)  \big (\mu_1 (t) - \mu_2(t) \big) \right) \\
& ~~=  \ \hat {\mathbb{W}}_n (s,t) + o_\mathbb{P}(1)
\end{align*}
 uniformly with respect to $(s,t) \in [0,1]^2$. An application of \eqref{315} shows that the weak convergence
\[
\sqrt{n} \Big(\hat {\mathbb{U}}_n (s,t) - \big( s \wedge s^* - s s^*   \big)  \big (\mu_1 (t) - \mu_2(t) \big) \colon (s,t) \in [0,1]^2\Big)
\rightsquigarrow
\big( {\mathbb{W}}(s,t) \colon (s,t) \in [0,1]^2\big)
\]
follows. Now the representation
\begin{align*}
\mathbb{D}_n
% &= \sqrt{n} \Big( \mathbb{\hat  M}_n   - s^*(1-s^*)  \| \mu_1 - \mu_2 \|_\infty \Big) \\
&= \sqrt{n} \Big ( \sup_{(s,t) \in [0,1]^2} \big| \hat {\mathbb{U}}_n (s,t) \big|
- \sup_{(s,t) \in [0,1]^2} \left| \big( s \wedge s^* - s s^*   \big)  \big (\mu_1 (t) - \mu_2(t) \big) \right| \Big)
\end{align*}
and Theorem \ref{thm0} yield the assertion of Theorem \ref{thm2}.
\end{proof}

\begin{proof}[\bf Proof of Theorem \ref{rate}]
In order to prove the assertion, use Corollary 2 in Hariz et al.\ \cite{hariz2007}.
Let $\mathcal{M}$ be the space of all signed finite measures on $C([0,1])$,
define  $\pi_t$ as the canonical projection   $ C([0,1])  \ni x \to x(t)$ and consider the class
$\mathcal{F}= \{\pi_t\colon t\in [0,1]\}$.
Note that
\begin{align*}
N(\nu) = \sup_{t\in[0,1]} \Big| \int_{C([0,1])} \pi_t(x) ~ \nu(dx) \Big|
\end{align*}
defines a semi-norm on $\mathcal{M}$. In particular, if $P=  \mP^{X_1}$ and $Q=\mP^{X_n}$
are the distributions on $C([0,1])$  before and after the change-point, it holds that
\begin{align*}
\int_{C([0,1])} \pi_t(x) (P-Q)(dx)
&= \int_{C([0,1])} \pi_t(x) (\mP^{X_1}-\mP^{X_n})(dx)
= \int_\Omega \pi_t(X_1) d\mP - \int_\Omega \pi_t(X_n) d\mP \\
&= \int_\Omega X_1(t) d\mP - \int_\Omega X_n(t) d\mP
= \mE [X_1(t)] - \mE [X_n(t)]
= \mu_1(t) - \mu_2(t),
\end{align*}
and therefore $N(P-Q) = \|\mu_1 - \mu_2\| > 0$.
The  estimator  of the change-point can now be rewritten as
$
\tilde{s} = {n}^{-1} \min \big( \arg \max_{1\leq k <n} \{ N(D_k) \} \big),
$
where
\begin{align*}
D_k =  \frac{k}{n} \Big(1-\frac{k}{n}\Big)
	\bigg( \frac{1}{k} \sum_{i=1}^k \delta_{X_{n,i}} - \frac{1}{n-k} \sum_{i=k+1}^n
	\delta_{X_{n,i}}
	\bigg)
\end{align*}
and $\delta_x$ denotes the Dirac measure at the point $x \in C([0,1])$.
The assertion of Theorem \ref{rate} follows from Corollary 2 and Remark 2 in Hariz et al.
\cite{hariz2007}
if the following conditions can be verified.

There exist  constants $C>0$ and $\xi >0$ such that, for any $p$,
\begin{align} \label{rateCond1}
\sup_{t\in [0,1]} \sup_{k+p \leq n} \mE \bigg[\Big( \sum_{i=k}^{k+p}
	\big(\pi_t(X_{n,i}) - \mE[\pi_t(X_{n,i})]\big)\Big)^2\bigg] \leq Cp^{2-\xi} .
\end{align}
For any $\varepsilon >0$,
\begin{align} \label{rateCond2}
N_{[~]}(\varepsilon,\mathcal{F}, \|\cdot\|_{\mathcal{G}}) < \infty ~,
\end{align}
where  $\|\cdot\|_{\mathcal{G}}$ is a norm on a space $\mathcal{G}$ (containing $\mathcal{F}$), which
satisfies  $(|P(|\pi_t|)| + |Q(|\pi_t|)|) \leq \|\pi_t\|_{\mathcal{G}}$
for any $t\in [0,1]$ and $N_{[~]}(\varepsilon,\mathcal{F}, \|\cdot\|_{\mathcal{G}})$ denotes the
\textit{bracketing number}. Moreover, for any $\nu\in\mathcal{M}$ and
$f:C([0,1]) \to \mathbb{R}$, define
\begin{align*}
\nu(f) = \int_{C([0,1])} f(x) \,\nu(dx) .
\end{align*}
%Furthermore, we have to verify that
%\begin{align} \label{rateCond4}
%N(\nu) \leq  \sup_{t \in [0,1]} |\nu(\pi_t)|~.
%\end{align}

%If $T:C([0,1])\to \mathbb{R}$ is a continuous linear operator and $f:C([0,1])\to C([0,1])$
%is Bochner integrable with respect to $\nu$, then $T(f)$ is integrable and
%\begin{align*}
%T\Big( \int_{C([0,1])} f(x) ~\nu(dx) \Big) = \int_{C([0,1])} T(f(x)) ~\nu(dx)
%\end{align*}
%(cf. \cite{einar1965}). If the identity map on $C([0,1])$ is integrable,
%i.e. $\int_{C([0,1])} \|x\|_{\infty} ~\nu(dx) < \infty$, we obtain
%\begin{align*}
%N(\nu) = \sup_{t\in [0,1]} \Big| \pi_t \Big( \int_{C([0,1])} x ~\nu(dx) \Big) \Big|
%	= \sup_{t\in [0,1]} \Big| \int_{C([0,1])} \pi_t(x) ~\nu(dx) \Big|
%	= \sup_{t\in [0,1]} |\nu(\pi_t)| ~,
%\end{align*}
%which means that \eqref{rateCond4} is satisfied \textbf{(kann es passieren, dass die rechte
%Seite existiert aber die linke nicht?)}.

Since Assumption (A3) is satisfied with a bounded random variable $M$, only consider the subspace of $C([0,1])$ that consists of
all functions that are Lipschitz continuous with a uniform constant $c$. Therefore,
\begin{align*}
|\pi_s(x)-\pi_t(x)| = |x(s)-x(t)| \leq c |s-t|
\end{align*}
for any $s,t\in[0,1]$. It follows from Theorem 2.7.11 in Van der Vaart and Wellner  \cite{wellner1996} that
\begin{align*}
N_{[~]} (2 \tilde{\varepsilon} \|c\|_{\mathcal{G}}, \mathcal{F}, \|\cdot\|_{\mathcal{G}})
	\leq N(\tilde{\varepsilon}, [0,1], |\cdot|) .
\end{align*}
In the equation above $N(\tilde{\varepsilon}, [0,1], |\cdot|)$ denotes the covering number, that
is the minimal number of balls of radius $\tilde{\varepsilon}$ needed to cover the unit interval $[0,1]$.
Note that
\begin{align*}
N(\tilde{\varepsilon}, [0,1], |\cdot|) = \Big\lceil \frac{1}{2 \tilde{\varepsilon}} \Big\rceil
< \infty
\end{align*}
and therefore \eqref{rateCond2} is satisfied.

Using Assumptions (A1), (A4) and (3.17) in Dehling and Philipp \cite{dehling2002}, leads to
\begin{align*}
\sup_{t\in [0,1]} & \sup_{k+p \leq n} \mathbb{E} \bigg[
\Big( \sum_{i=k}^{k+p}
	\big(\pi_t(X_{n,i}) - \mathbb{E}[\pi_t(X_{n,i})]\big)\Big)^2\bigg]
= \sup_{t\in [0,1]} \sup_{k+p \leq n} \mathbb{E} \bigg[\Big(\sum_{i=k}^{k+p}
	\big(X_{n,i}(t) - \mathbb{E}[X_{n,i}(t)]\big)\Big)^2\bigg] \\
&= \sup_{t\in [0,1]} \sup_{k+p \leq n} \sum_{i,j=k}^{k+p} \mathrm{Cov}(X_{n,i}(t), X_{n,j}(t))
\lesssim \sup_{k+p \leq n} \sum_{i,j=k}^{k+p} \varphi(|i-j|)^{1/2} \\
&\lesssim  \sum_{i=0}^{p} (p+1-i) \varphi(i)^{1/2}
\lesssim  (p+1) \sum_{i=0}^{\infty} a^{i/2}
	\lesssim  p
\end{align*}
which means that \eqref{rateCond1} is satisfied for $\xi =1$.
\end{proof}

\begin{proof}[\bf Proof of Theorem \ref{bTheorem}]
The weak convergence
\begin{align} \label{61}
(\hat{\mathbb{V}}_n,\hat{B}_n^{(1)},\dots,\hat{B}_n^{(R)})
\rightsquigarrow (\mathbb{V}, \mathbb{V}^{(1)},\dots,\mathbb{V}^{(R)})
\end{align}
in $C([0,1]^2)^{R+1}$ will be verified, where $\hat{\mathbb{V}}_n$ is defined in \eqref{39} and $\mathbb{V}^{(1)},\dots,\mathbb{V}^{(R)}$ are independent copies of\ $\mathbb{V}$ defined in \eqref{312a}.
The assertion then follows from the continuous mapping theorem.
We only consider the case \(d_\infty > 0\) since, for \(d_\infty = 0\), the assertion follows by similar arguments.
The proof of \eqref{61} is complicated and consists of a series of steps, which are described first.

\begin{itemize}
\item[(1)]  Replace the estimates $\hat \mu_1$ and $\hat \mu_2$ in \eqref{bProcess}
 by the true functions $\mu_1$ and $\mu_2$ and consider the process
\begin{align*}
\tilde{B}_n^{(k)}(s,t) =& \frac{1}{\sqrt{n}} \sum_{i=1}^{\lfloor sn \rfloor}
\frac{1}{\sqrt{l}}
\Big ( \sum_{j=i}^{i+l-1} \big(Y_{n,j}(t) - \mu_1(t)\big) \Big ) \xi_i^{(k)} \\
	&+ \sqrt{n}\Big (s - \frac{\lfloor sn \rfloor}{n} \Big )\frac{1}{\sqrt{l}}
	\Big ( \sum_{j=\lfloor sn \rfloor +1}^{\lfloor sn \rfloor+l}
	\big(Y_{n,j}(t)	- \mu_1(t)\big) \Big ) \xi_{\lfloor sn \rfloor +1}^{(k)} ,
\end{align*}
where  $Y_{n,j} = X_{n,j} - (\mu_2 - \mu_1) \mathds{1}\{j> \lfloor s^*n \rfloor \}$. For this process, show the weak convergence
 \begin{align}  \label{bLemma1}
(\hat{\mathbb{V}}_n,\tilde{B}_n^{(1)},\dots,\tilde{B}_n^{(R)})
\rightsquigarrow (\mathbb{V}, \mathbb{V}^{(1)},\dots,\mathbb{V}^{(R)}) .
\end{align}
\item[(2)] Next show that
\begin{align} \label{bLemma2}
\sup_{(s,t) \in [0,1]^2} |\tilde {B}_n^{(k)} (s,t) - \bar{B}_n^{(k)}(s,t)| = o_{\mathbb{P}} (1)~,
\end{align}
where the process $\bar{B}_n^{(k)}$ is defined by
\begin{align} \label{Bbar}
\begin{split}
\bar{B}_n^{(k)}(s,t) =& \frac{1}{\sqrt{n}} \sum_{i=1}^{\lfloor sn \rfloor} \frac{1}{\sqrt{l}}
\Big( \sum_{j=i}^{i+l-1} Y_{n,j}(t) - \frac{l}{n} \sum_{j=1}^n Y_{n,j}(t) \Big) \xi_i^{(k)} \\
&+ \sqrt{n}\Big(s - \frac{\lfloor sn \rfloor}{n} \Big)\frac{1}{\sqrt{l}}
\Big( \sum_{j=\lfloor sn \rfloor +1}^{\lfloor sn \rfloor+l} Y_{n,j}(t)
- \frac{l}{n} \sum_{j=1}^n Y_{n,j}(t) \Big) \xi_{\lfloor sn \rfloor +1}^{(k)}.
\end{split}
\end{align}
\item[(3)] Finally establish the assertion \eqref{61}, proving for $k=1,\dots,R$, that
\begin{align} \label{bLemma3}
\sup_{(s,t) \in [0,1]^2} |\hat{B}_n^{(k)} (s,t) - \bar{B}_n^{(k)}(s,t)| = o_{\mathbb{P}} (1)~.
\end{align}
\end{itemize}
Combining (1)--(3) completes the proof.
\end{proof}

\begin{proof}[\bf Proof of \eqref{bLemma1}]
The weak convergence of the process will be shown through proving the weak convergence of its finite-dimensional distributions and asymptotic tightness.

\medskip \noindent
{\it  (A) Convergence of the finite-dimensional distributions}. This part uses similar arguments as given in the proof of Theorem \ref{mixingCLT} and detailed arguments are only given when substantial differences occur. For  the sake of brevity and simpler notations, consider the case $R=1$ and prove
\begin{align} \label{Zn'toZ}
\hat  Z_n =
\sum_{j=1}^q c_j \hat{\mathbb{V}}_n (s_j,u_j)
	+ \sum_{j=1}^q d_j \tilde{B}_n^{(1)}(t_j,v_j)
\Rightarrow Z = \sum_{j=1}^q c_j \mathbb{V} (s_j,u_j)
	+ \sum_{j=1}^q d_j \mathbb{V}^{(1)}(t_j,v_j)
\end{align}
for any  $q\in\mathbb{N}$,  $(u_1,s_1,v_1,t_1),\dots,(u_q,s_q,v_q,t_q)\in[0,1]^4$   and
arbitrary constants $c_1,d_1,\dots,c_q,d_q \in\mathbb{R}$.
The convergence of the finite-dimensional distributions of  the process $(\hat{\mathbb{V}}_n ,
	\tilde{B}_n^{(1)}) $ then follows from an application of the Cram{\'e}r--Wold device.  For this purpose define
	\begin{align*}
 Z_n = \sum_{j=1}^q c_j \check{\mathbb{V}}_n (s_j,u_j)
	+ \sum_{j=1}^q d_j \check{B}_n^{(1)}(t_j,v_j),
\end{align*}
where
\begin{align*}
\check{\mathbb{V}}_n(s,t)
	= \frac{1}{\sqrt{n}} \sum_{j=1}^{\lfloor sn \rfloor} (X_{n,j} - \mu^{(j)}) ~~,~~
\check{B}^{(1)}_n(s,t) = \frac{1}{\sqrt{n}} \sum_{i=1}^{\lfloor sn \rfloor} \frac{1}{\sqrt{l}}
\Big ( \sum_{j=i}^{i+l-1} \big(X_{n,j}(t) - \mu^{(j)}(t) \big)\Big  ) \xi_i^{(1)} .
\end{align*}
Using Assumption (A2) and the fact that $(\xi_i^{(1)})_{i\in\mathbb{N}}$
is independent of $(X_{n,j} \colon n\in\mathbb{N}, j=1,\dots,n)$, it can be seen that $\Vert  \hat  Z_n -  Z_n \Vert_2   \lesssim  ~ \frac{1}{\sqrt{n}} + \frac{\sqrt{l}}{\sqrt{n}}
	\to 0 $.
Therefore \eqref{Zn'toZ} follows from
\begin{align} \label{Zn'toZa}
Z_n   = \frac{1}{\sqrt{n}} \sum_{i=1}^n \Big (Z_{i,n} + Z_{i,n}^{(1)} \Big )  \Rightarrow Z,
\end{align}
where (note that $X_{n,j}-\mu^{(j)} = Y_{n,j} - \mu_1$)
\begin{align} \label{zin}
Z_{i,n}
&=\sum_{j=1}^q c_j (Y_{n,i}(u_j)-\mu_1(u_j)) \mathds{1}\{i\leq \lfloor s_j n \rfloor \} , \\
Z_{i,n}^{(1)} &= \xi_i^{(1)} \sum_{j=1}^q d_j \frac{1}{\sqrt{l}}
\Big( \sum_{k=i}^{i+l-1} (Y_{n,k}(v_j)-\mu_1(v_j)) \Big)
\mathds{1}\{i\leq\lfloor t_j n \rfloor \}~. \label{zin1}
\end{align}
A blocking technique will again be utilized to investigate the weak convergence of this sum (see the proof of Theorem \ref{mixingCLT}). To this end, let $ \eta_b$  and $ \eta_s$ denote two constants such
that
\begin{align} \label{etaass}
\frac{\beta(2+\nu)+1}{2+2\nu}<\eta_b < \eta_s <\frac{1}{2},
\end{align}
and  define the length of the small and the big subblocks as
$s_n = \lfloor n^{1/2-\eta_s} \rfloor$ and $b_n = \lfloor n^{1/2-\eta_b} \rfloor$, respectively.
Note that $(\beta(2+\nu)+1)/(2+2\nu) < 1/2$  (since we assumed that  $\beta < 1/3$) and
that we have $k_n = \lfloor n/(b_n+s_n) \rfloor$ blocks in total. Now introducing the sums
\begin{align} \label{bnj}
B_{j,n} = \sum_{i=(j-1)(b_n+s_n)+1}^{(j-1)(b_n+s_n)+b_n} (Z_{i,n}+Z_{i,n}^{(1)})
\quad \text{and} \quad
S_{j,n} = \sum_{i=(j-1)(b_n+s_n)+b_n+1}^{j(b_n+s_n)} (Z_{i,n}+Z_{i,n}^{(1)}) ~,
\end{align}
of the terms in the $j$-th big and small subblock, respectively ($j=1,\dots,k_n$), the representations
\eqref{reprZN}  and \eqref{varZnCLT} in the proof of Theorem \ref{mixingCLT} are obtained, adopting appropriate redefinitions of the quantities $B_{j,n} $, $S_{j,n} $ and $R_n$ involved in this representation.
%\begin{align*}
%Z_n = \frac{1}{\sqrt{n}} \sum_{j=1}^{k_n} B_{j,n} + \frac{1}{\sqrt{n}} \sum_{j=1}^{k_n} S_{j,n}
%+\frac{1}{\sqrt{n}} R_n ~,
%\end{align*}
%where $R_n = \sum_{i=k_n(b_n+s_n)+1}^n (Z_{i,n} + Z_{i,n}^{(1)}) $. Since all random variables
%in this representation  are centered we have
%\begin{align} \label{varZn}
%\begin{split}
%\text{Var}(Z_n) = & ~ \text{Var}\Big(\frac{1}{\sqrt{n}}\sum_{j=1}^{k_n} B_{j,n}\Big)
%+ \frac{2}{n} \sum_{j,j'=1}^{k_n} \mathbb{E}[B_{j,n} S_{j',n}]
%+ \frac{2}{n} \sum_{j=1}^{k_n} \mathbb{E}[B_{j,n} R_n]  \\
%& +\frac{1}{n} \sum_{j,j'=1}^{k_n} \mathbb{E}[S_{j,n} S_{j',n}]
%+ \frac{2}{n} \sum_{j=1}^{k_n} \mathbb{E}[S_{j,n} R_n]
%+ \frac{1}{n}\mathbb{E}[R_n^2] ~.
%\end{split}
%\end{align}
It can be shown that
\begin{align} \label{h1}
|Z_n - 1/\sqrt{n} \sum_{j=1}^{k_n} B_{j,n}| & = |1/\sqrt{n} \sum_{j=1}^{k_n} S_{j,n}
+1/\sqrt{n} R_n| = o_\mathbb{P}(1)  \\
\mbox{Var}(Z_n) &= \text{Var}(\sum_{j=1}^{k_n} B_{j,n}) + o(1) \label{h2},
\end{align}
proving   that each term  on the right-hand side of   \eqref{varZnCLT}, with the exception of the first,
converges to zero in probability. Exemplarily, consider the fourth term to indicate the differences
to the proof of Theorem \ref{mixingCLT}. Note that
\begin{align*}
\mathbb{E}[S_{j,n} S_{j',n}]
= \sum_{i=(j-1)(b_n+s_n)+b_n+1}^{j(b_n+s_n)} \sum_{i'=(j'-1)(b_n+s_n)+b_n+1}^{j'(b_n+s_n)}
\mathbb{E}[Z_{i,n} Z_{i',n}] + \mathbb{E}[Z_{i,n}^{(1)} Z_{i',n}^{(1)}]
\end{align*}
and it follows from (3.17) in Dehling and Philipp \cite{dehling2002} together with Assumption (A1) that
\begin{align} \label{bound1}
|\mathbb{E}[Z_{i,n} Z_{i',n}]|
& \leq \sum_{j,j'=1}^q |c_j c_{j'}| |\text{Cov}(Y_i (u_j), Y_{i'}(u_{j'}))|
%\\
%& \leq 2 \varphi(|i-i'|)^{1/2} \sum_{j,j'=1}^q |c_j c_{j'}| ~
%\mathbb{E}[Y_i(u_j)^2]^{1/2} \mathbb{E}[Y_{i'}(u_{j'})^2]^{1/2} \\
%&
\lesssim \varphi(|i-i'|)^{1/2},
\end{align}
as the $\sigma$-field generated by a $C([0,1])$-valued random variable $X$ always contains the $\sigma$-field generated by
$X(t)$ for a fixed $t\in[0,1]$.   Moreover,  $\mathbb{E} [(Z_{i,n}^{(1)})^2] \leq $ const $ < \infty$, which follows from the representation
\begin{align*}
\mathbb{E} [(Z_{i,n}^{(1)})^2 ]
=  \frac{1}{l}  \sum_{j,j' = 1}^q & d_j d_{j'} \sum_{i',i''=i}^{i+l-1} \text{Cov}(Y_{n,i'}(v_j), Y_{n,i''}(v_{j'}))
\mathds{1}\{ i\leq \lfloor t_j n \rfloor \} \mathds{1}\{ i\leq \lfloor t_{j'} n \rfloor \}
\end{align*}
and the fact that the covariance in the last expression can be estimated by
(here $C$ denotes a constant)
\begin{align*}
 \frac{1}{l} \sum_{i',i''=i}^{i+l-1} |  \text{Cov}(Y_{n,i'}(v_j), Y_{n,i''}(v_{j'})) |  & \leq  C
+ \frac{2}{l} \sum_{i'=i}^{i+l-2} \sum_{i''=i'+1}^{i+l-1}
| \text{Cov}(Y_{n,i'}(v_j), Y_{n,i''}(v_{j'})) |   \\
 \lesssim  C + & \frac{2}{l} \sum_{i'=i}^{i+l-2} \sum_{i''=i'+1}^{i+l-1}
\varphi(|i'-i''|)^{1/2}  \leq C +  \frac{2}{l} \sum_{i'=1}^{l-1} (l-i') a^{i'/2} < \infty,
\end{align*}
using Assumptions (A1), (A4) and (3.17) in Dehling and Philipp \cite{dehling2002}).
Combining this result with \eqref{bound1} leads to
\begin{align} \label{bound2.1}
%\begin{split}
\mathbb{E}[S_{j,n}^2]
%&= \sum_{i=(j-1)(b_n+s_n)+b_n+1}^{j(b_n+s_n)} \mathbb{E}[(Z_{i,n}^{(1)})^2]
%+ \sum_{i,i'=(j-1)(b_n+s_n)+b_n+1}^{j(b_n+s_n)} \mathbb{E}[Z_{i,n} Z_{i',n}] \\
 \lesssim  \sum_{i,i'=(j-1)(b_n+s_n)+b_n+1}^{j(b_n+s_n)} \varphi(|i-i'|)^{1/2} + s_n
 \leq 2 \sum_{i=0}^{s_n-1} (s_n -i) \varphi(i)^{1/2}  + s_n
% \leq 2 s_n \sum_{i=0}^{\infty} a^{i/2}
= O(s_n) ~.
%\end{split}
\end{align}
 With similar arguments
it follows that $\mathbb{E}[S_{j,n}S_{j',n}] = O(s_n^2 \varphi(b_n)^{1/2})$ for $j\neq j'$
since there is at least one big subblock between the observations.
Hence,
\begin{align*}
\frac{1}{n} \sum_{j,j'=1}^{k_n} \mathbb{E}[S_{j,n}S_{j',n}]
= O(n^{-1} k_n s_n ) + O(n^{-1} k_n^2 s_n^2 \varphi(b_n)^{1/2})
= O(b_n^{-1} s_n) + O(n b_n^{-2} s_n^2 a^{b_n / 2}) =o(1)
\end{align*}
 as $b_n^{-1} s_n \to 0$ and $n a^{b_n / 2} \to 0$. It can be shown by similar arguments that the
 second, third, fifth and sixth term in \eqref{varZnCLT} are of order $o(1)$ and \eqref{h1} and \eqref{h2}
 therefore follow.

Equation \eqref{h1} implies  that it suffices to show that $n^{-1/2} \sum_{j=1}^{k_n} B_{j,n}$ converges in
distribution to $Z$ to prove \eqref{Zn'toZa}. Using similar arguments as in the proof of Theorem \ref{mixingCLT},
it can be shown that
\begin{align*}
\Big| \mathbb{E} \Big[~ \prod_{j=1}^{k_n} \psi_{j,n}(t) \Big]
- \prod_{j=1}^{k_n} \mathbb{E} \left[\psi_{j,n}(t)\right] \Big|
\lesssim k_n \varphi(s_n) = O(k_n a^{s_n})  = o(1) ~,
\end{align*}
where $\mathbb{E} [~ \prod_{j=1}^{k_n} \psi_{j,n}(t)]$ is the characteristic
function of $n^{-1/2} \sum_{j=1}^{k_n} B_{j,n}$,
$ \prod_{j=1}^{k_n} \mathbb{E} \left[\psi_{j,n}(t)\right]$
 is the characteristic
function  of  $n^{-1/2} \sum_{j=1}^{k_n} B'_{j,n}$  and
 $B'_{1,n},\dots,B'_{k_n,n}$ are  independent random variables such that
$B_{j,n}$ and $B'_{j,n}$ have the same distribution for any $j=1,\ldots k_n$.
Therefore \eqref{Zn'toZa} follows from
\begin{align}\label{h3}
n^{-1/2} \sum_{j=1}^{k_n} B'_{j,n}   \Rightarrow
Z,
\end{align}
which can be established by the  Lindeberg--Feller central limit theorem for triangular arrays.
Similar arguments as given in the discussion following  \eqref{varasy} give
\begin{align}\label{h3a}
\mbox{Var} \bigg (\frac{1}{\sqrt{n}} \sum_{j=1}^{k_n} B'_{j,n} \bigg )
= \mbox{Var}\bigg  (\frac{1}{\sqrt{n}} \sum_{j=1}^{k_n} B_{j,n}\bigg ) + o(1) = \mbox{Var}(Z_n) + o(1).
\end{align}
We now show that Var$(Z_n)$ converges to the  variance of the random variable $Z$ defined in  \eqref{Zn'toZ}, which is given by
\begin{align} \label{varZ}
\begin{split}
\text{Var}(Z)
%&= \text{Var}\Big(\sum_{j=1}^q c_j \mathbb{V}(s_j,u_j)
%+ \sum_{j=1}^q d_j \mathbb{V}^{(1)}(t_j,v_j)\Big) \\
&= \text{Var}\Big(\sum_{j=1}^q c_j \mathbb{V}(s_j,u_j)\Big)
+ \text{Var}\Big(\sum_{j=1}^q d_j \mathbb{V}^{(1)}(t_j,v_j)\Big) \\
%&= \sum_{j,j'=1}^q c_j c_{j'} \mathrm{Cov}(\mathbb{V}(s_j,u_j),\mathbb{V}(s_{j'},u_{j'}))
%+ \sum_{j,j'=1}^q d_j d_{j'} \mathrm{Cov}(\mathbb{V}^{(1)}(t_j,v_j),\mathbb{V}^{(1)}(t_{j'},v_{j'})) \\
&= \sum_{j,j'=1}^q c_j c_{j'} (s_j\wedge s_{j'})  C(u_j,u_{j'})
+ \sum_{j,j'=1}^q d_j d_{j'} (t_j\wedge t_{j'})  C(v_j,v_{j'}) ~.
\end{split}
\end{align}
On the other hand (observing that
$\mathrm{Cov}(Z_{i,n},Z^{(1)}_{i',n}) = 0$),
\begin{align} \label{h4}
\text{Var}(Z_n)
%&= \text{Var} \Big(\frac{1}{\sqrt{n}} \sum_{i=1}^n (Z_{i,n} + Z^{(1)}_{i,n}) \Big) \\
%&= \frac{1}{n} \sum_{i,i'=1}^n \mathrm{Cov}(Z_{i,n} + Z^{(1)}_{i,n},Z_{i',n} + Z^{(1)}_{i',n}) \\
&=  \frac{1}{n} \sum_{i,i'=1}^n \Big(\mathbb{E}[Z_{i,n} Z_{i',n}]
	+ \mathbb{E}[Z^{(1)}_{i,n} Z^{(1)}_{i',n}]\Big)  =  \frac{1}{n} \sum_{i,i'=1}^n \mathbb{E}[Z_{i,n} Z_{i',n}]
	+   \frac{1}{n} \sum_{i=1}^n \mathbb{E}[(Z^{(1)}_{i,n})^2]
~.
\end{align}
The second term in this expression satisfies
\begin{align*}
\frac{1}{n} \sum_{i=1}^n \mathbb{E}[(Z^{(1)}_{i,n})^2]
%=&  \frac{1}{n} \sum_{i=1}^n \sum_{j,j'=1}^q d_j d_{j'}
%  \mathbb{E} \Big[\frac{1}{l} \sum_{i'_1,i'_2=i}^{i+l-1} (Y_{i'_1}(v_j) - \mu_1(v_j))
% (Y_{i'_2}(v_{j'}) - \mu_1(v_{j'}))
%   \mathds{1}\big\{i \leq \lfloor (t_j \wedge t_{j'})n \rfloor \big\}
%  \Big] \\
=& \frac{1}{n} \sum_{i=1}^n \sum_{j,j'=1}^q d_j d_{j'}
 \frac{1}{l} \sum_{i'_1,i'_2=i}^{i+l-1} \mathrm{Cov}(Y_{i'_1}(v_j),Y_{i'_2}(v_{j'}))
  \mathds{1}\big\{i \leq \lfloor (t_j \wedge t_{j'})n \rfloor \big\} \\
% =& \sum_{j,j'=1}^q d_j d_{j'} \frac{1}{n}
%   \sum_{i=1}^{\lfloor (t_j \wedge  t_{j'})n \rfloor}
% \frac{1}{l} \sum_{i',i''=i}^{i+l-1} \gamma(i'-i'',v_j,v_{j'}) \\
=& \sum_{j,j'=1}^q d_j d_{j'}
   \frac{\lfloor (t_j \wedge  t_{j'})n \rfloor}{n} ~
 \sum_{i'=-(l-1)}^{l-1} \frac{l-|i'|}{l} ~ \gamma(i',v_j,v_{j'}) ~ \\
&\to \sum_{j,j'=1}^q d_j d_{j'}
  (t_j \wedge  t_{j'}) \sum_{i=-\infty}^{\infty}\gamma(i,v_j,v_{j'})
= \sum_{j,j'=1}^q d_j d_{j'}
  (t_j \wedge  t_{j'}) C(v_j,v_{j'}) ,
\end{align*}
where the dominated convergence theorem was used in the last step. For the first term in \eqref{h4}, assume without loss of generality that $s_j \leq s_{j'}$  and note that
\begin{align*}
\frac{1}{n} \sum_{i,i'=1}^n \mathbb{E}[Z_{i,n} Z_{i',n}]
% &= \frac{1}{n} \sum_{j,j'=1}^q c_j c_{j'}  \sum_{i,i'=1}^n
% \mathrm{Cov} (Y_i(u_j),Y_{i'}(u_{j'}))
%  \mathds{1}\{i\leq \lfloor s_j n \rfloor \} \mathds{1}\{i'\leq \lfloor s_{j'} n \rfloor \} \\
&= \frac{1}{n} \sum_{j,j'=1}^q c_j c_{j'}  \sum_{i=1}^{\lfloor s_j n \rfloor}
 \sum_{i'= 1}^{\lfloor s_{j'} n \rfloor}
 \mathrm{Cov} (Y_i(u_j),Y_{i'}(u_{j'})) \\
=& \sum_{j,j'=1}^q c_j c_{j'} \bigg( \frac{1}{n} \sum_{i=1}^{\lfloor s_j n \rfloor} \sum_{i'= 1}^{\lfloor s_{j} n \rfloor}
 \gamma(i-i',u_j,u_{j'})
+ \frac{1}{n} \sum_{i=1}^{\lfloor s_j n \rfloor}
\sum_{i'= \lfloor s_j n \rfloor +1}^{\lfloor s_{j'} n \rfloor}
 \mathrm{Cov} (Y_i(u_j),Y_{i'}(u_{j'})) \bigg) \\
 &=  \sum_{j,j'=1}^q c_j c_{j'} s_j \sum_{i=-\infty}^\infty \gamma(i,u_j,u_{j'}) + o(1)
= \sum_{j,j'=1}^q c_j c_{j'} s_j ~ C(u_j,u_{j'}) + o(1),
\end{align*}
where the dominated convergence theorem was used again for the first term
and the  bound
\begin{align*}
\frac{1}{n} \sum_{i=1}^{\lfloor s_j n \rfloor}
\sum_{i'= \lfloor s_j n \rfloor +1}^{\lfloor s_{j'} n \rfloor}
 |\mathrm{Cov} (Y_i(u_j),Y_{i'}(u_{j'}))|
 \lesssim
 % \frac{1}{n} \sum_{i=1}^{\lfloor s_j n \rfloor}
%\sum_{i'= \lfloor s_j n \rfloor +1}^{\lfloor s_{j'} n \rfloor} \varphi(|i-i'|)^{1/2} \\ \leq
 \frac{1}{n} \sum_{i=1}^{\lfloor s_{j'} n \rfloor -1}
  i \varphi(i)^{1/2}
\leq \frac{1}{n} \sum_{i=1}^{\lfloor s_{j'} n \rfloor -1}
  i a^{i/2} = o(1)
\end{align*}
for  the second term. Observing \eqref{h3a}--\eqref{h4}, it follows that Var$(n^{-1/2} \sum_{j=1}^{k_n} B'_{j,n} ) =$Var$(Z) +o(1)  $ as postulated.

For a proof of the Lindeberg condition, use H\"older's inequality (with $p=1+\nu /2$ and $q=(2+\nu)/\nu$,
where $\nu$ is the same as in Assumption (A1))  and Markov's
inequality to obtain
\begin{align}  \nonumber
\Delta_n &= \frac{1}{n} \sum_{j=1}^{k_n} \mathbb{E}\big[(B_{j,n}')^2
	\mathds{1} \{|B_{j,n}'| > \sqrt{n}\delta\} \big]  =
\frac{1}{n}  \sum_{j=1}^{k_n} \mathbb{E}\big[B_{j,n}^2
	\mathds{1} \{|B_{j,n}| > \sqrt{n}\delta\} \big]  \\
%	\\
 &\leq
\frac{1}{n} \sum_{j=1}^{k_n} \mathbb{E} \big[ |B_{j,n}|^{2+\nu}\big]^{2/(2+\nu)}
	\mathbb{P}\big(|B_{j,n}| > \sqrt{n}\delta\big)^{\nu/(2+\nu)}
	%\\ &=
 %\frac{1}{n} \sum_{j=1}^{k_n} \mathbb{E} \big[ |B_{j,n}|^{2+\nu}\big]^{2/(2+\nu)}
%	\mathbb{P}\big(|B_{j,n}|^{2+\nu} > n^{(2+\nu)/2}\delta^{2+\nu}\big)^{\nu/(2+\nu)} \\
%&
%& \leq
%\frac{1}{n} \sum_{j=1}^{k_n} \mathbb{E} \big[ |B_{j,n}|^{2+\nu}\big]^{2/(2+\nu)}
%	\mathbb{E}\big[|B_{j,n}|^{2+\nu}\big]^{\nu/(2+\nu)} (n^{1/2}\delta)^{-\nu} \\ &=
\leq \frac{1}{n} \sum_{j=1}^{k_n} \mathbb{E} \big[ |B_{j,n}|^{2+\nu}\big] (n^{1/2}\delta)^{-\nu} .
\label{linde}
\end{align}
Now, observing the definition  of   $Z_{i,n}$,  $Z_{i,n}^{(1)}$ in \eqref{zin}, \eqref{zin1}, respectively,  and  Assumption (A1),  it can be seen that  $\max_{1\leq i\leq n} \mathbb{E}\big[|Z_{i,n}|^{2+\nu}\big]^{1/(2+\nu)}
 < \infty$ and
$
 \max_{1\leq i\leq n} \mathbb{E}\big[|Z^{(1)}_{i,n}|^{2+\nu}\big]^{1/(2+\nu)}
 \lesssim \sqrt{l} $,  and the Lindeberg condition  follows from \eqref{linde}
 observing that  the representation  of $B_{j,n} $ in \eqref{bnj} and   Minkowski's inequality give
$\mathbb{E} \big[ |B_{j,n}|^{2+\nu}\big]^{1/(2+\nu)}=O(b_n l^{1/2})$, that is,
\begin{align*}
\Delta_n
= O \Big ( \frac{k_n b_n^{2+\nu} l^{(2+\nu)/2}}{ n^{1+ \nu/2}} \Big)= O \Big ( \frac{ b_n^{1+\nu} l^{(2+\nu)/2} }{n^{\nu/2}} \Big)
&=  O(n^{1/2-\eta_b(1+\nu)} l^{(2+\nu)/2})
= O(n^{1/2-\eta_b(1+\nu)} n^{\beta(2+\nu)/2}) ,
\end{align*}
which   converges to zero by the assumption  \eqref{etaass}.
\medskip

\noindent
{\it  (B) Asymptotic tightness:}   Since $(\hat{\mathbb{V}}_n)_{n\in\mathbb{N}}$ converges weakly  to $\mathbb{V}$, the process
 $(\hat{\mathbb{V}}_n)_{n\in\mathbb{N}}$ is  asymptotically tight and it remains to
 show that $(\tilde{B}_n^{(k)})_{n\in\mathbb{N}}$ is
asymptotically tight (for any $k=1,\dots,R$, note that
marginal asymptotic tightness implies joint asymptotic tightness).

Let $s,t\in[0,1]$ be arbitrary and define $\varepsilon_{n,j} = Y_{n,j} -\mu_1$ for
any $n\in\mathbb{N}$ and $j=1,\dots,n$, then, since $\xi_1^{(k)},\dots,\xi_n^{(k)}$ are independent of $\varepsilon_{n,1},\dots,\varepsilon_{n,n}$
with $\mathbb{E} [\xi_j^{(k)}] = 0$ and Var$(\xi_j^{(k)}) = 1$,
\begin{align} \nonumber
\Vert\tilde{B}^{(k)}_n (1,s) - \tilde{B}^{(k)}_n (1,t)\Vert_2^2 =&
% \mathbb{E}\Big[ \Big\vert \frac{1}{\sqrt{n}} \sum_{i=1}^n \frac{1}{\sqrt{l}}
%	\Big\{ \sum_{j=i}^{i+l-1} \big(Y_{n,j}(s) - \mu_1(s) -Y_{n,j}(t)+\mu_1(t) \big) \Big\}
%	\xi_i^{(k)} \Big\vert^2 \Big] \\
%=& \frac{1}{n} \sum_{i,i'=1}^n \mathbb{E} \Big[ \frac{1}{l}
%	\Big\{\sum_{j=i}^{i+l-1} \big(Y_{n,j}(s) - \mu_1(s) -Y_{n,j}(t)+\mu_1(t) \big) \Big\} \\
%	&\hspace{50pt}\times \Big\{\sum_{j'=i'}^{i'+l-1} \big(Y_{n,j'}(s) - \mu_1(s) -Y_{n,j'}(t)
%			+\mu_1(t) \big) \Big\} \xi_i^{(k)} \xi_{i'}^{(k)} \Big] \\
%=& \frac{1}{n} \sum_{i,i'=1}^n \mathbb{E} \Big[ \frac{1}{l}
%	\Big\{\sum_{j=i}^{i+l-1} \big(\varepsilon_{n,j}(s) - \varepsilon_{n,j}(t) \big) \Big\}
%	 \Big\{\sum_{j'=i'}^{i'+l-1}
%			\big(\varepsilon_{n,j'}(s) - \varepsilon_{n,j'}(t) \big) \Big\}
%			\xi_i^{(k)} \xi_{i'}^{(k)} \Big] ~.
%\end{align*}
%
%\begin{align}
%\begin{split}
\frac{1}{n}  \sum_{i=1}^n \mathbb{E} \Big[ \frac{1}{l}
	\sum_{j=i}^{i+l-1} \sum_{j'=i}^{i+l-1}
	\big(\varepsilon_{n,j}(s) - \varepsilon_{n,j}(t) \big)
	\big(\varepsilon_{n,j'}(s) - \varepsilon_{n,j'}(t) \big) \Big] \\
% \leq &
% \frac{1}{n} \sum_{i=1}^n \frac{2}{l} \sum_{j=0}^{l-1}
%	(l-j) \mathbb{E}\big[ \big(\varepsilon_{n,1}(s) - \varepsilon_{n,1}(t) \big)
%	\big(\varepsilon_{n,1+j}(s) - \varepsilon_{n,1+j}(t) \big) \big] \\
 \leq&  \frac{2}{l} \sum_{j=0}^{l-1}
	(l-j) \mathbb{E}\big[ \big(\varepsilon_{n,1}(s) - \varepsilon_{n,1}(t) \big)
	\big(\varepsilon_{n,1+j}(s) - \varepsilon_{n,1+j}(t) \big) \big] ,
%\end{split}
 \label{tight0}
\end{align}
utilizing the fact that each row of the array
$(\varepsilon_{n,j} \colon n\in\mathbb{N}, j=1,\dots,n)$ is stationary.
Assumption (A3) implies
%\begin{align*}
%| & \varepsilon_{n,i}(s) - \varepsilon_{n,i}(t) |^2 \\
%&= | Y_{n,i}(s) - \mu_1(s) - Y_{n,i}(t) + \mu_2(t)) |^2 \\
%&= | X_{n,i}(s) - X_{n,i}(t) + (\mu^{(i)}(s) - \mu^{(i)}(t)) |^2 \\
%&\leq | X_{n,i}(s) - X_{n,i}(t)|^2
%	+ 2| X_{n,i}(s) - X_{n,i}(t)| ~ |\mu^{(i)}(s) - \mu^{(i)}(t) |
%	+ |\mu^{(i)}(s) - \mu^{(i)}(t) |^2 \\
%&\leq M^2 | s-t|^2 + 2 M \mathbb{E}[M] ~ | s-t |^2 + \mathbb{E}[M]^2 ~ | s-t |^2 \\
%\end{align*}
%and therefore we have  \begin{align*}
$\mathbb{E} \big[ | \varepsilon_{n,i}(s) - \varepsilon_{n,i}(t)|^2 \big]^{1/2}
	\lesssim | s-t |,$
%\end{align*}
and (3.17) in Dehling and Philipp \cite{dehling2002} yields
\begin{align*}
\Vert\tilde{B}^{(k)}_n (1,s) - \tilde{B}^{(k)}_n (1,t)\Vert_2^2
&\lesssim \frac{2}{l} \sum_{j=0}^{l-1}
	(l-j) |s-t|^2 ~ a^{j/2} \lesssim |s-t|^2 \sum_{j=0}^{\infty} a^{j/2}
	\lesssim |s-t|^2 ~.
\end{align*}
Now consider the metric  $\rho(s,t) = |s-t|$ on the interval $[0,1]$ and define
$D(\eta,\rho) = \big\lceil \frac{1}{\eta} \big\rceil $
as the corresponding  packing number, then Theorem 2.2.4 in  Van der Vaart and Wellner
\cite{wellner1996} shows that
\begin{align*}
\Big \Vert \sup_{\rho(s,t)\leq \delta} \vert \tilde{B}^{(k)}_n (1,s)
	- \tilde{B}^{(k)}_n (1,t) \vert \Big\Vert_2
&\leq K' \Big [ \int_0^\eta \sqrt{D(\nu,\rho)}~d\nu + \delta D(\eta,\rho) \Big  ]
\lesssim
% \Big[ \int_0^\eta \frac{1}{\sqrt{\nu}} ~d\nu
%	+  \frac{\delta}{\eta} \Big] \\
%&=
2 \sqrt{\eta} +  \frac{\delta}{\eta} ,
\end{align*}
and Markov's inequality now yields, for any $\varepsilon>0$,
\begin{align*}
\mathbb{P} \Big  ( \sup_{\rho(s,t)\leq \delta} \big\vert \tilde{B}^{(k)}_n (1,s)
	- \tilde{B}^{(k)}_n (1,t) \big\vert > \varepsilon \Big )
%&\leq \frac{1}{\varepsilon^2} \Big\Vert \sup_{\rho(s,t)\leq \delta}
%	\vert \tilde{B}^{(k)}_n (1,s) - \tilde{B}^{(k)}_n (1,t) \vert \Big\Vert_2^2 \\
&\lesssim \frac{1}{\varepsilon^2}
	\Big[   2 \sqrt{\eta} + \frac{\delta}{\eta}\Big ].
\end{align*}
Since $\eta >0$ is arbitrary, it follows that
%\begin{align*}
%\lim_{\delta\searrow 0} ~ \limsup_{n\to\infty} ~ \mathbb{P} \Big (
%	\sup_{\rho(s,t)\leq \delta} \big\vert \tilde{B}^{(k)}_n (1,s)
%	- \tilde{B}^{(k)}_n (1,t) \big\vert > \varepsilon \Big  ) = 0 ~.
%\end{align*} Therefore
the process $(\tilde{B}^{(k)}_n(1,\cdot))_{n\in\mathbb{N}}$ is asymptotically uniformly
$\rho$-equicontinu\-ous in probability. Moreover, the finite-dimensional distributions
of $(\tilde{B}^{(k)}_n(1,\cdot))_{n\in\mathbb{N}}$ converge weakly to the finite-dimensional
distribution of $\mathbb{V}^{(k)} (1,\cdot)$ and therefore it follows that
 $\tilde{B}^{(k)}_n(1,\cdot) \rightsquigarrow \mathbb{V}^{(k)} (1,\cdot)$
in $C([0,1])$ (see the discussion at the beginning of the proof of Theorem \ref{mixingCLT}).

The asymptotic tightness of the process $\tilde{B}^{(k)}_n$ in $C([0,1]^2)$ is now a consequence of Corollary 3.5 in Samur \cite{samur1987}.
To be precise, note that $\tilde{B}^{(k)}_n(1,t)=n^{-1/2}\sum_{i=1}^n D_{n,i}^{(k)}(t)$,
where
\begin{align*}
D_{n,i}^{(k)}(t) = \frac{1}{\sqrt{l}} \Big (\sum_{j=i}^{i+l-1}
	\big(Y_{n,j}(t) - \mu_1(t)\big)\Big )\xi_i^{(k)}
= \frac{1}{\sqrt{l}} \Big (\sum_{j=i}^{i+l-1}
	\big(X_{n,j}(t) - \mu^{(j)}(t)\big)\Big )\xi_i^{(k)}.
\end{align*}
The array  $(D^{(k)}_{n,i} \colon n\in\mathbb{N}, i=1,\dots,n)$ is $\varphi$-mixing and
(by Assumption (A1) and Markov's inequality)
\begin{align*}
& n \mathbb{P} \Big( \frac{1}{\sqrt{n}} \|D^{(k)}_{n,1}\| > \varepsilon\Big)
% \leq \frac{1}{\varepsilon^4} \frac{1}{n} \mathbb{E} \big[ \|D^{(k)}_{n,1}\|^4 \big]
%= \frac{1}{\varepsilon^4} \frac{1}{n} \mathbb{E} \Big[ \Big\| \frac{1}{\sqrt{l}}
%	\Big( \sum_{j=1}^m (Y_{n,j}-\mu_1)\Big) \xi_1^{(k)} \Big\|^4 \Big] \\
%&\hspace{15pt}
\lesssim \frac{1}{\varepsilon^4} \frac{1}{n} \frac{1}{l^2}
	\mathbb{E} \Big[ \Big( \sum_{j=1}^l \|X_{n,j}-\mu^{(j)}\|\Big) ^4 \Big]
\lesssim \frac{1}{\varepsilon^4} \frac{l^2}{n}
= \frac{1}{\varepsilon^4} n^{2\beta -1} .
\end{align*}
Therefore, since $\beta < 1/3$ by assumption, $\lim_{n\to\infty} n
	\mathbb{P} ( n^{-1/2} \|D^{(k)}_{n,1}\| > \varepsilon) = 0 $.
By the previous discussion and \eqref{bWIP}, use Corollary 3.5 in Samur \cite{samur1987} to
obtain $(\tilde{B}^{(k)}_n)_{n\in\mathbb{N}} \rightsquigarrow \mathbb{V}$ in $C([0,1]^2)$, which  finally implies that
$(\tilde{B}^{(k)}_n)_{n\in\mathbb{N}}$ is asymptotically tight.
\end{proof}

\begin{proof}[\bf Proof of \eqref{bLemma2}]
Note that
\begin{align*}
\Vert  \bar{B}^{(k)}_n - \tilde{B}^{(k)}_n \Vert
\leq & \Big\Vert \frac{1}{\sqrt{n}} \sum_{i=1}^{\lfloor sn \rfloor} \frac{1}{\sqrt{l}}
	\Big( l\mu_1(t) - \frac{l}{n} \sum_{j=1}^n Y_{n,j}(t) \Big) \xi_i^{(k)} \Big\Vert \\
	&+ \sup_{s,t\in[0,1]} \Big\vert \sqrt{n}\Big(s-\frac{\lfloor sn \rfloor}{n} \Big)
	\frac{1}{\sqrt{l}}\Big( l\mu_1(t) - \frac{l}{n} \sum_{j=1}^n Y_{n,j}(t)\Big)
	\xi_{\lfloor sn \rfloor +1}^{(k)} \Big\vert  \\
&\leq \sup_{s\in [0,1]} \Big\vert \frac{1}{\sqrt{n}} \sum_{i=1}^{\lfloor sn \rfloor}
	\xi_i^{(k)} \Big\vert
	\times \frac{\sqrt{l}}{\sqrt{n}}\sup_{t\in[0,1]} \Big\vert
	\frac{1}{\sqrt{n}} \sum_{j=1}^n (X_{n,j}(t)-\mu^{(j)}(t)) \Big\vert \\
		&+ \frac{1}{\sqrt{n}} \max_{i=1}^n \vert \xi_i^{(k)} \vert
	\times \frac{\sqrt{l}}{\sqrt{n}} \sup_{t\in[0,1]}\Big\vert
	\frac{1}{\sqrt{n}} \sum_{j=1}^n \big(X_{n,j}(t) - \mu^{(j)}(t) \big) \Big\vert.
\end{align*}
Therein, $\sup_{s\in [0,1]} \vert n^{-1/2} \sum_{i=1}^{\lfloor sn \rfloor}
	\xi_i^{(k)} \vert = O_\mathbb{P} (1)$, $\max_{i=1}^n \vert \xi_i^{(k)} \vert = O_\mathbb{P} (\sqrt{\log n})$,
	while all other terms are of order $ O_\mathbb{P} ( \sqrt{l/n}) = o_\mathbb{P}(1)$, by Theorem \ref{mixingCLT}.
\end{proof}

\begin{proof}[\bf Proof of  \eqref{bLemma3}]
Recall the definition of  $ \bar{B}_n^{(k)}(s,t)$ in \eqref{Bbar} and define
\begin{align*}
\check{B}_n^{(k)}(s,t) =& \frac{1}{\sqrt{n}} \sum_{i=1}^{\lfloor sn \rfloor}\frac{1}{\sqrt{l}}
	\Big( \sum_{j=i}^{i+l-1} \check{Y}_{n,j}(t)
	- \frac{l}{n} \sum_{j=1}^n \check{Y}_{n,j}(t) \Big) \xi_i^{(k)} \\
&+ \sqrt{n}\Big(s - \frac{\lfloor sn \rfloor}{n} \Big)\frac{1}{\sqrt{l}}
\Big( \sum_{j=\lfloor sn \rfloor +1}^{\lfloor sn \rfloor+l} \check{Y}_{n,j}(t)
- \frac{l}{n} \sum_{j=1}^n \check{Y}_{n,j}(t) \Big) \xi_{\lfloor sn \rfloor +1}^{(k)},
\end{align*}
where $\check{Y}_{n,j} = X_{n,j}-(\hat{\mu}_2 - \hat{\mu}_1) \mathds{1}\{j> \lfloor s^*n \rfloor \}$. Then,
\begin{align} \label{439}
 \sup_{(s,t) \in [0,1]^2} |\check{B}_n^{(k)} (s,t) - \bar{B}_n^{(k)}(s,t)|  \leq
U_n^{(1,1)} + U_n^{(1,2)} + U_n^{(2)}  = o_{\mathbb{P}} (1),
\end{align}
where
\begin{align} \label{440}
 U_n^{(1,1)} & =  \sup_{(s,t) \in [0,1]^2} \Big|
	\frac{1}{\sqrt{n}} \sum_{i=1}^{\lfloor sn \rfloor} \frac{1}{\sqrt{l}}
	\Big( \sum_{j=i}^{i+l-1} \varepsilon_n(t) \mathds{1}\{j> \lfloor s^*n \rfloor \}
	\Big) \xi_i^{(k)} \Big|, \\
 U_n^{(1,2)} & = \sup_{(s,t) \in [0,1]^2} \Big|
	\frac{1}{\sqrt{n}} \sum_{i=1}^{\lfloor sn \rfloor} \frac{1}{\sqrt{l}} \Big(
	\frac{l}{n} \sum_{j=1}^n \varepsilon_n(t)\mathds{1}\{j> \lfloor s^*n \rfloor \}
	\Big) \xi_i^{(k)} \Big|, \label{440a} \\
\nonumber
U_n^{(2)}&= \sup_{(s,t) \in [0,1]^2} \Big|
\frac{\sqrt{n}}{\sqrt{l}} \Big(s - \frac{\lfloor sn \rfloor}{n} \Big)
	\Big( \sum_{j=\lfloor sn \rfloor +1}^{\lfloor sn \rfloor+l} \varepsilon_n(t)
	\mathds{1}\{j> \lfloor s^*n \rfloor \}
	- \frac{l}{n} \sum_{j=1}^n \varepsilon_n(t)\mathds{1}\{j> \lfloor s^*n \rfloor \}
	\Big) \xi_{\lfloor sn \rfloor +1}^{(k)} \Big|,
\end{align}
and  $\varepsilon_n = \hat{\mu}_2 - \mu_2 - (\hat{\mu}_1 - \mu_1)$. To prove \eqref{439}, it will be shown
that all  terms  on the right-hand side of \eqref{439} converge to zero in probability, concentrating on $U_n^{(1,1)}$  and $U_n^{(1,2)}$
for the sake of brevity
(the term  $U_n^{(2)}$ can be treated similarly).  At the end of this proof it will be verified that
\begin{align} \label{442}
\sup_{t\in[0,1]} \big| \varepsilon_n(t) \big| = O_{\mathbb{P}}\bigg(\frac{1}{\sqrt{n}}\bigg).
\end{align}
Direct calculations observing \eqref{442} yield
\begin{align} \label{441}
 U_n^{(1,2)}
 & \leq \sup_{s \in [0,1]} \Big| \frac{1}{\sqrt{n}}
\sum_{i=1}^{\lfloor sn \rfloor} \xi_i^{(k)} \Big| ~
\times \sqrt{l} \sup_{t\in[0,1]} \big| \varepsilon_n(t) \big|
= O_\mathbb{P} \bigg ( \sqrt{l\over n} \bigg)
= o_\mathbb{P} (1).
\end{align}
For the first term on the right side of \eqref{439} similar
arguments  yield
\begin{align*}
 U_n^{(1,1)} & = \sup_{t \in [0,1]} \Big|
	\frac{1}{\sqrt{n}} \sum_{i=1}^{\lfloor s^*n \rfloor} \frac{1}{\sqrt{l}}
	\Big( \sum_{j=i}^{i+l-1} \varepsilon_n(t) \mathds{1}\{j> \lfloor s^*n \rfloor \}
	\Big) \xi_i^{(k)} \Big| + o_{\mathbb{P}}(1) \\
&\leq
	 \sum_{i=\lfloor s^*n \rfloor -l +1}^{\lfloor s^*n \rfloor} |\xi_i^{(k)}|
	\frac{\sqrt{l}}{\sqrt{n}} \sup_{t \in [0,1]} \big|\varepsilon_n(t)\big| + o_{\mathbb{P}}(1)
	= O_{\mathbb{P}} \Big ( { l^{3/2} \over  n} \Big ) +o_{\mathbb{P}}(1)
	= O_{\mathbb{P}} ( n^{3/2\beta-1})+o_{\mathbb{P}}(1) = o_{\mathbb{P}}(1)
	,
\end{align*}
which follows using \eqref{442},  $\beta < 1/3$ (by assumption) and $\sum_{i=\lfloor s^*n \rfloor -l +1}^{\lfloor s^*n \rfloor} |\xi_i^{(k)}| = O_{\mathbb{P}}(l)$, the latter relation implied by Markov's inequality.

Therefore \eqref{439} holds and observing the definition of $\hat{B}_n^{(k)}$ in \eqref{bProcess}, it is next shown that
\begin{align} \label{443}
& \sup_{(s,t) \in [0,1]^2} |\hat{B}_n^{(k)} (s,t) - \check{B}_n^{(k)}(s,t)|  \leq  Z_n^{(1)} +Z_n^{(2)} +Z_n^{(3)} +Z_n^{(4)}   = o_{\mathbb{P}}(1),
\end{align}
where
\begin{align*}
Z_n^{(1)}  &  = \sup_{(s,t) \in [0,1]^2} \Big|
	\frac{1}{\sqrt{n}} \sum_{i=1}^{\lfloor sn \rfloor} \frac{1}{\sqrt{l}}
	 \sum_{j=i}^{i+l-1} \big(\hat{\mu}_1(t)-\hat{\mu}_2(t)\big)
	\big(\mathds{1}\{j> \lfloor s^*n \rfloor \}
	- \mathds{1}\{j> \lfloor \hat{s}n \rfloor \} \big)  \xi_i^{(k)}  \Big| \\
Z_n^{(2)}  &  = \sup_{(s,t) \in [0,1]^2} \Big|
	\frac{1}{\sqrt{n}} \sum_{i=1}^{\lfloor sn \rfloor} \frac{1}{\sqrt{l}}   \frac{l}{n} \sum_{j=1}^n \big(\hat{\mu}_1(t)-\hat{\mu}_2(t)\big)
	\big(\mathds{1}\{j> \lfloor s^*n \rfloor \}
	- \mathds{1}\{j> \lfloor \hat{s}n \rfloor \} \big)
	 \xi_i^{(k)}  \Big|, \\
Z_n^{(3)} & =  \sup_{(s,t) \in [0,1]^2} \Big| \Big(s - \frac{\lfloor sn \rfloor}{n} \Big)\frac{\sqrt{n}}{\sqrt{l}}
        \sum_{j=\lfloor sn \rfloor +1}^{\lfloor sn \rfloor+l}
	\big(\hat{\mu}_1(t)-\hat{\mu}_2(t)\big)
	\big(\mathds{1}\{j> \lfloor s^*n \rfloor \}
	- \mathds{1}\{j> \lfloor \hat{s}n \rfloor \} \big) \xi_{\lfloor sn \rfloor +1}^{(k)}\Big|, \\
Z_n^{(4)}  &=  \sup_{(s,t) \in [0,1]^2} \Big| \Big(s - \frac{\lfloor sn \rfloor}{n} \Big)\frac{\sqrt{n}}{\sqrt{l}}  \frac{l}{n} \sum_{j=1}^n \big(\hat{\mu}_1(t)-\hat{\mu}_2(t)\big)
	\big(\mathds{1}\{j> \lfloor s^*n \rfloor \}
	- \mathds{1}\{j> \lfloor \hat{s}n \rfloor \} \big)
	\Big) \xi_{\lfloor sn \rfloor +1}^{(k)} \Big|.
\end{align*}
As  $\big|\mathds{1}\{j> \lfloor s^*n \rfloor \} - \mathds{1}\{j> \lfloor \hat{s}n \rfloor \} \big|
= \mathds{1}\{\lfloor (\hat{s} \wedge s^*)n \rfloor < j
	\leq \lfloor (\hat{s} \vee s^*)n \rfloor\} $,
it follows that
\begin{align} \label{444}
Z_n^{(2)}&\leq \sup_{(s,t) \in [0,1]^2} \Big(
	\frac{1}{\sqrt{n}} \sum_{i=1}^{\lfloor sn \rfloor} \frac{|\xi_i^{(k)}|}{\sqrt{l}}
	\frac{l}{n} \sum_{j=1}^n \big|\hat{\mu}_1(t)-\hat{\mu}_2(t)\big| ~
	\mathds{1}\{\lfloor (\hat{s} \wedge s^*)n \rfloor < j
	\leq \lfloor (\hat{s} \vee s^*)n \rfloor\}
	  \Big) \\
%&\hspace{15pt} \leq \|\hat{\mu}_1-\hat{\mu}_2\|_\infty
%	\Big(\frac{1}{\sqrt{n}} \sum_{i=1}^{n} |\xi_i^{(k)}| \Big)
%	\frac{\sqrt{l}}{n}
%	\sum_{j=1}^n \mathds{1}\{\lfloor (\hat{s} \wedge s^*)n \rfloor < j
%	\leq \lfloor (\hat{s} \vee s^*)n \rfloor\} \\
%&\hspace{15pt} = \|\hat{\mu}_1-\hat{\mu}_2\|_\infty
%	\Big(\frac{1}{\sqrt{n}} \sum_{i=1}^{n} |\xi_i^{(k)}| \Big)
%	\frac{\sqrt{l}}{n} ~ \Big(\lfloor (\hat{s} \vee s^*)n \rfloor
%	- \lfloor (\hat{s} \wedge s^*)n \rfloor \Big) \\
&\leq 	\|\hat{\mu}_1-\hat{\mu}_2\|_\infty
	\Big(\frac{1}{\sqrt{n}} \sum_{i=1}^{n} |\xi_i^{(k)}| \Big)
	\frac{\sqrt{l}}{n} ~\big(n|\hat{s} - s^*|+1\big)  = O_{\mathbb{P}} \Big (\sqrt{l \over n} \Big ) =    o_{\mathbb{P}}( 1),
	\nonumber
\end{align}
using that $\frac{1}{\sqrt{n}} \sum_{i=1}^{n} |\xi_i^{(k)}| = O_{\mathbb{P}}(\sqrt{n})$ (by Markov's inequality)
and $\|\hat{\mu}_1 - \hat{\mu}_2 \|  = O_{\mathbb{P}}(1)$ (implied by \eqref{442}).

For the first term on the right-side of \eqref{443}, use similar arguments as in \eqref{444} to obtain
\begin{align*}
	Z_n^{(1)}  & \leq \big\|\hat{\mu}_1-\hat{\mu}_2\big\|_\infty
	\frac{1}{\sqrt{n}} \sum_{i=1}^{n} \frac{1}{\sqrt{l}}
	\Big( \sum_{j=i}^{i+l-1}  ~
	\mathds{1}\{\lfloor (\hat{s} \wedge s^*)n \rfloor < j
	\leq \lfloor (\hat{s} \vee s^*)n \rfloor\} \Big) |\xi_i^{(k)}| \\
&\hspace{15pt} = \big\|\hat{\mu}_1-\hat{\mu}_2\big\|_\infty
	\frac{1}{\sqrt{n}}
	\sum_{i=\lfloor (\hat{s} \wedge s^*)n \rfloor -l+2}^{\lfloor (\hat{s} \vee s^*)n \rfloor}
	\frac{1}{\sqrt{l}}
	\Big( \sum_{j=i}^{i+l-1}  ~
	\mathds{1}\{\lfloor (\hat{s} \wedge s^*)n \rfloor < j
	\leq \lfloor (\hat{s} \vee s^*)n \rfloor\} \Big) |\xi_i^{(k)}| \\
&\hspace{15pt} \leq \big\|\hat{\mu}_1-\hat{\mu}_2\big\|_\infty
	\frac{\sqrt{l}}{\sqrt{n}}
	\sum_{i=\lfloor (\hat{s} \wedge s^*)n \rfloor -l+2}^{\lfloor (\hat{s} \vee s^*)n \rfloor}
	|\xi_i^{(k)}| = O_{\mathbb{P}}(l^{3/2}n^{-1/2}) = O_{\mathbb{P}}(n^{ 3 \beta /2-1/2}) = o_{\mathbb{P}} (1),
\end{align*}
since $\sum_{i=\lfloor (\hat{s} \wedge s^*)n \rfloor -l+2}^{\lfloor (\hat{s} \vee s^*)n \rfloor}
	|\xi_i^{(k)}| = O_{\mathbb{P}}(l) $, which follows from $|\hat{s}-s^*| = O_{\mathbb{P}}(n^{-1}) $.
	Similarly, one can show $ Z_n^{(3)} = o_{\mathbb{P}} (1) $, $ Z_n^{(4)} = o_{\mathbb{P}} (1) $
	and therefore  \eqref{443} holds,
which implies (observing \eqref{439})  the assertion \eqref{bLemma3}.

In order to prove the remaining statement \eqref{442}, note that
\begin{align} \label{expectation1}
\sup_{t\in[0,1]} \big| \hat{\mu}_1(t) - \mu_1(t) \big| \leq
\frac{1}{\sqrt{n}} \frac{n}{\lfloor \hat{s}n \rfloor}  \big(  Q_n^{(1)}  +Q_n^{(2)}  \big),
\end{align}
where
\begin{align*}
Q_n^{(1)} &  =
	\sup_{t\in[0,1]} \Big| \frac{1}{\sqrt{n}}
	\sum_{j=1}^{\lfloor (\hat{s}\wedge s^*) n \rfloor}
	\big( X_{n,j}(t) - \mu_1(t) \big) \Big|, \\
Q_n^{(2)} & =
	\sup_{t\in[0,1]} \Big| \frac{1}{\sqrt{n}}
	\sum_{j=\lfloor (\hat{s}\wedge s^*) n \rfloor +1}^{\lfloor (\hat{s}\vee s^*) n \rfloor}
	\big( X_{n,j}(t) - \mu_1(t) \big) \Big| ~.
\end{align*}
Recall the definition of $\mathbb{\hat V}_n$ in \eqref{39}, then Theorem \ref{WIP}
and  the extended continuous mapping theorem (see Theorem 1.11.1 in  Van der Vaart and Wellner \cite{wellner1996})
yield
\begin{align} \label{expectation2}
Q_n^{(1)}&=
	\sup_{t\in[0,1]} \big|
	\hat{\mathbb{V}}_n \big(\lfloor (\hat{s}\wedge s^*) n \rfloor / n ,t\big) \big|
	= O_{\mathbb{P}}(1) ~.
\end{align}
On the other hand,
\begin{align*}
 \lim_{p\to\infty} \limsup_{n\to\infty} ~ \mathbb{P} \big(
	Q_n^{(2)} > p \big)  & \leq \lim_{p\to\infty} \limsup_{n\to\infty} ~ \mathbb{P} \big(
	Q_n^{(2)}  > p ~ , ~
	\big|\lfloor \hat{s}n \rfloor - \lfloor s^*n \rfloor \big| \leq p  \big)  \\
 & \qquad+ \lim_{p\to\infty} \limsup_{n\to\infty} ~ \mathbb{P} \big(
	\big|\lfloor \hat{s}n \rfloor - \lfloor s^*n \rfloor \big| > p  \big).
\end{align*}
The second term in this inequality is of order $o_{\mathbb{P}} (1)$,
since $|\hat{s}-s^*|=O_{\mathbb{P}}(n^{-1})$.
For the first term use Markov's inequality and Assumption (A1) to obtain
\begin{align*}
\lim_{p\to\infty} \limsup_{n\to\infty} ~ \mathbb{P} \big(
	Q_n^{(2)}   > p ~ , ~
	\big|\lfloor \hat{s}n \rfloor - \lfloor s^*n \rfloor \big| \leq p  \big)
&  \leq \lim_{p\to\infty} \limsup_{n\to\infty} ~ \mathbb{P} \Big(
	\frac{1}{\sqrt{n}}
	\sum_{j=\lfloor s^* n \rfloor - p}^{\lfloor s^*n \rfloor + p}
	\big( \|X_{n,j}\|+ \|\mu_1\| \big) > p \Big) \\
& \lesssim \lim_{p\to\infty} \limsup_{n\to\infty} ~
	\frac{1}{\sqrt{n}} = 0 ~.
\end{align*}
This means that $Q_n^{(2)} = O_{\mathbb{P}}(1) $.
Observing the calculations in \eqref{expectation1} and \eqref{expectation2} leads to
\begin{align} \label{449}
\sup_{t\in[0,1]} \big| \hat{\mu}_1(t) - \mu_1(t) \big| = O_{\mathbb{P}}\bigg(\frac{1}{\sqrt{n}}\bigg)
\qquad\mbox{and}\qquad
\sup_{t\in[0,1]} \big| \hat{\mu}_2(t) - \mu_2(t) \big| = O_{\mathbb{P}}\bigg(\frac{1}{\sqrt{n}}\bigg),
\end{align}
the second estimate following similarly. This  yields
 \eqref{442} and completes the proof of  \eqref{bLemma3}.
\end{proof}

\begin{proof}[\bf Proof of Theorem \ref{jointConvergence}]
In order to prove the assertion define
\begin{align*}
  \bold{D}_n = |\hat{d}_\infty - d_\infty| + \sup_{t\in[0,1]}|\hat{\mu}_1(t) - \hat{\mu}_2(t)
   - (\mu_1(t) - \mu_2(t))|.
\end{align*}
 Observing \eqref{449}  and  Corollary \ref{setT}, yields
\(\mathbb{P}({c_n}/{\sqrt{n}} \geq \bold{D}_n)\to 1\) as \(n\to\infty\), and
 the same  arguments as  given in the proof of Lemma \ref{setConvergence}
  show
%\begin{align} \label{CPsetConvergence}
$  d_h(\hat{\mathcal{E}}^\pm_n,\mathcal{E}^\pm ) \xrightarrow{\mathbb{P}^*} 0, $
% \end{align}
where the sets $\hat{\mathcal{E}}^+_n$ and $ \hat{\mathcal{E}}^-_n$ are now defined in \eqref{estimatedSetsCP}.
Theorem \ref{bTheorem} and   Lemma \ref{maxSetConvergence}
 lemma imply
\begin{align*}
(\tilde{D}_n(\mathcal{E}),~ T_n^{(1)},\dots,T_n^{(R)})
\Rightarrow (T(\mathcal{E)},~ T^{(1)},\dots,T^{(R)})
\end{align*}
in \(\mathbb{R}^{R+1}\), where
\begin{align*}
  \tilde{D}_n(\mathcal{E}) = \frac{1}{\hat{s}(1-\hat{s})}
    \max\big\{\max_{t\in\mathcal{E}^+} \hat{W}_n(\hat{s},t),
    \max_{t\in\mathcal{E}^-} (-\hat{W}_n(\hat{s},t)) \big\}
\end{align*}
and the random variable \(\hat{W}_n\) is defined by \eqref{313}. From the discussion following Theorem \ref{rate} and
the proofs of
Lemma \ref{lem:proofs:weak_conv:1} and Lemma \ref{lem:proofs:weak_conv:2}, $
  \tilde{D}_n(\mathcal{E}) = \mathbb{D}_n/[{\hat{s}(1-\hat{s})}] + o_{\mathbb{P}}(1)
    = \sqrt{n}(\hat d_\infty - d_\infty) + o_{\mathbb{P}}(1)
$, where \(\mathbb{D}_n\) is defined by \eqref{deps}. This implies the assertion.
\end{proof}

\end{document}